\documentclass[11pt]{article}
\usepackage[margin=1in]{geometry} 
\usepackage{amssymb,amsfonts,amsmath,bbm,mathrsfs,stmaryrd}
\usepackage{xcolor}

\usepackage[colorlinks,
             linkcolor=black!75!red,
             citecolor=blue,
             pdftitle={The fundamental progroupoid of an affine two-scheme},
             pdfauthor={Theo Johnson-Freyd and Alex Chirvasitu},
             pdfsubject={category theory},
             pdfkeywords={ct.category theory, rt.representation theory, ag.algebraic geometry},
             pdfproducer={pdfLaTeX},
             pdfpagemode=None,
             bookmarksopen=true
             bookmarksnumbered=true]{hyperref}

\usepackage{tikz}
\usetikzlibrary{arrows,decorations.pathreplacing,decorations.markings,shapes.geometric,through,fit,shapes.symbols,positioning}

\usepackage[amsmath,thmmarks,hyperref]{ntheorem}
\usepackage{cleveref}

\crefname{section}{Section}{Sections}
\crefformat{section}{#2Section~#1#3} 
\Crefformat{section}{#2Section~#1#3} 

\crefname{subsection}{\S}{\S\S}
\crefformat{subsection}{#2\S#1#3} 
\Crefformat{subsection}{#2\S#1#3} 

\theoremstyle{change}

\newtheorem{lemma}{Lemma}[subsection]
\newtheorem{proposition}[lemma]{Proposition}
\newtheorem{corollary}[lemma]{Corollary}
\newtheorem{theorem}[lemma]{Theorem}

\theoremstyle{nonumberplain}

\theoremstyle{change}
\theorembodyfont{\upshape}
\theoremsymbol{\ensuremath{\blacklozenge}}

\newtheorem{definition}[lemma]{Definition}
\newtheorem{example}[lemma]{Example}
\newtheorem{remark}[lemma]{Remark}

\crefname{definition}{definition}{definitions}
\crefformat{definition}{#2definition~#1#3} 
\Crefformat{definition}{#2Definition~#1#3} 

\crefname{lemma}{lemma}{lemmas}
\crefformat{lemma}{#2lemma~#1#3} 
\Crefformat{lemma}{#2Lemma~#1#3} 

\crefname{proposition}{proposition}{propositions}
\crefformat{proposition}{#2proposition~#1#3} 
\Crefformat{proposition}{#2Proposition~#1#3} 

\crefname{example}{example}{examples}
\crefformat{example}{#2example~#1#3} 
\Crefformat{example}{#2Example~#1#3} 

\crefname{remark}{remark}{remarks}
\crefformat{remark}{#2remark~#1#3} 
\Crefformat{remark}{#2Remark~#1#3} 

\crefname{corollary}{corollary}{corollaries}
\crefformat{corollary}{#2corollary~#1#3} 
\Crefformat{corollary}{#2Corollary~#1#3} 

\crefname{theorem}{theorem}{theorems}
\crefformat{theorem}{#2theorem~#1#3} 
\Crefformat{theorem}{#2Theorem~#1#3} 

\crefname{equation}{}{}
\crefformat{equation}{(#2#1#3)} 
\Crefformat{equation}{(#2#1#3)}

\theoremstyle{nonumberplain}
\theoremsymbol{\ensuremath{\blacksquare}}

\newtheorem{proof}{Proof}
\newtheorem{proof of iso}{Proof of \cref{prop.iso}}
\newtheorem{proof of Magid}{Proof of \cref{prop.Magid}}
\newtheorem{proof of equalizers bis}{Proof of \cref{prop.equalizers bis}}
\newtheorem{sketch}{Sketch of proof}
\newtheorem{defnS}{Definition of $S$}

\DeclareMathOperator{\id}{id}

\newtheorem{ST=1}{$ST\cong \id$}
\newtheorem{TS=1}{$TS\cong \id$}

\newcommand\Oo{{\mathcal O}}

\newcommand\CC{{\mathbb C}}

\newcommand\NN{{\mathbb N}}

\newcommand\ZZ{{\mathbb Z}}
\newcommand\Fp{{\mathbb F}_p}

\newcommand{\et}{{\rm\acute{e}t}}

\newcommand\mono{\hookrightarrow}

\newcommand\from{\leftarrow}
\newcommand\longto{\longrightarrow}
\newcommand{\longrightleftarrows}{\raisebox{3pt}{\ensuremath{\longrightarrow}}\hspace{-17.5pt}\raisebox{-1pt}{\ensuremath{\longleftarrow}}}

\newcommand\verylongarrow[1]{\!{\tikz[anchor=base,auto,baseline=(a.base)] \draw[->] node (a) {\phantom{A}} ++(1cm,0) node (b) {\phantom{A}\!\!\!} (a.mid) -- node{$\scriptstyle #1$} (b.mid);}}
\newcommand\veryverylongarrow[1]{\!{\tikz[anchor=base,auto,baseline=(a.base)] \draw[->] node (a) {\phantom{A}} ++(1.5cm,0) node (b) {\phantom{A}\!\!\!} (a.mid) -- node{$\scriptstyle #1$} (b.mid);}}

\DeclareMathOperator{\Hom}{Hom}

\DeclareMathOperator{\spec}{spec}
\DeclareMathOperator{\maxspec}{maxspec}
\DeclareMathOperator{\Spec}{2Spec}
\DeclareMathOperator{\QCoh}{\textsc{QCoh}}
\DeclareMathOperator{\sh}{\textsc{Sheaves}}
\DeclareMathOperator{\psh}{\textsc{Presheaves}}
\DeclareMathOperator{\Forget}{Forget}
\DeclareMathOperator{\shify}{Sheafify}
\DeclareMathOperator{\Lan}{Lan}
\DeclareMathOperator{\Functors}{Functors}
\DeclareMathOperator{\Compat}{Compat}
\DeclareMathOperator{\Coconts}{Cocont's}
\DeclareMathOperator{\End}{End}
\DeclareMathOperator{\TK}{TK}
\DeclareMathOperator{\Gal}{Gal}
\DeclareMathOperator{\PSpec}{PSpec}
\DeclareMathOperator{\Pic}{Pic}
\DeclareMathOperator{\GL}{GL}

\DeclareMathOperator*{\colim}{colim}

\DeclareMathOperator*{\revarinjlim}{{\tikz[baseline=(lim.base)] \draw[->] node[inner sep = 0pt] (lim) {$\lim$} (lim.south west) ++(1pt,-3pt) coordinate (a) (lim.south east) ++(-1pt,-3pt) coordinate (b) (a) -- (b) ;}}
\renewcommand\varinjlim\revarinjlim

\DeclareMathOperator*{\revarprojlim}{{\tikz[baseline=(lim.base)] \draw[<-] node[inner sep = 0pt] (lim) {$\lim$} (lim.south west) ++(1pt,-3pt) coordinate (a) (lim.south east) ++(-1pt,-3pt) coordinate (b) (a) -- (b) ;}}
\renewcommand\varprojlim\revarprojlim

\newcommand\Fun\Functors

\newcommand{\pt}{\{\textrm{pt.}\}}

\newcommand{\define}[1]{{\em #1}}

\newcommand{\cat}[1]{\textsc{#1}}

\newcommand{\HOM}{\underline{\Hom}}

\newcommand{\qedhere}{\mbox{}\hfill\ensuremath{\blacksquare}}

\newcommand{\op}{{\rm op}}

\newcommand{\sch}[1]{\Spec({^{#1}\cat{Set}})}

\renewcommand{\square}{\mathrel{\Box}}

\newcommand\B{{\rm B}}

\newcommand\Section[1]{\section{#1}\setcounter{lemma}{0}}

\title{The fundamental pro-groupoid of an affine 2-scheme}
\author{Alex Chirvasitu\footnote{\url{chirvasitua@math.berkeley.edu}}\; and Theo Johnson-Freyd\footnote{\url{theojf@math.berkeley.edu}}}

\begin{document}
\maketitle

\begin{abstract}
  A natural question in the theory of Tannakian categories is: What if you don't remember $\Forget$?  Working over an arbitrary commutative ring $R$, we prove that an answer to this question is given by the functor represented by the \'etale fundamental groupoid $\pi_1(\spec(R))$, i.e.\ the separable absolute Galois group of $R$ when it is a field.  This gives a new definition for \'etale $\pi_1(\spec(R))$ in terms of the category of $R$-modules rather than the category of \'etale covers.  More generally, we introduce a new notion of ``commutative 2-ring'' that includes both Grothendieck topoi and symmetric monoidal categories of modules, and define a notion of $\pi_1$ for the corresponding ``affine 2-schemes.''  These results help to simplify and clarify some of the peculiarities of the \'etale fundamental group.  For example, \'etale fundamental groups are not ``true'' groups but only profinite groups, and one cannot hope to recover more: the ``Tannakian'' functor represented by the \'etale fundamental group of a scheme preserves finite products but not all products.
\end{abstract}


\tableofcontents

\Section{Introduction}

We begin by motivating our paper with a natural question in the theory of Tannakian categories.  We then provide an overview of our general framework, and finish the introduction with an outline of the paper.

\subsection{Motivation: What happens to the Tannaka theorem when there are multiple ``forgetful'' functors around?} \label{section.introquestion}

	Let $G$ be a monoid in \cat{Set} (for example, a discrete group) and $k$ a field.  Denote by $^G\cat{Vect}^k$ the category of $G$-representations on $k$-vector spaces.  It is a symmetric monoidal category, and the forgetful functor $\Forget: {^G\cat{Vect}^k} \to \cat{Vect}^k$ that forgets the $G$ structure is a symmetric monoidal functor.  We denote by $\End_\otimes(\Forget)$ the monoid of \emph{monoidal} natural endomorphisms of this functor.  Then the following is well known, and follows almost immediately from the Yoneda lemma:

\begin{theorem}[Tannaka-Krein Reconstruction] \label{thm.TK}
  The homomorphism $G \to \End_\otimes(\Forget)$ given by sending $g\in G$ to the endomorphism ``act by $g$'' of $\Forget(X)$ is an isomorphism.
  Any group (or, even, any monoid) is recoverable from its monoidal category of representations along with the data of the forgetful functor. \qedhere
\end{theorem}
For reference, see for example \cite{MR654325,MR1106898,MR1173027}.
\begin{remark}
  Most discussions of Tannaka-Krein theory try to recover a group only from its finite-dimensional representations.  The Tannaka-Krein theorem would be false as stated if only the finite-dimensional representations were given; in that case, $\End_\otimes(\Forget)$ would be the profinite completion of $G$.  Throughout this paper, when we speak of ``modules of a ring'' or ``representations of a group'' we always mean to include all modules/representations, including the infinite-dimensional ones.
\end{remark}

	Given the Tannaka-Krein theorem, a natural question arises: What can be recovered from the symmetric monoidal category if the forgetful functor is not recorded?  When $k = \CC$, an answer is given by \cite{MR654325}:
\begin{theorem} \label{thm.delignemilneresult}
  Let $G$ be a discrete group.  Up to isomorphism, the forgetful functor is the unique symmetric monoidal $\CC$-linear cocontinuous functor ${^G\cat{Vect}^\CC} \to \cat{Vect}^\CC$. \qedhere
\end{theorem}
(Recall that a functor is \define{cocontinuous} if it preserves all small colimits.)

	In general, one should consider the following object, which we will generalize in subsequent sections:
\begin{definition} \label{defn.TKcat}
  Let $G$ be a small category, and $R$ a commutative ring, and denote by $^G\cat{Mod}^R$ (resp.\ $\cat{Mod}^R$) the category of $G$-representations on $R$-modules (resp.\ the category of $R$-modules).  Then $\cat{Mod}^R$ is symmetric monoidal on account of $R$ being commutative, and $^G\cat{Mod}^R$ inherits a symmetric monoidal structure of ``pointwise tensor product.''   The \define{Tannaka-Krein category} $\TK(G,R)$ has as its objects the symmetric monoidal cocontinuous $R$-linear functors ${^G\cat{Mod}^R} \to {\cat{Mod}^R}$, and as its morphisms the monoidal natural transformations.
\end{definition}

In general, the Yoneda lemma shows that, when $R$ is a field, the canonical functor $G \to \TK(G,R)$ (sending an object of $G$ to the evaluation functor at that object and a morphism of $G$ to its induced morphism) is full and faithful.  \Cref{thm.delignemilneresult} says that when $R = \CC$ and $G$ is a groupoid, the canonical functor is an equivalence of groupoids.  The following result is known to experts, but perhaps not as well known as it should be (we could not find a direct reference, and so we provide a proof in \cref{prop.pi1_of_fields}):
\begin{theorem}
  Let $R = k$ be a field, and let $\Gal(k)$ denote its absolute Galois group.  Think of $\Gal(k)$ as a groupoid with one object, with its profinite topology.  Given a groupoid $G$ with the discrete topology, let $\Hom(\Gal(k),G)$ denote the groupoid whose objects are continuous (in the topological sense) functors $\Gal(k) \to G$ and whose morphisms are natural transformations.  Then there is a natural equivalence of groupoids $\Hom(\Gal(k),G) \cong \TK(G,k)$.
\end{theorem}

One of our main results (\cref{prop.Magid}) is to extend this theorem to rings, where we interpret the profinite groupoid $\Gal(R)$ in the sense of \cite{Magid1974,Borceux2001}.  This result determines $\Gal(R)$ up to equivalence, by identifying the functor it represents.

It is also well known that the notion of ``Galois group of a ring'' agrees with the notion of ``\'etale fundamental group of an affine scheme'' (there are many accounts of Galois theory and \'etale fundamental groups; see for example \cite{MR559531}).  Our result, then, gives a new but equivalent definition and construction of the \'etale fundamental group: usually the \'etale fundamental group of a scheme is constructed from the topos of \'etale covers of the scheme, whereas we construct it from the abelian category of quasicoherent modules.

\subsection{Generalization: 2-rings and 2-affine algebraic geometry}

	To formulate our results, we introduce a new notion of \define{commutative 2-ring}: a symmetric closed monoidal presentable category.  Examples include: Grothendieck topoi; for any scheme $S$, the category of quasicoherent $\Oo(S)$-modules; natural constructions arising from algebraic stacks.  Morphisms of commutative 2-rings generalize the geometric morphisms of topoi.  One should then understand our results as pertaining to \define{affine 2-schemes}, and in particular to their ``\'etale fundamental groups.''

 Recall one of the basic motivations of affine algebraic geometry: a ``space'' is determined by its algebra of ``functions'' (real-valued, say); conversely, any ring should be thought of as the ``algebra of functions'' on some generalized ``space.''  One quickly discovers, however, that there are well-behaved generalized spaces (non-affine schemes, geometric stacks) that are not determined by their algebras of (global) functions.  These spaces are, however, determined by their categories of quasicoherent modules (\cite{MR0232821,Lurie:fk}), and one should think of $\QCoh(S)$ as ``the algebra of $\cat{AbGp}$-valued functions on $S$.'' Similarly, a Grothendieck topos should be thought of as ``the algebra of \cat{Set}-valued functions'' on some ``space.'' In this sense, all of (1-)algebraic geometry is ``2-affine.''
 Following the usual logic of algebraic geometry, we encourage the reader to think of commutative 2-rings as algebras of (\cat{Set}- or \cat{AbGp}-valued) functions.  The corresponding category $\cat{Af2Sch}$ of affine 2-schemes is rich enough to talk about all of ``1-mathematics''; for example, we show in \cref{prop.1into2} that the category \cat{AfSch} of (1-)affine schemes embeds fully faithfully into \cat{Af2Sch}, and we prove the similar result for the 2-category \cat{Gpd} of small groupoids in \cref{lemma.scheming_fully_faithful}. (Lurie's Tannaka-Krein theorem for geometric stacks  \cite{Lurie:fk}, which we recall in \cref{prop.lurieTK}, is closely related.)

	Let $G$ be a small groupoid and $A$ a commutative 2-ring, for example $A = {^R}\cat{Mod}$ for a commutative ring $R$.  Then $^GA = \Fun(G, A)$ is again a commutative 2-ring.  From a geometric point of view, the Tannaka-Krein category from \cref{defn.TKcat} is precisely the category of maps $\Spec(A) \to \Spec({^GA})$ of ``2-schemes over $\Spec(A)$.''  Since $\cat{Gpd} \mono \cat{Af2Sch}$, we can think of $G$ as itself an affine 2-scheme, and we show in \cref{lemma.main} that $\Spec({^GA}) = G \times \Spec(A)$ in $\cat{Af2Sch}$.  Then rather than working with ``2-schemes over $\Spec(A)$,'' we can equivalently describe the Tannaka-Krein category as the category of maps $\Spec(A) \to G$.  When the groupoid $G$ has a unique object, it should be thought of as $\B(\text{some group})$, whence it is not surprising that the maps $\Spec(A) \to G$ are precisely the ``$G$-torsors on $\Spec A$'' (\cref{prop.GSet -> A are torsors}).

	From this point of view, one expects that the $G$-torsors on $\Spec A$ are controlled by a ``fundamental group'' $\pi_1(\Spec A)$.  Indeed, this is our main result (\cref{thm.Gal}):
\begin{theorem}
	There is a functor, which we call $\pi_1(-)$, from ``good'' (\cref{defn.good2rings}) affine 2-schemes to pro-groupoids, such that, for a small groupoid $G$ and a good affine 2-scheme $X$, there is a natural equivalence
\[ \Hom_{\cat{Af2Sch}}(X, G) \cong \Hom_{\cat{ProGpd}}(\pi_1(X), G). \]
\end{theorem}

\begin{remark}
  There is a 0-categorical version of our story, which we will briefly sketch.  We encourage the reader to work out the details as a warm-up exercise.
  
  Let $S$ be a set and $R$ a ring; then there is an $R$-algebra $^SR$ of $R$-valued functions on $S$.  One can then ask about the set of $R$-algebra homomorphisms $^SR \to R$.    It is then expected that $\Hom_{\cat{Schemes}/R}(\spec (R), \spec(^SR))$ would be controlled by some (pro-)set $\pi_0(R)$.  This is true provided $S$ is not too large compared with $R$.
  
  The case when $R=k$ is a field is sufficient to illustrate the size concerns.  Then there is an embedding $S \mono \Hom_{k\cat{-Alg}}({^Sk}, k)$, which sends a point $s\in S$ to the homomorphism that evaluates any function in $^Sk$ at $s$.  The question parallel to \cref{thm.delignemilneresult} is whether this injection is necessarily an isomorphism of sets.  The answer surprised us: it is ``yes'' provided that either $S$ and $k$ are both finite or $\lvert S \rvert \leq \lvert k \rvert$, and ``no'' when $S$ is infinite and $\lvert S \rvert > \lvert k \rvert$, but constructing the extra $k$-points in $\spec(^Sk)$ requires the axiom of choice.  Notice that these extra points evaporate upon field-extending to sufficiently large fields.  (Rather, they remain as points in $\maxspec(^Sk)$, but the residue fields blow up.)
  
  These size issues are not directly relevant to the paper at hand, because we only consider small groupoids $G$, whereas our commutative 2-rings are all roughly the size of \cat{Set}.  But they do highlight that one must take size concerns into account if trying to apply the results in this paper to essentially-large groupoids.
\end{remark}

\subsection{Overview of the paper}

	In \cref{section.prelim} we recall the requisite background on presentable categories (following \cite{Lurie2009}, we drop the word ``locally'' from ``locally presentable'').  A \define{presentable category} is a locally small category that is a cocompletion of a small category, all of whose objects are ``little''; see \cref{defn.locpres}.  An equivalent definition (\cref{thm.locpres1,thm.locpres2}) is that a presentable category is any category of sheaves on a (colimit) sketch --- a \define{colimit sketch} is a small category with some distinguished diagrams that are trying to be colimit cocones, and a \define{sheaf} thereon is a \cat{Set}-valued contravariant functor that realizes each distinguished diagram as a limit in \cat{Set}.  In particular, Grothendieck topoi are all presentable, but so are abelian categories of quasicoherent modules on a scheme.

Presentable categories have many good properties, enough so that they form the objects of a 2-category with many formal similarities to the (1-)category of (1-)abelian groups.  In particular, any cocontinuous functor out of a presentable category is necessarily a left adjoint (we suggest the adjective \define{saft}, for ``special adjoint functor theorem,'' for any category with this property); and presentable categories are complete, cocomplete, well-powered, and well-copowered (\cref{lemma.adjoints}).  Moreover, the 2-category of presentable categories (and left adjoints as 1-morphisms) admits a natural tensor structure satisfying a universal property analogous to the one for abelian groups, and hom-categories between presentable categories are necessarily presentable (\cref{cor.2abgp}).  We denote by \cat{2AbGp} the symmetric monoidal 2-category of presentable categories, left adjoints, and natural transformations (see \cref{defn.2abgp,remark.2abgp} for equivalent descriptions).  Finally, we define a \define{commutative 2-ring} to be a symmetric monoidal object in \cat{2AbGp}; the 1-opposite to the 2-category of commutative 2-rings is our 2-category \cat{2AfSch} of \define{affine 2-schemes}.

	In \cref{section.pi1construction} we use the language of commutative 2-rings to rephrase our question from \cref{section.introquestion}.  Namely, for any small groupoid $G$ and any commutative 2-ring $A$, the category $^GA$ of functors $G \to A$ is a commutative 2-ring, and indeed an $A$-2-algebra; the generalization of the Tannaka-Krein category from \cref{defn.TKcat} is the category $\TK(G,A) = \Hom_{A\cat{-Com2Alg}}({^GA},A)$
 of all morphisms $^GA \to A$ of commutative $A$-algebras.  By \cref{lemma.main}, this is equivalent to the category $\Hom_{\cat{Com2Ring}}({^G\cat{Set}},A)$ of morphisms $^G\cat{Set} \to A$ of commutative 2-rings.  In \cref{section.torsor} we give an ``internal'' description of such morphisms as \define{right $G$-torsors in $A$} (\cref{prop.GSet -> A are torsors}).  In the rest of \cref{section.pi1construction}, we prove that, when $G$ is a small groupoid, the category $\TK(G,A) = \cat{Tors}(G,A)$ of right $G$-torsors in $A$ is an essentially-small groupoid (\cref{thm.iso,prop.smallness}).  Along the way (\cref{section.cartesian}), we see that for the purpose of understanding $G$-torsors in a commutative 2-ring $A$, we can safely replace $A$ by a \define{cartesian} ($\otimes = \times$) 2-ring $\cat{Cog}(A)$, which is the category of cocommutative coalgebras in $A$.

	In \cref{section.mainresults} we prove our main result: that when a commutative 2-ring $A$ satisfies a mild technical condition having to do with the behaviour of coproducts (\cref{defn.good2rings}; the condition is satisfied in the main examples of abelian categories and Grothendieck topoi, and is designed to rule out bad 2-rings like complete lattices), the functor $G \mapsto \cat{Tors}(G,A)$ is pro-representable.  This requires that we first review (\cref{section.2limits}) enough of the theory of 2-limits and pro-objects in 2-categories to make sense of the result.  We call the representing pro-object of the functor $\cat{Tors}(-,A)$ the \define{fundamental pro-groupoid of $A$}, and denote it by $\pi_1(A)$.

	Finally, in \cref{section.examplesandconnections} we compare our $\pi_1(A)$ to older notions of Galois and (\'etale) fundamental groups.  We find (\cref{prop.Magid}) that when $A$ is the category of modules of a commutative ring $R$, then the profinitization $\pi_1^f(A)$ of $\pi_1(A)$ is the absolute separable Galois group of $R$ in the sense of \cite{Magid1974}.  When the ring $R$ is a field, this profinitization is no loss of data, and $\pi_1(A) = \Gal(R)$ on the nose (\cref{prop.pi1_of_fields}); we do not know if this generalizes to the non-field case.  Finally, in \cref{section.comparisonArtin} we consider the case when the commutative 2-ring $A$ is the category of sheaves on a connected site $C$; then we find perfect agreement between our $\pi_1(A)$ and the $\pi_1^\et(C)$ of  \cite{Artin1986} (note that this $\pi_1^\et$ generically contains more data than the profinite $\pi_1^\et$ usually considered in algebraic geometry; in particular, when $C$ is a good topological space, $\pi_1^\et(C)$ is its honest fundamental group, not the profinitization thereof).

	We conclude with some discussion in \cref{section.furtherresearch} of generalizations and directions for further research.


\subsection{Acknowledgements}

	We would like to thank D.\ Ben-Zvi, A.\ Geraschenko, M.\ Olsson, N.\ Reshetikhin, N.\ Rozenblyum, V.\ Serganova, C.\ Teleman, and H. Williams for their many helpful comments and discussions.  
	M.\ Brandenburg caught a number of errors in an early version of this paper.
	We would particularly like to thank the anonymous referee for the detailed and engaging review.  The correct definition of torsors in a general commutative 2-ring (and in particular the definition of the morphisms $\tau_{x}$ defined above \cref{lemma.Galoiscond}) is due to the referee, and we have also included a few comments and examples the referee suggested.
	Some of this work was completed while the second author was a visitor at Northwestern University, whom he thanks for the hospitality.  This work is supported by the NSF grant DMS-0901431.

\Section{Preliminaries} \label{section.prelim}

	In this section we define our basic objects of study.  Our motivation is a desire to think of both Grothendieck topoi and monoidal abelian categories as types of ``commutative 2-rings'' with which we can set up a theory of ``affine 2-algebraic geometry.''  In particular, we need a notion of ``2-abelian group,'' and find that a convenient definition is that of a \define{presentable category}.  We will recall all definitions here, but leave many proofs to the references.

	By ``there are setly many X'' we mean that the collection of isomorphism classes of Xs, which might a priori be too big to be a set, is in fact a (small) set: for example, a category is \define{locally small} if there are only setly many morphisms between any pair of objects.  The notions of ``limit'' and ``colimit'' of a diagram are standard.  A category is \define{(co)complete} if every diagram with setly many arrows has a (co)limit.  A functor is \define{(co)continuous} if it preserves all existent (co)limits of diagrams with setly many arrows.

	We write ``$=$'' when we mean that there is a sufficiently canonical natural equivalence. Given categories $X,Y$, we write $\Fun(X,Y)$ or $^XY$ for the category of functors $X\to Y$, while the full subcategory on the cocontinuous functors is $\Coconts(X,Y)$. We also use the notation $Y^X$ for the category $\Fun(X^\op,Y)$, and $^XY^Z$ for $\Fun(X\times Z^\op,Y)$. We write \cat{CAT} for the (strict) 2-category of locally small categories.  Many of our 2-categories embed into \cat{CAT}, but when we say ``2-category'' (and in particular ``monoidal category'') we mean that the associators and so on may be nontrivial.  Nevertheless, we suppress all associators and similar objects from our notation.  See the start of \cref{section.2limits} for further discussion and references on higher categories.

	We typically identify a groupoid with its set of arrows, and write $G_0$ for its set of objects. In other words, we might write $g\in G$ for an arrow $g$ of the groupoid $G$. Despite this, $a\in A$ might also mean that $a$ is an object of the category $A$; this never applies to groupoids however, and the significance of the symbol $\in$ will be clear from the context. The source and target maps $G\to G_0$ of a groupoid are written $r_G$ and $\ell_G$ respectively, or, when there is no danger of confusion, simply $r$ and $\ell$. These stand for ``right'' and ``left,'' and the reader will find it useful to imagine arrows going from right to left (so that they compose properly, i.e.\ $fg$ is $g$ followed by $f$). In accordance with this convention, we write $G(t,s)$ for the set of arrows of a groupoid $G$ from $s$ to $t$.


\subsection{Presentable categories}

Most of the results in this section and the next are known, and indeed due to \cite{MR0327863} (see also \cite{MR0292908,MR1294136}).  In many places we outline the proofs for readers who, like ourselves, are not fully-versed in this material. We also refer the reader to \cite{MR1291599,MR1313497,MR1315049} for details.

\begin{definition}\label{defn.locpres}
  Fix a cardinal $\kappa$; a set is \define{$\kappa$-small} if its cardinality is strictly less than $\kappa$, and a category is $\kappa$-small if its collection of arrows is. A \define{$\kappa$-directed category} is a poset in which every $\kappa$-small set of objects has an upper bound. A \define{$\kappa$-directed colimit} is the colimit of a diagram for which the indexing category is $\kappa$-directed. An object $x$ of a category is \define{$\kappa$-little} if $\Hom(x,-)$ preserves $\kappa$-directed colimits.  An object is \define{little} if it is $\kappa$-little for some $\kappa$.
  
  A category $C$ is \define{presentable} if it is locally small and cocomplete and there exists a set $S$ of objects, all of which are little, such that the cocompletion in $C$ of (the full subcategory on) $S$ is the entire category $C$.  A category is \define{$\kappa$-presentable} if the generating set can be taken to consist entirely of $\kappa$-little objects.
\end{definition}
\begin{remark}
  Following \cite{Lurie2009}, we prefer to substitute ``presentable'' for the usual term ``locally presentable'' (\cite{MR1294136}). Also, what we call ``little'' is usually called ``small,'' ``presentable,'' or ``compact'' in the literature, but we prefer not to overburden a word if it has an underused synonym.
\end{remark}
\begin{remark}  
  Let $\lambda > \kappa$ be cardinals.  Then every $\lambda$-directed colimit is also $\kappa$-directed, and so every $\kappa$-little object is also $\lambda$-little.  Thus every $\kappa$-presentable category is also $\lambda$-presentable.
\end{remark}

\begin{remark}
  One can also define a \define{$\kappa$-filtered category} as one for which every $\kappa$-small subcategory can be completed to some cocone (not necessarily satisfying any universal property), and hence arrive at a corresponding notion of \define{$\kappa$-filtered colimit}.  
  But for every $\kappa$-filtered diagram there is a canonical $\kappa$-directed diagram with the same colimit, and so replacing the word ``directed'' by ``filtered'' throughout \cref{defn.locpres} does not change the notion of ``$\kappa$-little,'' and hence does not change the notion of ``presentable'' \cite[Theorem 1.5]{MR1294136}.
\end{remark}

\begin{example}\label{eg.locpres}
  (a) Most categories ``in nature'' are presentable.  Given any ring $R$, the category of $R$-modules is $\aleph_0$-presentable.  Varieties of algebras are presentable: the category of groups, for example, is $\aleph_0$-presentable.  In general, a variety of algebras is $\kappa$-presentable where $\kappa$ is larger than the arity of any operation.  For example, the category of Banach spaces and contractions is $\aleph_1$-presentable but not $\aleph_0$-presentable, because of the presence of countable sums.  
    
  (b) Every Grothendieck topos is presentable. 
    
  (c) The category of topological spaces is not presentable.  The opposite category to a presentable category is not presentable, unless it is (equivalent to) a partially ordered set.
\end{example}

We refer the reader to \cite{MR1294136} for the following results:
\begin{lemma}\label{lemma.presprop}
  (a) Let $C$ be $\kappa$-presentable.  Then there are setly many $\kappa$-little objects in $C$, and in particular one can take $S$ in \cref{defn.locpres} to consist of representatives of all  isomorphism classes of $\kappa$-little objects.
  
  (b) Any $\kappa$-small colimit of $\kappa$-little objects is $\kappa$-little.  Hence every object of a presentable category is little.
  
  (c) Every presentable category is both complete and cocomplete and both well-powered and well-copowered.
  
  (d) Let $C$ be a presentable category and $D$ any locally small category.  Then a functor $C \to D$ is a left adjoint iff it is cocontinuous.  If $D$ is also presentable, then a functor $C \to D$ is a right adjoint iff it is continuous and commutes with $\kappa$-directed colimits for some $\kappa$.  \label{lemma.adjoints}
\qedhere \end{lemma}
\begin{remark}
  A category $C$ for which every cocontinuous functor with domain $C$ is a left adjoint is called ``compact'' in \cite{MR850527}.  We suggest instead the word \define{saft}, short for ``special adjoint functor theorem'': on the one hand, ``compact'' is an overloaded word, and on the other hand, this property is the conclusion of any result that deserves the name ``special adjoint functor theorem'' \cite{MR0281769}.
  
  A saft category is necessarily complete.  Although most saft categories (and all ``total'' categories) are also cocomplete (in particular, cocompleteness is a common condition in special adjoint functor theorems), saftness does not imply cocompleteness \cite{MR0470020}.
\end{remark}

\begin{definition}\label{defn.2abgp}
  The 2-category \cat{2AbGp} of \define{2-abelian groups} has as objects the presentable categories.  A \define{1-morphism of 2-abelian groups} $f:A \to B$ is the data consisting of a left adjoint $f^*:A \to B$ and its right adjoint $f_*:B \to A$.  Given $f,g: A \to B$ two 1-morphisms, a \define{2-morphism in the category of 2-abelian groups} is a natural transformation $\eta: f^* \to g^*$; this uniquely determines an adjoint natural transformation $\eta^\dagger: g_* \to f_*$.
\end{definition}

\begin{remark}\label{remark.2abgp}
  One can present the 2-category \cat{2AbGp} up to equivalence in various ways.  By \cref{lemma.adjoints}, an equivalent category on the same objects has $\Hom(A,B) = \Coconts(A, B)$. 
  
  We write $\cat{2AbGp}^{1\,\op,\,2\,\op}$ for the 2-category that is opposite to $\cat{2AbGp}$ at both the 1- and 2-levels.  It is equivalent to the 2-category with objects the presentable categories and $\Hom(A,B) = $ the full subcategory of $\Fun(A, B)$ on the continuous functors that commute with $\kappa$-directed colimits for some $\kappa$.
\end{remark}

\begin{remark}
  Let us now outline some motivation for the \cref{defn.2abgp}.  Whatever a ``2-abelian group'' is, it should be a locally small category, and it is natural, in light of both finite sets and finite-dimensional vector spaces, to replace addition with the formation of colimits.  So a ``2-abelian group'' should certainly be cocomplete, and a morphism of 2-abelian groups should certainly be cocontinuous.   But remember that an addition-preserving function between abelian semigroups need not be particularly nice, whereas an addition-preserving function between abelian groups is a much richer thing.  Similarly, we argue that the right notion of ``2-abelian group'' includes the property of being saft.
  
  But in fact we ask for even more.  The categories of cocontinuous functors between saft categories are not necessarily locally small, and in particular not necessarily themselves saft.  Rather, we would like a theory in which the hom categories between 2-abelian groups are again 2-abelian groups.  For local smallness of categories of cocontinuous functors, it is enough to ask that the domain be a cocompletion of a small category.  Presentable categories satisfy all our desired conditions and include all our desired examples.
\end{remark}

Our choice of \cat{2AbGp} has many of the properties of the usual category of (1-)abelian groups.  For example, we have:
\begin{proposition}\label{prop.limitsin2abgp}
  The 2-category \cat{2AbGp} is complete and cocomplete. One can compute \linebreak\mbox{(2-)}limits and colimits in \cat{2AbGp} as follows:
  
  Let \cat{CAT} denote the 2-category of locally small categories.  The faithful embeddings of $\cat{2AbGp}$ and $\cat{2AbGp}^{1\,\op,\,2\,\op}$ into $\cat{CAT}$ from \cref{remark.2abgp} preserve limits.  So to compute a limit in \cat{2AbGp}, just compute the corresponding limit of cocontinuous functors; to compute a colimit, look ``in the mirror'' and compute the limit of right adjoints. 
\end{proposition}
\begin{sketch}
	The statement about limit calculation is a simple tweak on \cite[Exercise 2.n]{MR1294136} (which says essentially that a limit of accessible categories is again accessible), by adding the observation that a limit of a diagram in $\cat{CAT}$ consisting of cocontinuous functors between cocomplete categories is again cocomplete, and the canonical cone out of the limit consists of cocontinuous functors. 
	
	A similar proof works for the second assertion, concerning colimits, by using the characterization of right adjoints (\cref{lemma.presprop}) and the observation that the limit in $\cat{CAT}$ of a diagram consisting of such functors between cocomplete categories with $\kappa$-filtered limits is a category which again has all colimits and $\kappa$-filtered limits (and the canonical cone out of the limit consists of cocontinuous functors which preserve $\kappa$-filtered limits).   
\end{sketch}


\subsection{Sheaves on sketches}

We now develop a ``Morita'' theory equivalent to \cat{2AbGp}, extending the Morita theory of rings.  The basic object of study are colimit sketches.  The results in this section are not new: the main result (\cref{thm.locpres1,thm.locpres2}) dates back to \cite{MR0158921,MR0209335} and holds that presentable categories are precisely the categories of models of limit sketches.

\begin{definition}
  A \define{colimit sketch} consists of the following data: A category $X$, a collection $D$ of ``distinguished diagrams'' in $X$, and an assignment $c$ of each diagram in $D$ to a cocone over it, called a ``distinguished cocone,'' which should be thought of as a prescribed ``colimit'' over it.  In the usual way with algebraic theories, we will abuse notation, and refer to the sketch $(X,D,c)$ by the same name as we use for its underlying category $X$.  A sketch is \define{($\kappa$-)small} if its underlying category is.

  Let $X$ be a sketch and $C$ a category.  A functor $X\to C$ is \define{compatible} if   it takes distinguished cocones to colimit cocones --- if it realizes the prescriptions.  We write $\Compat(X,C)$ for the full subcategory of $\Functors(X, C)$ on  the compatible functors.  A \define{sheaf} on a sketch $X$ is a functor $X^\op \to \cat{Set}$ that takes the distinguished cocones in $X$ to limit cones in \cat{Set}.  
  We set $\sh(X) = \Compat(X,\cat{Set}^\op)^\op$, which is a full subcategory of $\psh(X) = \Fun(X^\op, \cat{Set})$.
  
  The \define{Yoneda embedding} is the functor $\gamma_X : X \to \psh(X)$ given by $x \mapsto \Hom(-,x)$.  It  factors through $\Forget: \sh(X) \to \psh(X)$ iff every distinguished cocone in $X$ is already a colimit cocone.
\end{definition}

\begin{example}
  Let $C$ be a locally small category.  Let $C_\kappa$ be the full subcategory of $C$ on the $\kappa$-little objects.   We make it into a sketch by declaring that the distinguished cocones are precisely the $\kappa$-directed colimits that exist in $C_\kappa$.
\end{example}

	The following is \cite[Theorem 1.46]{MR1294136} (an error in the proof was corrected in \cite{MR1699341}):
\begin{theorem}\label{thm.locpres1}
  A locally small category $C$ is $\kappa$-presentable if and only if the Yoneda embedding $C \to \sh(C_\kappa)$ is an equivalence.\qedhere
\end{theorem}

	The converse is:
\begin{theorem}\label{thm.locpres2}
  Let $X$ be a small sketch.  Then $\sh(X)$ is presentable.  The functor $\Forget: \sh(X) \to \psh(X)$ has a left adjoint $\shify : \psh(X) \to \sh(X)$, and the composition ${\shify} \circ \gamma : X \to \sh(X)$ is compatible.  Moreover, $\sh(X)$ is the universal cocomplete category receiving a compatible functor from $X$, in the following sense: if $C$ is any cocomplete category, then pulling back along ${\shify} \circ \gamma$ yields a natural-in-$C$ equivalence of categories:
  \[ \Compat(X, C) = \Coconts(\sh(X), C) \]
\end{theorem}

	We include the proof to illustrate the construction of 2-abelian groups.  In particular, in \cref{cor.2abgp,prop.boxtimes_G} and elsewhere, we will build presentable categories as categories of sheaves on sketches, and we want the reader to have some direct understanding of such constructions.

\begin{proof}
  It is well-known that $\psh(X)$ is presentable, and by definition $\sh(X)$ is a full subcategory.  A straightforward exercise, turning on the slogan ``limits commute,'' shows that $\sh(X)$ is closed under limits in $\psh(X)$; also straightforward is that $\sh(X)$ is closed under $\kappa$-directed colimits in $\psh(X)$ for some $\kappa$, because $X$ is small and $\kappa$-directed colimits commute with $\kappa$-small limits in \cat{Set}.  Hence by \cite{MR1189203} the functor $\Forget: \sh(X) \mono \psh(X)$ has a left adjoint $\shify$ satisfying $\shify \circ \Forget = \id_{\sh(X)}$.  Since $\psh(X)$ is the cocompletion of the image of $\gamma: X \to \psh(X)$ and $\shify$ is cocontinuous (being a left adjoint), $\sh(X)$ is the cocompletion of the image of ${\shify} \circ \gamma : X \to \sh(X)$, and hence presentable.
  
  Let $D$ be a distinguished diagram in $X$ and $D \to d$ its distinguished cocone, and pick $y\in \sh(X)$.  Then:
  \begin{align*}
    \{\text{cocones } \shify(\gamma(D)) \to y\} & = \{\text{cocones }\gamma(D) \to y\} && \text{(by the adjunction)} \\
     & = \{\text{cones }\pt \to y(D) \text{ in }\cat{Set} \} && \text{(by Yoneda)} \\
     & = y(d) && \text{(because $y$ is a sheaf)} \\
     & = \Hom_{\psh(X)}(\gamma(d),y) && \text{(Yoneda)} \\
     & = \Hom_{\sh(X)}(\shify(\gamma(d)),y). && \text{(by the adjunction)}
  \end{align*}
  This proves that ${\shify} \circ \gamma$ is compatible.
  
  Thus we have the ``pulling back along ${\shify} \circ \gamma$'' functor:
  \[ \Coconts(\sh(X), C) \veryverylongarrow{\circ ({\shify} \circ \gamma)} \Compat(X, C) \]
  We wish to show that it is an equivalence. 
  
  The equivalence is well-known when $X$ has no distinguished diagrams:
  \[ \gamma: \Coconts(\psh(X), C) \overset\sim\longleftrightarrow \Functors(X, C) : \Lan_\gamma\]
  Let $f : X \to C$ be any functor.  Then $f_! = \Lan_\gamma f: \psh(X) \to C$ is a left adjoint; its right adjoint is $f^* : c \mapsto \Hom_C(f(-),c) \in \psh(X)$.  Moreover, $f^*$ lands in $\sh(X)$ iff $f$ is compatible.
  
  Suppose that $g: \psh(X) \to C$ is any cocontinuous functor whose right adjoint $g^\dagger$ lands in $\sh(X)$:
  \[ \tikz[baseline=(C.base)] {
  \path node (C) {$C \,\overset{g^\dagger}{\raisebox{3pt}{\ensuremath\longto}}\, \sh(X) \overset{\Forget}{\underset{\shify}{\longrightleftarrows}} \psh(X) $}
  (C.south west) ++(15pt,10pt) coordinate (sw)
  (C.south east) ++(-50pt,10pt) coordinate (se)
  ;
  \draw[->] (se) .. controls +(-1,-.5) and +(1,-.5) .. node[auto] {$\scriptstyle g$} (sw)
  ; } \]
  The upper arrows are right adjoints and the lower arrows are left adjoints.
  Then, since $\Forget$ is full, the left adjoint of $g^\dagger: C \to \sh(X)$ is necessarily $g \circ {\Forget}$.  On the other hand, the composition of adjoints is the adjoint to the composition.  So we have that $g = (g \circ {\Forget}) \circ {\shify}$ is a factorization into left adjoints.
  
  On the other hand, suppose that $g$ factors as $g = g' \circ \shify$, with $g'$ cocontinuous.  Then its adjoint $g^\dagger$ factors as $g^\dagger = \Forget \circ (g')^\dagger$, where $(g')^\dagger$ is the right adjoint to $g'$.  With the previous paragraph, we have shown that for a cocontinuous functor $g : \psh(X) \to C$, the following are equivalent: (i) $g$ is a composition of $\shify$ and a cocontinuous functor $\sh(X) \to C$; (ii) the right adjoint to $g$ factors through $\Forget$; (iii) $g \circ {\Forget} : \sh(X) \to C$ is cocontinuous and $g = g \circ {\Forget} \circ {\shify}$.
  
  All together, we have adjoint functors:
  \begin{gather*}
    p = {\Lan_\gamma} \circ {\Forget} : \Compat(X, C) \longto \Coconts(\sh(X), C) \\
    q = \mbox{} \circ {\shify} \circ \gamma : \Coconts(\sh(X), C) \longto \Compat(X, C)
  \end{gather*}
  We study their compositions:
  \begin{equation*}
    qp(f)  = (\Lan_\gamma f) \circ {\Forget} \circ {\shify} \circ \gamma 
      = (\Lan_\gamma f) \circ \gamma 
      = f
  \end{equation*}
  The second isomorphism is from the equivalence above, with $g = f_! = \Lan_\gamma f$ for $f$ compatible; the last is by the equivalence when $X$ has no distinguished diagrams.
  \begin{equation*}
    pq(g)  = \Lan_\gamma(g \circ {\shify} \circ \gamma) \circ {\Forget} 
      = g \circ {\shify} \circ {\Forget} 
      = g
  \end{equation*}
  The second isomorphism is because $ g \circ {\shify} : \psh(X) \to C$ is cocontinuous, and we again harness the statement when $X$ has no distinguished diagrams.
\end{proof}

\begin{corollary}\label{cor.2abgp}
  Let $A,B$ be presentable categories.  Then $\Hom_{2\cat{AbGp}}(A, B)$ is presentable.
  Moreover, the 2-category $2\cat{AbGp}$  is symmetric monoidal closed: for every pair of presentable categories $A,B$, there is a presentable category $A \boxtimes B$ so that for every cocomplete $C$ there is a natural equivalence of categories:
  \begin{multline*} \Hom_{2\cat{AbGp}}(A\boxtimes B,C) = \Hom_{2\cat{AbGp}}(A,\Hom_{2\cat{AbGp}}(B,C)) = \\ = \{\text{functors }A \times B \to C \text{ that are cocontinuous in each variable} \\ \text {while holding the other variable fixed}\}\end{multline*}
\end{corollary}

\begin{proof}
  Let $A',B'$ be any small sketches so that $A = \sh(A')$ and $B = \sh(B')$; such  sketches exist by \cref{thm.locpres1}.
    By \cref{thm.locpres2}, $\Hom_{2\cat{AbGp}}(A,B) = \Compat(A',B)$ is a full subcategory of $\Fun(A',B) = \sh((A'^{\op})\times B',\cat{Set})$, where $(A')^{\op}\times B'$ is given a sketch structure in which the distinguished cocones are those of the form (object in $A'$)$\times$(distinguished cocone in $B'$).  In particular, $\Fun(A',B)$ is presentable.  Let $D$ be the set of distinguished cocones in $A'$.  Then $\prod_{d\in D} B$ is presentable by \cref{prop.limitsin2abgp}.  For each diagram $d\in D$, let $c(d)$ denote its tip.  Then $d$ determines two functors $\Fun(A',B) \to B$.  One of them sends $f\in \Fun(A',B)$ to $f(c(d))$.  The other sends $f$ to $\colim f(d)$.  By the universal property of $\colim$, there is a uniquely defined morphism in $B$ from $\colim f(d)$ to $f(c(d))$.  Packaging these maps together for all $d\in D$, we can build a diagram in \cat{CAT} of the form:
    \begin{equation}\label{eqn.homABispresentable}
     \tikz[anchor=base,baseline=(G.base)]{
    	\path (0,0) node (G) {$\Fun(A',B)$} +(2.25,0) node {$\Uparrow$} +(4,0) node (H) {$\prod_{d\in D} B$};
    	\draw[->] (G) .. controls +(1,1) and +(-1,1) .. (H) node[pos=.5,auto] {$\scriptstyle f \mapsto f(c(d))$};
    	\draw[->] (G) .. controls +(1,-1) and +(-1,-1) .. (H) node[pos=.5,auto,swap] {$\scriptstyle f \mapsto \colim f(d)$};
  } \end{equation}
  Unpacking definitions shows that $\Compat(A',B)$ is precisely the limit in \cat{CAT} of \cref{eqn.homABispresentable}.  But colimits are computed pointwise, and colimits commute.  Thus both the top and bottom arrows in \cref{eqn.homABispresentable} are cocontinuous functors between presentable categories.  It follows from \cref{prop.limitsin2abgp} that $\Compat(A',B) = \Hom_{\cat{2AbGp}}(A,B)$ is the limit of \cref{eqn.homABispresentable} as a diagram in \cat{2AbGp}, and in particular $\Hom_{\cat{2AbGp}}(A,B)$ is presentable.

	If $C$ is a cocomplete category, then using \cref{thm.locpres2} we have a natural equivalence of categories:
  \begin{multline*}
    \Hom_{\cat{2AbGp}}(A,\Hom_{2\cat{AbGp}}(B,C)) = \Compat(A', \Compat(B',C)) = \\
     = \{\text{functors }A' \times B' \to C \text{ that are compatible in each variable} \\ \text {while holding the other variable fixed}\}
  \end{multline*}
  We put the structure of a sketch on the product $A' \times B'$ of underlying categories, by declaring that the distinguished cocones are those cocones of the form $(\text{distinguished cocone})\times(\text{constant})$ or $(\text{constant})\times(\text{distinguished cocone})$.  Then the right-hand side is precisely $\Compat(A'\times B', C)$, and we define $A \boxtimes B = \sh(A' \times B')$. \Cref{thm.locpres2} takes care of the rest.
\end{proof}

	Keeping the notations from the previous proof, note that there is a canonical separately-cocontinuous functor $A\times B\to A\boxtimes B$. We denote by $a\boxtimes b$ the image of the pair of objects $a\in A$, $b\in B$ through this functor. 

\begin{remark}\label{remark.Morita}
  One should think of a sketch $X$ as a collection of ``generators and relations'' for the 2-abelian group $\sh(X)$.  From this point of view, it is not surprising that non-isomorphic sketches can have equivalent categories of sheaves.  We can then present the 2-category $\cat{2AbGp}$ in a ``Morita'' or ``matrix'' style.  The objects are (colimit) sketches, and the morphisms $X \to Y$ are certain functors $X \times Y^\op \to \cat{Set}$.  One should think of such a functor as an $X,Y$-bimodule.  The 1-composition is a form of ``tensor product of bimodules,'' and can be presented using the language of coends \cite{MR651714}.
\end{remark}

\begin{example} \label{eg.GH_boxtimes}
  Let $G$ be a small category (we will be particularly interested in the case when $G$ is a groupoid).  Then $^G\cat{Set} = \Fun(G, \cat{Set}) = \sh(G^\op)$ is presentable, where $G^\op$ is the opposite category thought of as a sketch with no distinguished diagrams.  If $G,H$ are both small categories, then ${^G\cat{Set}} \boxtimes {^H\cat{Set}} = {^{G\times H}\cat{Set}}$, and $\Hom( {^G\cat{Set}}, {^H\cat{Set}} ) = {^{G^\op \times H}\cat{Set}}$.
  
  Similar results hold when $G,H$ are rings and $\cat{Set}$ is replaced by the category $\cat{AbGp}$ of abelian groups, although the computations are less obvious: writing $^G\cat{Mod}$ for the category of $G$-modules, we have ${^G\cat{Mod}} \boxtimes {^H\cat{Mod}} = {^{G\otimes H}\cat{Mod}}$, and $\Hom( {^G\cat{Mod}}, {^H\cat{Mod}} ) = {^{G^\op \otimes H}\cat{Mod}}$.
\end{example}

\begin{remark}
  In \cite{MR1106898}, Deligne defines a tensor product of abelian categories $A,B$, by asking that, for each abelian category $C$, the category of right exact functors $A\boxtimes B \to C$ be equivalent to the category of functors $A \times B \to C$ that are right exact in each variable. Between abelian categories, the right-exact functors are precisely the functors that preserve finite colimits; as such, Deligne's tensor product is a finitary (abelian) version of ours.  Since Deligne's primary applications are to categories that consist entirely of dualizable objects, he cannot demand cocontinuity in general --- his categories do not contain most infinite colimits.
  One can give a small abelian category the structure of a sketch by declaring all finite colimits to be distinguished, and then complete it to a presentable category; by the universal properties of the constructions, it is clear this operation of completion intertwines Deligne's tensor product with ours.
\end{remark}

We make the following observation, which will come in handy later:

\begin{lemma}\label{lemma.^GM}
    Let $M$ be a 2-abelian group, and $G$ a small category. We then have an equivalence $M\boxtimes{^G\cat{Set}}= {^GM}$.  
\end{lemma}
\begin{proof}
    We show that there is an equivalence 
\[
    \Hom_{\cat{2AbGp}}(^GM,N) \to \Hom_{\cat{2AbGp}}(M\boxtimes{^G\cat{Set}},N), 
\] 
natural in $N\in\cat{2AbGp}$. For categories $X$ and $Y$, denote by $L(X, Y)$ and $R(X, Y)$ the full subcategories of $\Fun(X, Y)$ consisting of left, and, respectively, right adjoints. Notice that we have an equivalence $L(X, Y)=R(Y, X)^\op$. We claim that there is a series of equivalences
\begin{multline}\label{mline.^GM}
    L({^GM},N) = R(N,{^GM})^\op = \bigl({^GR}(N,M)\bigr)^\op = L(M,N)^G = \\ = L(M,N^G) = L(M,L( {^G\cat{Set}}, N)) = L(M\boxtimes{^G\cat{Set}},N),
\end{multline}
all natural in $N$. We justify these equivalences as follows:

The first one is a consequence of the observation made in the sentence preceding \cref{mline.^GM}. 

The second, third, and fourth say that the superscript $^G$ (either left or right) commutes with $L(M,-)$ and $R(M,-)$. Recall that a functor between presentable categories is a left adjoint if and only if it preserves colimits (\cref{lemma.adjoints}). Since colimits in both $^GN$ and $\Fun(M,N)$ are computed pointwise, we conclude that a functor $M\to {^GN}$ is a left adjoint if and only if the corresponding functor $M\times G\to N$ preserves colimits in the first variable, if and only if the corresponding functor $G\to \Fun(M,N)$ takes values in $L(M,N)$. The argument for $R$ is similar, using the fact that a functor between presentable categories is a right adjoint if and only if it is continuous and it preserves $\kappa$-filtered colimits for some cardinal $\kappa$. 

The fifth equivalence is $N^G=\Hom_{\cat{2AbGp}}(^G\cat{Set},N)$, an immediate consequence of \cref{thm.locpres2} (simply use the fact that $^G\cat{Set} = \sh(G^\op)$, where $G^\op$ is a sketch with no distinguished diagrams). The sixth is the universal property of $\boxtimes$.   
\end{proof}
\begin{remark}
    By the universal property of the tensor product of 2-abelian groups, there is a canonical morphism $M\boxtimes {^G\cat{Set}}\to {^GM}$ (determined up to unique natural isomorphism). It is not hard to check that this canonical morphism is the equivalence exhibited in the proof of \cref{lemma.^GM}.   
\end{remark}


\subsection{2-rings and affine 2-schemes}

\begin{definition}
  A \define{2-ring} is an algebra object in \cat{2AbGp}.  Equivalently, a 2-ring is a presentable monoidal category for which the monoidal structure distributes over colimits.  2-rings are the objects of a 2-category \cat{2Ring}: 1-morphisms in \cat{2Ring} are 1-morphisms of the underlying 2-abelian groups equipped with the structure of monoidal functors, and 2-morphisms in \cat{2Ring} are monoidal natural transformations between 1-morphisms.
  
  A \define{commutative 2-ring} is a 2-ring equipped with a symmetric structure.  Commutative 2-rings form a 2-category \cat{Com2Ring} whose 1-morphisms are those 1-morphisms of the underlying 2-rings that intertwine the symmetric structures.  The 2-morphisms are left unchanged, so that for any two commutative 2-rings $A,B$, the category $\Hom_{\cat{Com2Ring}}(A,B)$ is a full subcategory of $\Hom_{\cat{2Ring}}(A,B)$.
  
  The category \cat{2AfSch} of \define{affine 2-schemes} is 1-opposite to the category of commutative 2-rings.  Given a commutative 2-ring $A$, we write $\Spec(A)$ for the corresponding affine 2-scheme. 
\end{definition}

\begin{remark}
  Let $A$ be a commutative 2-ring, with multiplication $\otimes: A\boxtimes A \to A$.  Then for each $a\in A$, the functor $\otimes a: A\to A$ is cocontinuous, and so has a right adjoint, $\HOM_A(a,-)$, and this construction is natural in $a$.  The associativity of $\otimes$ implies that $A$ is enriched over itself in the sense of \cite{MR651714}.
\end{remark}

\begin{definition}\label{defn.modalg}
  Let $A$ be a 2-ring.  A \define{left $A$-2-module} is an $A$-module object in \cat{2AbGp}; i.e.\ it is a 2-abelian group $M$, a map $A\boxtimes M \to M$, and an associator intertwining the two maps $A\boxtimes A \boxtimes M \to M$, satisfying some compatibility conditions.
  
  When $A$ is a commutative 2-ring, we denote by $A\cat{-2Mod}$ its 2-category of 2-modules.  It has a symmetric monoidal structure $\boxtimes_A$, similar to the analogous structure on the module category of a commutative 1-algebra.  A \define{(commutative) $A$-2-algebra} is a (commutative) monoid object in $A\cat{-2Mod}$; we write $A\cat{-Com2Alg}$ for the category of commutative $A$-2-algebras.
  
  A commutative $A$-2-algebra structure on a commutative 2-ring $B$ is the same as a commutative 2-ring morphism $A\to B$. Consequently, the 2-category $A\cat{-Com2Alg}$ can also be described as follows:  the objects are the commutative 2-ring morphisms $A\to B$; the morphisms between $A\overset f{\rightarrow}B$ and $A\overset g {\rightarrow}C$ are the pairs $(h,\zeta)$, where $h: B \to C$ is a 1-morphism of in \cat{Com2Ring} and $\zeta: g \Rightarrow h\circ f$ is a 2-isomorphism in \cat{Com2Ring}; a 2-morphism $(h,\zeta) \Rightarrow (h',\zeta')$ is a 2-morphism $\xi: h\Rightarrow h'$ such that $\zeta' = (\xi. f)\zeta$, where we denote  ``whiskering of a 2-morphism by a 1-morphism'' by a lower dot, and the ``2-composition'' by concatenation.
\end{definition}

\begin{example}\label{eg.2rings}
  The category $\cat{Set}$ is a commutative 2-ring with $\otimes = \times$.  It is the initial object in both $\cat{2Ring}$ and $\cat{Com2Ring}$; the unique 1-morphism sends a set with cardinality $\kappa$ to the $\kappa$-fold coproduct of the unit object.  The terminal 2-ring is \define{the zero ring}: the category with one object and only the identity morphism.

  More generally, let $G$ be any small category.  Then the 2-abelian group $^G\cat{Set} = \sh(G^\op)$ has the structure of a commutative 2-ring by setting $\otimes$ to be the categorical product --- indeed, $\sh(G^\op)$ is a Grothendieck topos, and any Grothendieck topos is a commutative 2-ring with $\otimes = \times$.  
\end{example}

\begin{remark}
	Regarding the previous example, there is a subtle point to keep in mind: the usual 2-category $\cat{Topos}$ consisting of (for us, Grothendieck) topoi, with left exact left adjoints as 1-arrows and natural transformations between those as 2-arrows (the 1-opposite of \cite[Definition VII.1.1]{MacLane1994}) embeds in $\cat{2Ring}$ (see the previous example), but not fully: in general, our morphisms of 2-rings, which in this case are just product-preserving left adjoints between topoi, need not be left exact (i.e. preserve all finite limits). The embedding is full, however, for topoi of the form $^G\cat{Set}$ for {\it groupoids} $G$. All in all we have 2-functors
	\[
		\cat{Cat}^{1\,\op}\to \cat{Topos} \hookrightarrow \cat{2Ring}, \qquad C\mapsto {^C}\cat{Set},
	\]
and the composition is fully faithful when restricted to $\cat{Gpd}^{1\,\op}$ (cf. \cref{lemma.scheming_fully_faithful}), but not in general; see \cref{section.nongroupoids} for an example. 
\end{remark}

\begin{example}\label{eg.Set[X]}
	Consider the monoid $\NN$ as a category with countably many objects and only identity morphisms. We can give the presentable category of sheaves on this countable discrete category a non-cartesian commutative 2-ring structure by $\otimes = $ convolution, or equivalently by extending $\otimes$ by cocontinuity from:
  \[ \hom(-,a) \otimes \hom(-,b) = \hom_{\NN}(-,a+b) \]
We denote this 2-ring by $\cat{Set}[X]$, as it is the 2-ring freely generated by one object $X=\hom(-,1)$.
\end{example}

\begin{example}
   The category $\cat{AbGp} = {^\ZZ\cat{Mod}}$ of (1-)abelian groups is a 2-ring with $\otimes = $ the usual tensor product of abelian groups.  If $R$ is a ring, then a 2-ring structure on the category $^R\cat{Mod}$ of left $R$-modules is the same as a \define{sesquialgebra} structure on $R$ in the sense of \cite{MR2304628}.  Every Hopf algebra is a sesquialgebra, and the corresponding 2-ring is commutative exactly when the Hopf algebra is triangular.  But more importantly for this paper, if $R$ is commutative as a ring, then $^R\cat{Mod}$ has a canonical commutative 2-ring structure given by $\otimes =  \otimes_{R}$.  Moreover, $^R\cat{Mod}$ is an $\cat{AbGp}$-2-algebra.
\end{example}

\begin{remark}\label{remark.boxtimes_is_coprod}
    We will not prove this in any detail, but there is a higher analogue of the fact that the usual tensor product of abelian groups is the binary coproduct in the category of commutative rings: $\boxtimes$ is the binary (2-)coproduct in $\cat{Com2Ring}$. 
    
    The idea of the proof is the usual one: first, there is a commutative 2-ring structure on $A\boxtimes B$ whenever $A,B$ are commutative 2-rings. Next, there is a 2-ring map $A\to A\boxtimes B$ obtained by sending $a\in A$ to $a\boxtimes 1_B$, and a similar map $B\to A\boxtimes B$. Finally, given maps $f:A\to C$ and $g:B\to C$ of commutative 2-rings one forms $f\otimes g:A\boxtimes B\to C$ by sending $a\boxtimes b$ to $f(a)\otimes g(b)$. We leave the task of constructing the monoidal structure on $f\otimes g$ to the reader, with the observation that just as in the case of ordinary commutative rings, the symmetric structure is crucial.  
\end{remark}

\begin{proposition}\label{prop.1into2}
  The assignment $R \mapsto \Spec({^R\cat{Mod}})$ of commutative rings to affine 2-schemes makes the category of (1-)affine schemes into a full sub-2-category of the category of affine 2-schemes.
\end{proposition}

\begin{proof}
  Let $R,S$ be two commutative rings.  The category $\Hom_{2\cat{AbGp}}({^R\cat{Mod}}, {^S\cat{Mod}})$ is naturally equivalent to the category ${^S \cat{Mod}^R}$ of left-$S$ right-$R$ bimodules (any bimodule $^SM^R$ gives a functor $M \underset R \otimes$, and all cocontinuous functors are of this form).  Equipping the functor $ M \underset R \otimes $ with a symmetric monoidal structure is the same as equipping the module $^S M ^R$ with:
  \begin{itemize}
    \item An isomorphism ${^S M} \cong {^S S}$, where $^S M$ is the left-$S$ module $M$ that has forgotten its $R$-module structure, and $^S S$ is the free rank-one left-$S$ module.  This isomorphism expresses that the functor sends the unit object $^RR \in {^R\cat{Mod}}$ to the unit object $^SS \in {^S\cat{Mod}}$.
    \item An associative, commutative isomorphism of left-$S$ right-$(R\otimes R)$ bimodules
    \[ {^S M ^R} \underset R \otimes {^R R ^{R\otimes R}} \cong {^S S ^{S\otimes S}} \underset {S\otimes S} \otimes ({^S M ^R} \otimes {^S M ^R}) \]
    where all the actions in ${^RR^{R\otimes R}}$, ${^SS^{S\otimes S}}$ are by multiplication (these bimodules exist because $R$ and $S$ are commutative).  This isomorphism expresses that the tensor product of modules maps to the tensor product of their images, and the associative and commutative properties express that the functor is monoidal and symmetrically so.
  \end{itemize}
  The isomorphism ${^S M} \cong {^S S}$ makes $S$ into a right-$R$ module compatible with the $S$ action; such a structure is the same as a homomorphism $R \to S$ of rings (an element of $R$ maps to its action on $1\in S$).  Upon identifying $M \cong S$, both sides of the second isomorphism are copies of $S$ on which $S$ acts on the left by multiplication, and $R\otimes R$ acts on the right via the homomorphism.  Any automorphism of this $^S S ^{R\otimes R}$ is of the form ``multiply by an invertible element of $S$.''  So the category of symmetric monoidal functors ${^R\cat{Mod}} \to {^S\cat{Mod}}$ is equivalent to a category whose objects are pairs: a homomorphism $R\to S$ and an invertible element of $S$.
  
  But a natural transformation between the functors $M \otimes$ and $N\otimes$, where $^SM^R$ and $^SN^R$ are bimodules, is the same as a homomorphism of bimodules.  In particular, any natural transformation between monoidal functors is of the form ``multiply by an element of $S$''; but such a natural transformation has a hard time being monoidal.  Let $\phi,\varphi : R\to S$ be homomorphisms and $\sigma,\varsigma \in S$ invertible elements, and let $(\phi,\sigma)$ and $(\varphi,\varsigma)$ be the corresponding monoidal functors.  Then there is no monoidal natural transformation $(\phi,\sigma) \to (\varphi,\varsigma)$ unless $\phi = \varphi$, and if we do have equality, then there is a unique such natural transformation, given by multiplication by $\varsigma\sigma^{-1}$.
\end{proof}

\begin{remark}\label{prop.lurieTK}
  A closely related and dramatically more powerful result is in \cite{Lurie:fk}.  An algebraic stack $X$ is \define{geometric} if it is quasicompact and the diagonal map $X \to X \times X$ is representable and affine.  There is a functor from the 2-category geometric stacks to \cat{2AfSch} that assigns to every stack its commutative 2-ring of quasicoherent sheaves.  This functor is faithful, but it is not known to be full --- morphisms between geometric stacks correspond precisely to those morphisms of affine 2-schemes that satisfy a technical condition called ``tameness'' --- but it is expected to be.  In personal correspondence \cite{Brandenburg}, Brandenburg has informed us of a result (joint with Dolan) that will appear in his thesis: all affine 2-scheme morphisms $\Spec(\QCoh(X)) \to \Spec(\QCoh(Y))$ are tame when $X$ and $Y$ are projective schemes, and Brandenburg has a rough outline of how to get to all geometric stacks.
\end{remark}

\begin{example}\label{eg.main}
  Generalizing \cref{eg.2rings} is our main object of study: Let $A$ be any commutative 2-ring (for example, the category of modules of a commutative ring) and $G$ any small category.  Then $^GA = \Fun(G, A) = \Coconts(\sh(G^\op), A)$ is a 2-ring, where the monoidal structure is given by ``pointwise tensor product'': if $\alpha,\beta : G \to A$ are functors, then for each object or arrow $g \in G$, we set $(\alpha \otimes \beta)(g) = \alpha(g) \otimes \beta(g)$.  Since colimits of functors are computed pointwise, this is a 2-ring structure, and the commutativity constraint is inherited from $A$.
\end{example}

The following lemma is of central importance:
\begin{lemma}\label{lemma.main}
  Let $G$ be a small category and $A$ a commutative 2-ring.  Then $^GA$ is an $A$-2-algebra.  Moreover, there is an equivalence of categories:
  \[ \TK(G,A) = \Hom_{A\cat{-Com2Alg}}({^GA}, A) \simeq \Hom_{\cat{Com2Ring}}({^G\cat{Set}}, A) \]
\end{lemma}
\begin{proof}
  The commutative 2-ring structure comes from that on $A$, as in \cref{eg.main}. 
    
  We know from \cref{lemma.^GM} that $^GA$ can be identified with $A\boxtimes{^G\cat{Set}}$ as 2-abelian groups, and from now on we make this identification. The $A$-2-algebra structure on $A\boxtimes{^G\cat{Set}}$ is given by the canonical ``inclusion'' of $A$ into the coproduct $A\boxtimes{^G\cat{Set}}$. Using this observation, the description of the 2-category $A\cat{-Com2Alg}$ given after \cref{defn.modalg}, and \cref{remark.boxtimes_is_coprod}, we have
  \begin{align*}
      \Hom_{A\cat{-ComAlg}}(A\boxtimes{^G\cat{Set}},A) &\simeq \Hom_{A\cat{-ComAlg}}(A,A)\times\Hom_{\cat{Com2Ring}}({^G\cat{Set}},A)\\ &\simeq                                      \Hom_{\cat{Com2Ring}}({^G\cat{Set}},A). 
  \end{align*}
\end{proof}


\subsection{Coends and tensor products of functors}\label{subsection.coends}

The formalism of tensor products $x\boxtimes_D y$ of two functors $x$ and $y$ defined on $D$ (where $x$ is contravariant while $y$ is covariant) will be very useful. The definition requires the notion of coend, as in \cite[IX]{MR1712872}, for example (especially section 6 of that chapter). This subsection will be devoted to introducing our notation and conventions, as well as exploring this construction through some preliminary results to be used below. 

Let $C$, $D$, $E$ be small categories, $M$ an arbitrary category, and $A$ a category which ``acts on $M$ on the right'' by means of a functor $\triangleleft:M\times A\to M$; we typically write $m\triangleleft a$ for $\triangleleft(m,a)$. Then, given $x\in{^C}M^D$, $y\in{^D}A^E$ and objects $c\in C$, $e\in E$, define $(x\boxtimes_D y)(c,e)$ to be the coend 
\[
	\int^{d\in D} x(c,d)\triangleleft y(d,e),
\]
assuming all such colimits exist (the categories on which we perform such constructions will always be cocomplete, and moreover, the actions will always be appropriately cocontinuous). It is clear that this definition is functorial in $c$ and $e$, and what we have defined is, in fact, a functor $\boxtimes_D:{^CM^D}\times{^DA^E}\to{^CM^E}$. Strictly speaking, $\boxtimes_D$ also depends on the action $\triangleleft$, but it will always be clear from the context which action we have in mind. When $D$ is the category $1$ (one object, one morphism), we simply write $\boxtimes$ for $\boxtimes_D$. Notice also that $C$ or $E$ or both could be the category $1$, in which case we omit them. 

\begin{remark}\label{remark.right action}
	More generally, one can define in this way a functor $\boxtimes_D:{^CM^{D\times P}}\times{^DA^{E\times Q}}\to{^CM^{E\times P\times Q}}$.
\end{remark}

	For us, $A$ will always be monoidal, and $\triangleleft$ will be unital and associative in the obvious sense (there is an associator realizing a natural isomorphism between the two possible ways of going from $M\times A\times A$ to $M$, etc.). From now on, whenever the $\boxtimes_D$ construction appears, we make the assumptions that $M$ is presentable, $A$ is a 2-ring, and the action $\triangleleft$ is separately cocontinuous. In many examples $A=\cat{Set}$, and the action of a set $S$ on an object $m\in M$ is simply the {\it copower} $m\cdot S$: the coproduct of copies of $m$ indexed by the set $X$. This is usually denoted by $S\cdot m$ in the literature (e.g. \cite[III.3]{MR1712872}), but we reverse the order to keep the notation consistent with the use of right actions. On the other hand, if $a\in A$, the $S$-indexed copower of $a$ might be regarded as the result of a right action of $a$ on $S\in\cat{Set}$, in which case it will be written $S\cdot a$.    

	We observed before (proof of \cref{lemma.^GM}) that as a consequence of \cref{thm.locpres2}, there is an equivalence of categories $\Hom_{\cat{2AbGp}}({^G\cat{Set}}, A) = A^G$.  In the proof of \cref{thm.locpres2} we used the language of Kan extensions, but when $G$ has no extra sketch structure, we can describe this equivalence much more explicitly.
	
\begin{proposition}\label{prop.boxtimes_G}
  Let $G$ be a small category.  The equivalence $\Hom_{\cat{2AbGp}}({^G\cat{Set}} \to A) \cong A^G$ from \cref{thm.locpres2} is realized as follows.  Any morphism $\alpha:{^G\cat{Set}} \to A$ determines a morphism $\alpha^G : {^G\cat{Set}^G} \to A^G$ by post-composition, and hence a distinguished object $\alpha({\Hom_{G}(-,-)}) \in A^G$.   Conversely, any object $x \in A^G$ determines a 1-morphism $x\boxtimes_G : {^G\cat{Set}} \to A$. \qedhere
\end{proposition} 

\begin{definition}\label{defn.reg}
  Let $A$ be a 2-abelian group and $G$ a small category.  The category $^GA$ consists of the \define{left $G$-representations in $A$}, and $A^G$ are the \define{right $G$-representations in $A$}.  
  The \define{regular representation} of a category $G$ is the object ${^GG^G} \in {^G\cat{Set}^G}$ corresponding to $\Hom_G : G^\op \times G \to \cat{Set}$.  Let $G_0$ denote the category with no nontrivial morphisms on the set of objects in $G$; the inclusion $G_0 \mono G$ induces forgetful functors $U: {^G\cat{Set}} \to {^{G_0}\cat{Set}}$ and $\cat{Set}^G \to \cat{Set}^{G_0}$.  Since $G_0$ has no nontrivial morphisms, there is a canonical identification $^{G_0}A = A^{G_0}$ for any 2-abelian group $A$, and we call its objects \define{$G_0$-sorted objects in $A$}.  The \define{left (resp.\ right) regular representation of $G$} is the object $^GG$ ($G^G$) in $^G\cat{Set}^{G_0}$ ($^{G_0}\cat{Set}^G$) formed from $^GG^G$ by forgetting the right (left) action.
\end{definition} 

	When it is clear from the context which action of $G$ on itself we have in mind, we might just write $G$ instead of $^GG$ or $G^G$. 

	Let $G$ be any small category.  Then $^G\cat{Set}$ is a 2-ring with $\otimes = \times$, and so every object in $^G\cat{Set}$ is a (coassociative, counital) coalgebra in a unique way, every object is cocommutative as a coalgebra, and every morphism is a morphism of coalgebras.  In particular, the regular representation $^GG^G \in {^G\cat{Set}^G} = {^{G\times G^\op}\cat{Set}}$ is a coalgebra.  Let $A$ be a 2-ring and $\alpha: {^G\cat{Set}} \to A$ a 1-morphism of 2-rings.  Then $\alpha^G : {^G\cat{Set}^G} \to A^G$ is also symmetric monoidal.  We write $\hat x = \hat x(\alpha)$ for the cocommutative coalgebra ${\alpha^G}({^GG^G}) \in A^G$, and $\varepsilon_{\hat x}$ and $\Delta_{\hat x}$ for its counit and comultiplication, respectively. We sometimes suppress the exponent in $\alpha^G$ and simply write $\hat x = \alpha(^GG^G)$. 

	For a cocomplete category $M$ and a small category $G$, any object $x\in M^{G_0}$ can be regarded as an object in $^{G_0}M^{G_0}$ in an obvious way: simply set $x(s,s)=x(s)$ and $x(t,s)$ to be the initial object $0_M$ when $s\ne t\in G_0$. If one then regards $G$ as an object in $^{G_0}\cat{Set}^G$, it makes sense to define $x\boxtimes_{G_0} G\in{^{G_0}M^G}$. 

	Define $\pi_0=\pi_0(G)\in{^{G_0}\cat{Set}^{G_0}}$ to be the equivalence relation on $G_0$ whose classes are the connected components of $G$. If $A$ is a 2-ring right-acting on $M$ and $y\in A^{G_0}$, one can regard $y$ as an object of $^{G_0}A^{G_0}$ (just as $x$ above), so we can construct $x\boxtimes_{G_0}^{\pi_0}y=x\boxtimes_{G_0}\pi_0\boxtimes_{G_0} y\in{^{G_0}M^{G_0}}$. If $y=U(\hat y)$ for some $\hat y\in A^G$, notice that the right action of $G$ on $y$ descends to a right action on $x\boxtimes_{G_0}^{\pi_0}y$, despite the fact that $\hat y$ cannot be regarded as an object of $^{G_0}A^G$, in general. We use the notation $x\boxtimes_{G_0}^{\pi_0}\hat y$ for the resulting object in $^{G_0}M^G$. Simply put, $x\boxtimes_{G_0}^{\pi_0}\hat y(t,s)$ is $x(t)\triangleleft\hat y(s)$ if $G(t,s)\ne\emptyset$, and $0_M$ otherwise. This gives a functor $x\boxtimes_{G_0}^{\pi_0}:A^G\to{^{G_0}M^G}$, and if $\hat x\in M^G$, one defines a functor $\hat x\boxtimes_G^{\pi_0}:A^G\to{M^{G\times G}}$. Similarly, there are functors $\boxtimes_{G_0}^{\pi_0}\hat y:M^{G_0}\to{^{G_0}M^G}$ and $\boxtimes_G^{\pi_0}\hat y:M^G\to M^{G\times G}$. These will come up again soon. 

	Now let $\alpha : {^G\cat{Set}} \to A$ be a 1-morphism of 2-rings, $\hat x = \alpha(^GG^G)$, and $x = U(\hat x) \in A^{G_0}$ its associated sorted object. Let $\tau_x : x\boxtimes_{G_0}G \to x\boxtimes \hat x$ be the arrow defined on the summand $x(\ell(g))=x(\ell(g))\cdot\{g\}$ of $(x\boxtimes_{G_0}G)(\ell(g),r(g))=x(\ell(g))\cdot G(\ell(g),r(g))$ by
	\begin{equation*}
		x(\ell(g)) \verylongarrow{\Delta}
  	x(\ell(g))\otimes x(\ell(g))\verylongarrow{\id\otimes x(g)} x(\ell(g))\otimes x(r(g)).
	\end{equation*}
Then $\tau_x$ factors through $x\boxtimes_{G_0}^{\pi_0}\hat x$, and from now on it is the resulting map $x\boxtimes_{G_0}G \to x\boxtimes_{G_0}^{\pi_0}\hat x$ that we will refer to as $\tau_x$. The following observation will be very important below:

\begin{lemma}\label{lemma.Galoiscond}
  With the above notations, the map $\tau_x:x\boxtimes_{G_0}G\to x\boxtimes_{G_0}^{\pi_0}\hat x$ is an isomorphism in $^{G_0}A^G$.
\end{lemma}
\begin{proof}
  When $A = {^G\cat{Set}}$ and $\alpha = \id$, the statement is equivalent to $G$ being a groupoid. This gives an isomorphism $\tau=\tau_{^GG}$ in $^{G_0}({^G\cat{Set}})^G$. Now notice that $\tau_x$ is the image of $\tau$ through the functor $^{G_0}\alpha^G:{^{G_0}}({^G\cat{Set}})^G\to {^{G_0}}A^G$ induced by $\alpha$, and hence must also be an isomorphism.   
\end{proof}


\Section{The fundamental pro-groupoid of an affine 2-scheme} \label{section.pi1construction}

\subsection{\texorpdfstring{Definition of $\pi_1$}{Definition of pi_1}}\label{section.pi1}

	We introduce the following notions, which play a central role in the paper:

\begin{definition}\label{defn.pi1}
	Let $A$ be a commutative 2-ring. The \define{fundamental} or \define{Galois pro-groupoid of $A$} is the 2-functor $\cat{Gpd}\to\cat{CAT}$ from small groupoids to possibly large but locally small categories defined, for any small groupoid $G$, by:
	\[
		\pi_1(A)(G) = \Hom_{A\cat{-ComAlg}}(^GA, A).
	\]         
	The \define{fundamental pro-finite groupoid of $A$} is the 2-functor $\pi_1^f(A):\cat{Gpd}_f\to\cat{CAT}$ from finite groupoids to categories defined by the same formula. 
\end{definition} 
\begin{remark}
    The naturality of \cref{defn.pi1} in $G$ is a simple matter, and is left to the reader. One could also allow $G$ to be an arbitrary small category, but we will typically consider only groupoids, as many of the methods used do not work in the general case.  We discuss what works and what fails in \cref{section.nongroupoids}.
\end{remark}
\begin{remark}\label{remark.CAT->GPD}
	We will see later (\cref{thm.iso} and discussion below) that $\pi_1(A)$ always lands in the 2-category $\cat{GPD}$ of possibly large but locally small \define{groupoids}, and that in fact we can restrict the codomain of $\pi_1(A)$ even further to $\cat{Gpd}$ (\cref{prop.smallness}). We will frequently abuse terminology and notation and denote these restricted versions of $\pi_1$ by the same symbols. 
\end{remark}

	The main result of the paper (\cref{thm.Gal}) is that under certain reasonable conditions on the 2-ring $A$, the 2-functors $\pi_1(A)$ and $\pi_1^f(A)$ are pro-representable. We will then refer to their representing objects (which are pro-objects of the (2,1)-categories $\cat{Gpd}$ and $\cat{Gpd}_f$) as \define{being} $\pi_1(A)$ and $\pi_1^f(A)$ respectively. This justifies the choice of name for the two versions of $\pi_1$. We now start working our way towards that goal, by making several simplifications and reformulations.


\subsection{\texorpdfstring{$G$-torsors in a 2-ring}{G-torsors in a 2-ring}}\label{section.torsor}

	Let $G$ be a small groupoid and $A$ a commutative 2-ring.  In view of \cref{lemma.main}, we could just as well have defined $\pi_1(A)$ by $G\mapsto\Hom_{\cat{Com2Ring}}({^G\cat{Set}},A)$. In this section we describe the 2-ring homomorphisms $^G\cat{Set} \to A$ as \define{right $G$-torsors in $A$}, and also show that morphisms of $G$-torsors are isomorphisms, and hence that $\pi_1(A)$ takes values in $\cat{GPD}$ (\cref{remark.CAT->GPD}).  We continue the notation of \cref{subsection.coends}.

\begin{definition}\label{defn.torsor}
  Let $G$ be a groupoid and $A$ a commutative 2-ring.  Let $\hat y \in A^G$ be a cocommutative coalgebra, and $y = U(\hat y)\in A^{G_0}$ its associated sorted object. Using the notations introduced before \cref{lemma.Galoiscond}, $\hat y$ is said to be a \define{right $G$-pseudotorsor in $A$} if $\tau_y$ is an isomorphism. If in addition the counit $\hat y\to 1_A$ is a colimiting cocone, we call $\hat y$ a \define{right G-torsor}. The category of $G$-torsors in $A$ with maps of coalgebras in $A^G$ as morphisms is denoted by $\cat{Tors}(G,A)$.
\end{definition}

\begin{remark}
	We are trying to imitate the usual notion of $G$-torsor for a group $G$ (see \cite[VIII.2]{MacLane1994}): the condition that $\tau_y$ be an isomorphism is the correct ``right principality'' condition and is the analogue of \cite[VIII.2 Definition 6 (ii)]{MacLane1994}, while the condition on the counit is an analogue of sorts for \cite[VIII.2 Definition 6 (i)]{MacLane1994}. In fact, it can be shown that when the 2-ring in question is a Grothendieck topos (cf. \cref{eg.2rings}) and $G$ is a group, our definition agrees with Mac Lane and Moerdijk's. 
\end{remark}

	The following result provides further justification for the second condition imposed on a torsor in \cref{defn.torsor}.

\begin{lemma}\label{lemma.colcocone}
  Let $A$ be a 2-ring, $G$ a small category, and $\alpha$ and $\hat x=\hat x(\alpha)$ as above. Regard the counit $\varepsilon_{\hat x}:\hat x\to 1_{A^G}$ as a cocone from the functor $\hat x:G^\op\to A$ to the monoidal unit $1_A\in A$. This cocone is a colimit.
\end{lemma}
\begin{proof}
  Since $\alpha$ is cocontinuous, it suffices to check the case when $A = {^G\cat{Set}}$ and $\alpha = \id$.  Then $\hat x = {^GG^G}$ is the diagram in $^G\cat{Set}$ whose $s$th entry is $G(-,s) \in {^G\cat{Set}}$.  Since colimits of functors are computed pointwise, it suffices to show that $\colim_{s\in G^\op} G(t,s) = 1$ for each $t\in G$. But this is clear from the fact that every element of $G(t,s)$ is the image of $\id:t\to t$ through a map $G(t,t)\to G(t,s)$ induced by a morphism in $G$.     
\end{proof}

	We now introduce some more notation. Recall that for a small groupoid $G$ (or more generally a small category), we are denoting by $\ell=\ell_G$ and $r=r_G$ the target and source map  from the set of arrows of $G$ (also denoted by $G$) to $G_0$. Define $\hat\Lambda\in {^G}\cat{Set}^{G\times G}$ by $\hat\Lambda=G^G\boxtimes_{G_0}G^G$. In other words, $\hat\Lambda(t;s,s')$ is the set of pairs of arrows with common target $t$ with sources $s$ and $s'$. The three $G$-actions are the ones coming from multiplication in $G$. We denote by $\Lambda\in {^G}\cat{Set}^{G_0\times G}$ the object obtained from $\hat\Lambda$ by forgetting the right action of $G$ on $s$ in $\Lambda(t;s,s')$. 

	Recall the notations $\boxtimes_{G_0}^{\pi_0}$ and $\boxtimes_G^{\pi_0}$ introduced before \cref{lemma.Galoiscond}. Let $A$ be a commutative 2-ring, $G$ a small groupoid, and $\hat x\in A^G$ a cocommutative coalgebra. The maps 
  \[
  	\coprod_{t\in G_0}\hat x(t)\cdot\hat\Lambda(t;s,s')\to \hat x(s)\otimes\hat x(s')
  \]
defined on the summand $\hat x(t)$ corresponding to the pair $(g,g')$ in $\hat\Lambda(t;s,s')$ by 
  \[
  	\hat x(t) \overset{\Delta_{\hat x}}{\xrightarrow{\hspace*{1.2cm}}} \hat x(t)\otimes\hat x(t) \overset{g\otimes g'}{\xrightarrow{\hspace*{1.2cm}}} \hat x(s)\otimes\hat x(s') 
  \]
induce a morphism $\eta_{\hat x}:\hat x\boxtimes_G\hat\Lambda \to \hat x\boxtimes_G^{\pi_0}\hat x$ in $A^{G\times G}$. 

\begin{lemma}\label{lemma.etaiso}
    With the notations above, if $\hat x$ is a $G$-pseudotorsor in $A$, $\eta_{\hat x}$ is an isomorphism in $A^{G\times G}$.     
\end{lemma}
\begin{proof}
	Let $x$ be, as before, the $G_0$-sorted object associated to $\hat x$. Forgetting one of the right actions of $G$, $\eta_{\hat x}$ induces a map $\eta_x:\hat x\boxtimes_G\Lambda\to x\boxtimes_{G_0}^{\pi_0}\hat x$ in $^{G_0}A^G$. In order to prove that $\eta_{\hat x}$ is an isomorphism in $A^{G\times G}$, it is enough to prove that $\eta_x$ is an isomorphism in $^{G_0}A^G$. Indeed, one checks easily that since $\eta_{\hat x}$ respects the right action of $G\times G$, an inverse for $\eta_x$ automatically does so too.

	We now claim that we have an identification $\hat x\boxtimes_G\Lambda \cong x\boxtimes_{G_0}G$ so as to make the following diagram commutative:
    \begin{equation}\label{eqn.taueta}
        \tikz[baseline=(current bounding box.center)]{
          \path (0,0) node (xDelta) {$\hat x\boxtimes_G\Lambda$} +(2,0) node (Gx) {$x\boxtimes_{G_0}G$} +(4,0) node (xx) {$x\boxtimes_{G_0}^{\pi_0}\hat x$.};
          \draw[->] (xDelta) -- (Gx) node[pos=.5,auto] {$\scriptstyle \cong$};
          \draw[->] (Gx) -- (xx) node[pos=.5,auto] {$\scriptstyle \tau_x$};
          \draw[->,bend right=15] (xDelta.south east) to node[pos=.5,auto,swap] {$\scriptstyle \eta_x$} (xx.south west); 
         }
     \end{equation}
In view of the previous paragraph, this claim would finish the proof. 

	To prove the claim, notice first that sending the pair $(g,g')\in\Lambda(t;s,s')$ of arrows in $G$ with common target $t$ and sources $s$ and $s'$ to the pair $(g,g^{-1}g')\in G_2$ (the set of composable pairs of arrows in $G$) is an isomorphism $\Lambda\cong G_2$ in $^G\cat{Set}^{G_0\times G}$. Here, $G_2$ is viewed as an object in $^G\cat{Set}^{G_0\times G}$ by $G_2(t;s,s')=$ pairs $(g,h)$ such that $g$ has source and target $s$ and $t$, while $h$ has source and target $s'$ and $s$ respectively. The left (right) action of $G$ on such a pair $(g,h)$ is by composition with $g$ (resp. $h$). The analogue for groups is the familiar observation that if $H$ is a group, then the diagonal action by left translation and the left translation action on just the left hand factor are isomorphic left $H$-set structures on $H\times H$. 

  Applying the functor $\hat x\boxtimes_G$ to the isomorphism $\Lambda\cong G_2$ just constructed, one has an isomorphism $\hat x\boxtimes_G\Lambda \cong \hat x\boxtimes_G G_2$ in $^{G_0}A^G$. Now note that $G_2\cong {^GG}\boxtimes_{G_0}G^G$ as objects of $^G\cat{Set}^{G_0\times G}$. We thus get the following series of isomorphisms in $^{G_0}A^G$:
	\[
  	\hat x\boxtimes_G\Lambda\ \ \cong\ \ \hat x\boxtimes_G G_2\ \ \cong\ \ \hat x\boxtimes_G[{^GG}\boxtimes_{G_0}G^G]\ \ \cong\ \ x\boxtimes_{G_0}G. 
  \]
It is straightforward now to unwrap the definitions and check that the composition of all of these makes \cref{eqn.taueta} commutative. 
\end{proof}

	The following proposition is in a sense well known, and has appeared in many guises and in many settings before; \cite[Theorem VIII.2.7]{MacLane1994} and \cite[Theorem 1.2]{Ulbrich1989} are two such examples, in the context of torsors in topoi and of Hopf-Galois objects, respectively. 
        
\begin{proposition}\label{prop.GSet -> A are torsors}
  Let $G$ be a small groupoid, $A$ a commutative 2-ring, and $\alpha: {^G\cat{Set}} \to A$ a 1-morphism of commutative 2-rings.  Then $\hat x = \alpha({^G}G^G) \in A^G$ is a right $G$-torsor, and $\alpha\mapsto\alpha({^G}G^G)$ defines an equivalence $\Hom_{\cat{Com2Ring}}(^G\cat{Set}, A)\simeq\cat{Tors}(G,A)$.
\end{proposition}
\begin{sketch}            
	We noted in the discussion preceding \cref{lemma.colcocone} that $\hat x$ is a cocommutative coalgebra. The fact that it is actually a torsor follows from \cref{lemma.colcocone} and \cref{lemma.Galoiscond}, each of which takes care of one of the two conditions a torsor has to satisfy. 

  Now let $\hat x$ be a $G$-torsor in $A$, and consider the functor $\hat x\boxtimes_G:{^G\cat{Set}}\to A$. We have to check that (a) the torsor structure on $\hat x$ makes $\hat x\boxtimes_G$ a symmetric monoidal functor in a natural way, (b) the canonical isomorphism $\hat x\cong \hat x(\hat x\boxtimes_G)$ in $A^G$ from \cref{prop.boxtimes_G} is actually an isomorphism of coalgebras in $A^G$, and similarly, (c) that the isomorphism $\alpha\cong \hat x(\alpha)\boxtimes_G$ in $\Hom_{\cat{2AbGp}}({^G}\cat{Set},A)$ is  monoidal. These statements are more or less routine. For this reason, we only partially explain (a) by indicating how the monoidal structure is constructed on the functor $\hat x\boxtimes_G$, and leave the rest to the reader.  
    
	To construct an isomorphism $\hat x\boxtimes_G(1_{^G\cat{Set}})\cong 1_A$, notice first that the former is, by the very definition of $\boxtimes_G$, the colimit of the functor $\hat x:G^\op\to A$. By one of the two torsor conditions for $\hat x$, $\varepsilon_{\hat x}\otimes_G\id$ will be an isomorphism.
    
    We now want an isomorphism $\hat x\boxtimes_G(S\otimes T)\cong(\hat x\boxtimes_G S)\otimes(\hat x\boxtimes_G T)$, functorial in $S,T\in{^G\cat{Set}}$. For objects $t,s,s'\in G_0$, a pair $(g,g')\in\hat\Lambda(t;s,s')$ and elements $a\in S(s)$, $a'\in T(s')$, send $((g,g'),(a,a'))$ to $(ga,g'a')\in (S\otimes T)(t)$. It is easily checked that this induces an isomorphism $\zeta:\hat\Lambda\boxtimes_{G\times G}(S\boxtimes T)\to S\otimes T$ in $^G\cat{Set}$.     
The composition 
	\begin{align*}
  	\hat x\boxtimes_G(S\otimes T) & \veryverylongarrow{\id\boxtimes_G\zeta^{-1}}
   	(\hat x\boxtimes_G\hat\Lambda)\boxtimes_{G\times G}(S\boxtimes T)\\ &\veryverylongarrow{\eta_{\hat x}\boxtimes_{G\times G}\id}
    (\hat x\boxtimes_G^{\pi_0}\hat x)\boxtimes_{G\times G}(S\boxtimes T)\\ &\veryverylongarrow{\cong}
    (\hat x\boxtimes_G S)\otimes(\hat x\boxtimes_G T)
	\end{align*}
is the desired isomorphism, completing the construction of the monoidal structure on the functor $\hat x\boxtimes_G$. The second isomorphism in this series comes from \cref{lemma.etaiso}, while the last one makes use of the symmetry in $A$.  
\end{sketch}
\begin{remark}
	This result says, essentially, that $^G\cat{Set}$ is the universal 2-ring containing a right $G$-torsor. 
\end{remark}

	We call $\alpha(^GG^G)$ \define{the right $G$-torsor associated to $\alpha$}, and henceforth we freely transition between commutative 2-ring morphisms $^G\cat{Set} \to A$ and right $G$-torsors in $A$.	As promised at the beginning of this subsection, our aim will now be to show that morphisms of torsors are automatically isomorphisms. For this, we need some preparations. 

\begin{definition}
  A functor is \define{conservative} if it reflects isomorphisms.
\end{definition}
\begin{lemma}\label{lemma.torsors are conservative}
 Let $G$ be a groupoid, $A$ a 2-ring, and $\alpha:{^G\cat{Set}}\to A$ a morphism of 2-rings.  Let $\hat x = \alpha(^GG^G) \in A^G$, and let $x\in A^{G_0}$ be its associated sorted object. Then the functor $x\boxtimes :A\to A^{G_0}$ is conservative. 
\end{lemma}

\begin{proof}
  The functor $x\boxtimes $ factors as $\hat x\boxtimes:A\to A^G$, followed by the forgetful functor $U:A^G\to A^{G_0}$. Since $U$ is conservative, it suffices to show that $\hat x\boxtimes $ is conservative. 

  Let $f:a\to b$ be an arrow in $A$, and assume that $\hat x\boxtimes f$ is an isomorphism in $A^G$.  Recall that $\varepsilon_{\hat x}:\hat x\to 1$ is a colimiting cone of $\hat x$.  Since tensoring with a fixed object preserves colimits, the cocone $\hat x\boxtimes a\to 1\otimes a = a$ is a colimit.  
Together with its analogue for $b$, this cone fits into the following commutative square:
 \[ \tikz{
  \path (0,0) node (xa) {$\hat x\boxtimes a$} +(2,0) node (xb) {$\hat x\boxtimes b$} +(0,-2) node (a) {$a$} +(2,-2) node (b) {$b$};
  \draw[->] (xa) -- (xb) node[pos=.5,auto] {$\scriptstyle \hat x\boxtimes f$};
  \draw[->] (a) -- (b)  node[pos=.5,auto] {$\scriptstyle f$};
  \draw[->] (xa) -- node[auto,swap] {$\scriptstyle  \varepsilon_{\hat x} \boxtimes a$} (a);
  \draw[->] (xb) -- node[auto] {$\scriptstyle  \varepsilon_{\hat x} \boxtimes b$} (b);
 } \]
Thus $f$ is the unique arrow between the colimits of $\hat x \boxtimes a$ and $\hat x\boxtimes b$ induced by the natural transformation $\hat x\boxtimes f$. But since the latter is an isomorphism, so is $f$.
\end{proof}

The following definition imitates the usual construction of the cotensor product of comodules from e.g.\ \cite[Section 10]{Brzezinski2003}.
\begin{definition}
  Let $x$ be a coalgebra in a monoidal category $V$. Denote by $V^x$ and $^xV$ the categories of right and respectively left $x$-comodules in $V$, and let $m:C\to V^x$ and $n:D\to {^xV}$ be functors, where $C,D$ are small categories. Regarding $x$ as a functor from the category $1$ to $V$, we can talk about $m\boxtimes x:C\to V$, $m\boxtimes x\boxtimes n:C\times D\to V$, and so on. We refer to $m$ (resp.\ $n$) as a $C$- (resp.\ $D$-) indexed right (resp.\ left) $x$-comodule, and we think of the usual comodules over $x$ as comodules indexed over the category $1$. Finally, let $\rho_m:m\to m\otimes x$ be the natural transformation inducing the right comodule structure on each $m(c)$, $c\in C$, and define $\rho_n$ similarly. The \define{cotensor product} $m\square_x n$ is  the following equalizer in $^{C\times D}V$, if it exists:
 \[\tikz[anchor=base]{
  \path (0,0) node (msqn) {$m\square_x n$} +(2,0) node (mn) {$m\boxtimes n$} +(5,0) node (mxn) {$m\boxtimes x\boxtimes n$};
  \draw[->] (msqn) -- (mn);
  \draw[transform canvas={yshift=0.8mm},->] (mn) -- (mxn) node[pos=.5,auto] {$\scriptstyle \rho_m\boxtimes\id_n$};
  \draw[transform canvas={yshift=-0.8mm},->] (mn) -- (mxn) node[pos=.5,auto,swap] {$\scriptstyle \id_m\boxtimes\rho_n$}; 
 }\]
We usually work with categories $V$ which have the necessary limits, and we assume from now on that choices of such equalizers have been made, so as to induce a functor $\square_x: {^C}(V^x)\times {^D}(^xV)\to {^{C\times D}}V$.
\end{definition}

\begin{lemma}\label{lemma.cotensor}
  Let $V$ be a symmetric monoidal category, and $f:x\to y$ a map of coalgebras in $V$.  We denote by $\rho_x:x\to x\otimes y$ the map $(\id_x\otimes f)\circ\Delta_x$, making $x$ into a right $y$-comodule.
  Then the diagram
 \[\tikz[anchor=base]{
  \path (0,0) node (x) {$x$} +(3,0) node (xy) {$x\otimes y$} +(7,0) node (xyy) {$x\otimes y\otimes y$};
  \draw[->] (x) -- (xy) node[pos=.5,auto] {$\scriptstyle \rho_x$};
  \draw[transform canvas={yshift=0.8mm},->] (xy) -- (xyy) node[pos=.5,auto] {$\scriptstyle \rho_x\otimes\id_y$};
  \draw[transform canvas={yshift=-0.8mm},->] (xy) -- (xyy) node[pos=.5,auto,swap] {$\scriptstyle \id_x\otimes\Delta_y$}; 
 }\]
 exhibits $x$ as the cotensor product $x\square_yy$. 
\end{lemma}
\begin{proof}
	For ease of notation, we relabel the maps in the above diagram as
		\[\tikz[anchor=base]{
  		\path (0,0) node (x) {$x$} +(3,0) node (xy) {$x\otimes y$} +(7,0) node (xyy) {$x\otimes y\otimes y$.};
  		\draw[->] (x) -- (xy) node[pos=.5,auto] {$\scriptstyle e$};
  		\draw[transform canvas={yshift=0.8mm},->] (xy) -- (xyy) node[pos=.5,auto] {$\scriptstyle \partial_1$};
  		\draw[transform canvas={yshift=-0.8mm},->] (xy) -- (xyy) node[pos=.5,auto,swap] {$\scriptstyle \partial_0$}; 
 		}\]
Now embed this diagram into
	\[\tikz[anchor=base]{
  	\path (0,0) node[inner ysep=0pt] (x) {$x$} +(3,0) node[inner ysep=0pt] (xy) {$x\otimes y$} +(7,0) node[inner ysep=0pt] (xyy) {$x\otimes y\otimes y$,};
  	\draw[->] (x) -- (xy) node[pos=.5,auto] {$\scriptstyle e$};
  	\draw[transform canvas={yshift=0.8mm},->] (xy) -- (xyy) node[pos=.5,auto,inner ysep=1pt] {$\scriptstyle \partial_1$};
  	\draw[transform canvas={yshift=-0.8mm},->] (xy) -- (xyy) node[pos=.5,auto,swap, inner ysep=1pt] {$\scriptstyle \partial_0$};
  	\draw[<-,bend right=15] (x.south east) to node[pos=.5,auto,swap] {$\scriptstyle u$} (xy.south west);
  	\draw[<-,bend right=55] (xy.south east) to node[pos=.5,auto,swap] {$\scriptstyle v$} (xyy.south west); 
 	}\]
where $u=\id_x\otimes\varepsilon_y$ and $v=\id_{x\otimes y}\otimes\varepsilon_y$.
Note that we have the identities
	\[
  	\partial_0e=\partial_1e,\quad ue=\id_{x\otimes y},\quad v\partial_0=\id_{x\otimes y\otimes y},\quad v\partial_1=eu. 
 	\]
These make the diagram (dual to) a split fork in the sense of \cite[Section VI.6]{MR1712872}, hence an equalizer diagram by the dual of \cite[Lemma in VI.6]{MR1712872}.   
\end{proof}
\begin{remark}\label{remark.cotensor}
	The proof of \cref{lemma.cotensor} is more important than the statement itself. We obtained $\rho_x:x\to x\otimes y$ as an equalizer in a split fork --- a \define{split equalizer} (the notion dual to Mac Lane's ``split coequalizers'' \cite[Section VI.6]{MR1712872}). Whenever this happens, we say that $x$ has been realized as a \define{split cotensor product} $x\square_yy$. 
\end{remark}

	Let $A$ be a 2-ring, $G$ a groupoid, and $f: \hat x \to \hat y$ a map of coalgebras in $A^G$.  This induces a map (also denoted by $f$) $x\to y$ of coalgebras in $^{G_0}A$, and hence the structure $\rho_x:x\to x\boxtimes_{G_0} y$ of a right $y$-comodule on $x$ (note that the tensor product in $A^{G_0}$ is exactly $\boxtimes_{G_0}$, introduced in \cref{subsection.coends}; since either $y$ or $x$ can be regarded as an object in $^{G_0}A^{G_0}$, $x\boxtimes_{G_0}y$ has a natural $G_0$-sorted object structure). Finally, note that $y\boxtimes_{G_0}^{\pi_0}\hat y\in{^{G_0}A^G}$ is naturally a left $y$-comodule indexed by $G^\op$. 
\begin{lemma}\label{lemma.cotensor1}
	$x\boxtimes_{G_0}^{\pi_0}\hat y$ can be realized as a split cotensor product $x\square_y(y\boxtimes_{G_0}^{\pi_0}\hat y)$. 
\end{lemma}
\begin{proof}
	Apply the functor $\boxtimes_{G_0}^{\pi_0}\hat y$ to the diagram in the statement of \cref{lemma.cotensor} for $V= {^{G_0}A}$ and $\otimes=\boxtimes_{G_0}$, and use the fact that functors preserve split forks. 
\end{proof}
 
	Just like $y\boxtimes_{G_0}^{\pi_0}\hat y$, $y\boxtimes_{G_0}G$ is a $G^\op$-indexed left $y$-comodule, so we can now talk about the cotensor product $x\square_y(y\boxtimes_{G_0}G)$. Since split forks are preserved by functors, \cref{remark.cotensor} shows that 
	\[\tikz[anchor=base]{
  	\path (0,0) node (x) {$x\boxtimes_{G_0}G$} +(5,0) node (xy) {$x\boxtimes_{G_0} y\boxtimes_{G_0} G$} +(12,0) node (xyy) {$x\boxtimes_{G_0} y\boxtimes_{G_0} y\boxtimes_{G_0} G$};
  	\draw[->] (x) -- (xy) node[pos=.5,auto] {$\scriptstyle \rho_x\boxtimes_{G_0} G$};
  	\draw[transform canvas={yshift=0.8mm},->] (xy) -- (xyy) node[pos=.5,auto] {$\scriptstyle \rho_x\boxtimes_{G_0} \id_y\boxtimes_{G_0} G$};
  	\draw[transform canvas={yshift=-0.8mm},->] (xy) -- (xyy) node[pos=.5,auto,swap] {$\scriptstyle \id_x\boxtimes_{G_0} \Delta_y\boxtimes_{G_0} G$}; 
 	}\]
is a split equalizer. All together, we have proved:
\begin{lemma}\label{lemma.cotensor2}
	The map $w:x\boxtimes_{G_0} G\to x\boxtimes_{G_0} y\boxtimes_{G_0} G$ whose component at $(t,s)\in G_0\times G_0$ is defined by
 		\[
  		(x\boxtimes_{G_0} G)(t,s)=\coprod_{g\in G(t,s)}x(t)\verylongarrow{\coprod_g\Delta_{x(t)}}\coprod_g x(t)\otimes x(t)\verylongarrow{\coprod_g\id\otimes f} x(t)\otimes \coprod_g y(t)=[x\boxtimes_{G_0} y\boxtimes_{G_0} G](t,s) 
 		\]
exhibits $x\boxtimes_{G_0} G$ as a split cotensor product $x\square_y(y\boxtimes_{G_0} G)$.\qedhere   
\end{lemma}

	Finally, we come to the main result of the section:
\begin{theorem}\label{thm.iso}
	Let $A$ be a commutative 2-ring, and $G$ a groupoid. If $\alpha$ and $\beta$ are 1-morphisms of 2-rings from $^G\cat{Set}$ to $A$ and $\chi:\alpha\Rightarrow \beta$ is a 2-morphism, then $\chi$ is an isomorphism.  
\end{theorem}
	With the preparations we have made, this is now a simple matter of imitating the usual proofs that ``maps of torsors are automatically isomorphisms'' (e.g. \cite[Proposition 1.6]{Bichon2010}). 
\begin{proof}
	We use the notation introduced above: $\hat x\in A^G$ and $x\in {^{G_0}A}$ are the torsor and respectively the $G_0$-sorted object associated to $\alpha$, and $\hat y$ and $y$ are similarly associated to $\beta$. The 2-morphism $\chi$ induces a map $f:\hat x\to\hat y$, which is an isomorphism in the category of coalgebras in $A^G$ if and only if $\chi$ is an isomorphism.  It suffices to check that $Uf: x\to y$ is an isomorphism in $A^{G_0}$.

	Consider the following diagram in $^{G_0}A^G$:
 		\[\tikz{
  		\path (0,0) node (Gx) {$x\boxtimes_{G_0} G$} +(5,0) node (xy) {$x\boxtimes_{G_0}^{\pi_0}\hat y$} +(0,-3) node (xGy) {$x\square_y(y\boxtimes_{G_0} G)$} +(5,-3) node (xyy) {$x\square_y(y\boxtimes_{G_0}^{\pi_0} \hat y)$};
  		\draw[->] (Gx) -- (xy) node[pos=.5,auto] {$\scriptstyle (id_x\boxtimes_{G_0}^{\pi_0} f)\circ\tau_x$};
  		\draw[->] (xGy) -- (xyy) node[pos=.5,auto] {$\scriptstyle \id_x\square_y\tau_y$};
  		\draw[->] (Gx) -- (xGy) node[pos=.5,auto,swap] {$\scriptstyle w$};
   		\draw[->] (xy) -- (xyy) node[pos=.5,auto] {$\scriptstyle \rho_x\boxtimes_{G_0}^{\pi_0}\id_{\hat y}$};
 		}\]
One checks easily that it is commutative. We know that $\tau_x$ and $\tau_y$ are isomorphisms from the discussion right before \cref{lemma.torsors are conservative}. We also know that the vertical arrows are isomorphisms: \cref{lemma.cotensor1} states this for the right hand arrow, and \cref{lemma.cotensor2} for the left hand one. It follows then that $\id_x\boxtimes_{G_0}^{\pi_0} f$ must be an isomorphism. But then every component $f_s:x(s)\to y(s)$ is an isomorphism by \cref{lemma.torsors are conservative}, or, in other words, $Uf$ is an isomorphism, and we are done.
\end{proof}
\begin{remark}\label{remark.pi1GPD}
	The theorem shows that the 2-functor $\pi_1:\cat{Gpd}\to\cat{CAT}$ actually takes values in the full sub-2-category $\cat{GPD}$ of $\cat{CAT}$ consisting of possibly large but locally small groupoids.  
\end{remark}


\subsection{Cartesian 2-rings} \label{section.cartesian}

	In this subsection we simplify matters further by restricting our attention to a special class of commutative 2-rings:    
\begin{definition}
    A 2-ring $A$ is \define{cartesian} if its monoidal structure is the one given by the binary product and terminal object. 
\end{definition}
\begin{example}
	(a) We observed before in \cref{eg.2rings} that Grothendieck topoi (for example, $^G\cat{Set}$) are cartesian 2-rings.   
    
	(b) Any complete Heyting algebra, i.e.\ a complete distributive lattice which is cartesian closed when regarded as a category with product given by the meet operation (\cite[III.2]{MacLane1994}), is a cartesian 2-ring. We will see later that such examples are ill-suited for our purposes, and they need to be ruled out in order to prove our main result (\cref{thm.Gal}). 
\end{example}

	Now let $A$ be a commutative 2-ring, and denote by $\cat{Cog}(A)$ the category of cocommutative coalgebras in $A$. The comultiplication and counit of a coalgebra $a$ in a monoidal category will be denoted by $\Delta_a$ and $\varepsilon_a$, respectively. We will see in this subsection that $\cat{Cog}(A)$ is a cartesian 2-ring, and it is  a ``good replacement'' for $A$ from our point of view. First, we recall the following simple result, entirely analogous to the fact that the tensor product is the coproduct in the category of commutative algebras over some ground field. 
    
\begin{lemma}
	Let $A$ be a symmetric monoidal category, and $\cat{Cog}(A)$ the category of cocommutative coalgebras in $A$. Then the tensor product of coalgebras is the binary product in $\cat{Cog}(A)$.\qedhere 
\end{lemma}

    In particular, $\cat{Cog}(A)$ has a natural structure of cartesian monoidal category for any symmetric monoidal category $A$. We now want to show that when $A$ is a commutative 2-ring, $\cat{Cog}(A)$ itself is a cartesian 2-ring. The cocontinuity of functors of the form $x\otimes$ for $x\in \cat{Cog}(A)$ is clear, since the forgetful functor $\cat{Cog}(A)\to A$ preserves colimits (in other words, the colimits of cocommutative coalgebras in $A$ are the colimits of the diagrams of underlying objects). What is not obvious, however, is that $\cat{Cog}(A)$ is a presentable category. The result is an immediate consequence of Porst's work on (co)reflections for categories of (co)algebras in \cite{Porst2008}.
    
\begin{proposition}\label{prop.A'_cartesian}
    For a commutative 2-ring $A$, its category of cocommutative coalgebras $\cat{Cog}(A)$ is a cartesian 2-ring. Moreover, the forgetful functor $\cat{Cog}(A)\to A$ is comonadic.  
\end{proposition}     
\begin{proof}
    The previous discussion reduces the first statement to showing that $\cat{Cog}(A)$ is presentable. A commutative 2-ring is an example of what Porst calls an admissible monoidal category \cite[Definition 2.1]{Porst2008}. But then the theory developed in Porst's paper applies to $A$, and we conclude that $\cat{Cog}(A)$ is indeed presentable, and also comonadic over $A$ \cite[Summary 4.3]{Porst2008}.  
\end{proof}

    We now show that $\cat{Cog}(A)\to A$ is a universal map from a cartesian 2-ring to $A$:

\begin{proposition}\label{prop.cartify}
    Let $A,B$ be commutative 2-rings with $B$ cartesian. Then, the functor 
    \[
        U_*:\Hom_{\cat{Com2Ring}}(B,\cat{Cog}(A))\to\Hom_{\cat{Com2Ring}}(B,A)
    \] 
induced by the forgetful functor $U:\cat{Cog}(A)\to A$ is an equivalence of categories. 
\end{proposition}          
\begin{proof}
    Recall that in a cartesian monoidal category, all objects are coalgebras in a unique way, and all morphisms are automatically coalgebra maps. 

    Now let $b\in B$ be an object, and $\alpha:B\to A$ a morphism of commutative 2-rings. The object $\alpha(b)\in A$ has the usual coalgebra structure, given by the compositions
    \begin{equation}\label{eqn.coalg}
        \tikz[baseline=(current bounding box.center)]{
          \path (0,0) node[anchor=east] (alphab) {$\Delta: \alpha(b)$} +(2,0) node[anchor=west] (alphabb) {$\alpha(b\otimes b)$} +(5,0) node[anchor=west] (alphabalphab) {$\alpha(b)\otimes \alpha(b),$}
              ++(0,-1) node[anchor=east] (alphabbis) {$\varepsilon: \alpha(b)$} +(2,0) node[anchor=west] (alpha1) {$\alpha(1_B)$} +(5,0) node[anchor=west] (1) {$1_A$,};
          \draw[->] (alphab) -- (alphabb) node[pos=.5,auto] {$\scriptstyle \alpha(\Delta_b)$};
          \draw[->] (alphabb) -- (alphabalphab) node[pos=.5,auto] {$\scriptstyle \cong$};
          \draw[->] (alphabbis) -- (alpha1) node[pos=.5,auto] {$\scriptstyle \alpha(\varepsilon_b)$};
          \draw[->] (alpha1) -- (1) node[pos=.5,auto] {$\scriptstyle \cong$};
         }
     \end{equation}
where the right hand arrows are the ones giving the monoidal structure on $\alpha$. One checks immediately that if $f$ is a morphism in $B$, $\alpha(f)$ preserves these coalgebra structures.  

  The faithfulness of $U_*$ follows from that of $U$. 
  
  Now let $\alpha,\beta:B\to \cat{Cog}(A)$ be monoidal functors, and $\xi:U_*(\alpha)\to U_*(\beta)$ a monoidal natural transformation. Let $b\in B$ be an object, and consider the diagram
    \[\tikz{
  \path (0,0) node (alphab) {$\alpha(b)$} +(3,0) node (alphabb) {$\alpha(b\otimes b)$} +(6,0) node (alphabalphab) {$\alpha(b)\otimes \alpha(b)$}
              +(0,-2) node (betab) {$\beta(b)$} +(3,-2) node (betabb) {$\beta(b\otimes b)$} +(6,-2) node (betabbetab) {$\beta(b)\otimes \beta(b)$,};
  \draw[->] (alphab) -- (alphabb) node[pos=.5,auto] {$\scriptstyle \alpha(\Delta_b)$};
  \draw[->] (alphabb) -- (alphabalphab) node[pos=.5,auto] {$\scriptstyle \cong$};
  \draw[->] (betab) -- (betabb) node[pos=.5,auto] {$\scriptstyle \beta(\Delta_b)$};
  \draw[->] (betabb) -- (betabbetab) node[pos=.5,auto] {$\scriptstyle \cong$};  
  \draw[->] (alphab) -- (betab) node[pos=.5,auto,swap] {$\scriptstyle \xi_b$};
  \draw[->] (alphabb) -- (betabb) node[pos=.5,auto,swap] {$\scriptstyle \xi_{b\otimes b}$};
  \draw[->] (alphabalphab) -- (betabbetab) node[pos=.5,auto,swap] {$\scriptstyle \xi_b\otimes\xi_b$};
 }\]
where the rows are as in \cref{eqn.coalg}. The left hand square commutes because $\Delta_b$ is a map of coalgebras (since $b$ is a {\it cocommutative} coalgebra) and $\xi$ is a natural transformation, while the right hand square commutes because $\xi$ is monoidal. The commutativity of the outer rectangle means that $\xi_b$ preserves the comultiplications on $\alpha(b)$ and $\beta(b)$ constructed by \cref{eqn.coalg}. A similar argument shows that the counits on $\alpha(b)$ and $\beta(b)$ are also intertwined by $\xi_b$, and hence that $\xi_b$ is a map of coalgebras. This shows that the natural transformation $\xi$ is the image through $U_*$ of a natural transformation between $\alpha,\beta$, and hence proves the fullness of $U_*$.  
\end{proof}

    By applying \cref{prop.A'_cartesian} to the cartesian 2-ring $B={^G}\cat{Set}$, we are free to work entirely within the framework of cartesian 2-rings when studying $\pi_1$:

\begin{corollary}\label{cor.cartify}
    Let $A$ be a commutative 2-ring, and $U:\cat{Cog}(A)\to A$ the forgetful functor. Then, the natural transformation $\pi_1(\cat{Cog}(A))\to\pi_1(A)$ defined for each small groupoid $G$ by
    \[
        \pi_1(\cat{Cog}(A))(G)=\Hom_{\cat{Com2Ring}}({^G}\cat{Set},\cat{Cog}(A)) \verylongarrow{U_*} \Hom_{\cat{Com2Ring}}({^G}\cat{Set},A)=\pi_1(A)(G) 
    \]
gives an equivalence $\pi_1(\cat{Cog}(A))\simeq\pi_1(A)$. A similar statement holds for $\pi_1^f$.   \qedhere
\end{corollary}


\subsection{Smallness}

    In this short section we prove that for a commutative 2-ring $A$, the 2-functor $\pi_1(A)(-):\cat{Gpd}\to \cat{GPD}$ defined above (see \cref{defn.pi1} and \cref{remark.pi1GPD}) actually takes values in the 2-category of essentially small groupoids. This will allow us to interpret it as a endo-2-functor of $\cat{Gpd}$, which in turn will bring us one step closer to our goal of proving the pro-representability of $\pi_1(A)$ under suitable conditions on $A$. We first need some notations. 

    Let $A$ be a cartesian commutative 2-ring, $G$ a small groupoid, and $y\in A$ an object. We denote by $G_y$ the groupoid with the same objects as $G$, and whose arrows $s\to t$ (for $s,t\in G_0$) are the elements of $\Hom_A(y,G(t,s)\cdot 1_A)$, with the obvious composition given by composition of arrows in $G$.   

\begin{lemma}\label{lemma.smallness}
    Let $A$ be a cartesian commutative 2-ring, $G$ a (small) groupoid, and $\hat x$ a right $G$-pseudotorsor in $A$. Then, for any object $y\in A$, the functor $G^{op}\to\cat{Set}$ defined by $G_0\ni s\mapsto\Hom_A(y,\hat x(s))$ defines a right $G_y$-pseudotorsor in $\cat{Set}$. 
\end{lemma}
\begin{proof}    
    Unpacking the definition of $\tau_x$ (see \cref{defn.torsor}), we get an isomorphism
\[
    x(t)\times(G(t,s)\cdot 1_A)\stackrel{\cong}{\longrightarrow} x(t)\times x(s)
\]
for every pair $(t,s)\in \pi_0(G)$. Now simply apply the functor $\Hom_A(y,-)$. 
\end{proof}

\begin{proposition}\label{prop.smallness}
	Let $A$ be a commutative 2-ring, and $G$ a small groupoid. Then there are setly many $G$-torsors in $A$.  
\end{proposition}

\begin{proof}
    We assume that $A$ is cartesian (we can always do this, according to \cref{cor.cartify}). Suppose $A$ is $\kappa$-presentable for some regular cardinal $\kappa$. Let $\hat x$ be a $G$-torsor in $A$, and $y\in A$ a $\kappa$-presentable object. 
    
    We know from \cref{lemma.smallness} that $s\mapsto\Hom_A(y,\hat x(s))$ is a right $G_y$-pseudotorsor in $\cat{Set}$, which we denote by $\tilde y$. The pseudotorsor condition (see \cref{defn.torsor}) now implies that for $s\in G_0$, $\tilde y(s)$ is either empty or a principal homogeneous space over the automorphism group of $s$ in the groupoid $G_y$. Since there are setly many $\kappa$-little objects $y$, there is a uniform bound on the sizes of all $G_y$. This implies that there is a uniform bound on the sets $\Hom_A(y,\hat x(s))$, and hence that the canonical cocone (\cite[0.4]{MR1294136}) on the $\kappa$-little objects with tip $\hat x(s)$ is $\lambda$-small for some cardinal $\lambda >\kappa$ which does not depend on $\hat x$. Since $\hat x(s)$ is the colimit of the canonical cocone (see \cite[Proposition 1.22]{MR1294136}) and the colimit of a $\lambda$-small diagram of $\kappa$-little objects is $\lambda$-little (\cite[Proposition 1.16]{MR1294136}), $\hat x(s)$ must be $\lambda$-little. But this means that there is a small set of choices of $\hat x(s)$ (up to isomorphisms), and the conclusion now follows easily.  
\end{proof}


\Section{The main result} \label{section.mainresults}

\subsection{2-limits and pro-representability} \label{section.2limits}

    The previous section shows that our (2-)functor $\pi_1:\cat{Gpd}\to\cat{CAT}$ can actually be interpreted as a 2-functor from the 2-category $\cat{Gpd}$ of small groupoids into itself. As indicated before, our final aim will be to show that this 2-functor and its cousin $\pi_1^f$ are actually pro-representable. To make the categorical terminology precise, we will refer mainly to \cite{Lurie2009}. 
    
    Lurie works with $(\infty,1)$-categories, which are simplicial sets with a certain property. The notion can be specialized to that of $(n,1)$-category (\cite[2.3.4]{Lurie2009}). For us, the theory is important through its applications to $(2,1)$-categories. One can define $(2,1)$-categories in a ``classical'' way too, without resorting to the theory of simplicial sets, by starting out with the usual notion of (weak) 2-category, also called a bicategory (as in \cite{Benabou1967}, for example), and then imposing the condition that 2-morphisms be invertible. We take it for granted that the two notions are ``equivalent'' in a more or less straightforward way, and we will henceforth apply the results of \cite{Lurie2009} directly to the $(2,1)$ level without making all the necessary translations explicitly. 
    
    Generalities on 2-categories, 2-functors, etc.\ can, by now, be found in countless sources; see \cite[1.5]{Leinster2004}, for instance, for a very concise and readable account. In the case of (2,1)-categories, limits in the sense of \cite[Chapter 4]{Lurie2009} are simply those of the underlying 2-categories, sometimes called bilimits, treated for example in \cite[Section 6]{Kelly1989} or \cite[Section 1]{Street1980}. When we happen to refer to ``limits (colimits) in a 2-category,'' bilimits (bicolimits) are always what we mean. We will describe the specific limits used in our arguments, such as equalizers or equifiers of groupoids. 
    
\begin{definition}
    A poset $J$ is said to be \define{left-filtered} (or, equivalently for posets, \define{left-directed}) if for any pair of elements $x,y\in J$ there is an element $z$ with $z\le x,y$. One similarly defines \define{right-filtered}/\define{right-directed} posets by reversing the inequalities. 
\end{definition}
\begin{definition}
    Let $C$ be a (2,1)-category. A \define{pro-object in $C$} is a 2-functor $x:J\to C$ for some left-filtered poset $J$, the latter being regarded as a (2,1) category in an obvious way.  
    
The pro-objects of $C$ form a (2,1)-category $\cat{Pro-}C$ with hom-categories defined as follows: if $x:J\to C$ and $y:I\to C$ are two pro-objects, set
    \[
        \Hom_{\cat{Pro-}C}(x,y) = \varprojlim_{i\in I}\varinjlim_{j\in J}\Hom_C(x(j),y(i)).
    \] 
\end{definition}    
\begin{remark}
    A (2,1)-category $C$ has a fully faithful embedding in $\cat{Pro-}C$: simply send an object $c\in C$ to the constant diagram $c$.
\end{remark}
\begin{definition}
    If $C$ is a (2,1)-category, a functor $F:C\to\cat{Gpd}$ is \define{pro-representable} if it is (equivalent to a functor) of the form $\Hom_{\cat{Pro-}C}(x,-):C\to\cat{Gpd}$ for some pro-object $x\in\cat{Pro-}C$. 
\end{definition}

    Finally, we adopt some non-standard terminology regarding finite limits:

\begin{definition}\label{defn.finitely}
    A (2,1)-category is  \define{finitely complete} if it admits finite products and equalizers. A functor between (2,1)-categories is \define{finitely continuous} if it preserves finite products and equalizers. 
\end{definition}

We then have the following immediate consequence of the results in \cite{Lurie2009}:

\begin{proposition}\label{prop.prorep}
    Let $C$ be a finitely complete (2,1)-category. A functor $F:C\to\cat{Gpd}$ is pro-representable if and only if it is finitely continuous.
\end{proposition}
\begin{proof}
    This is the analogue for (2,1)-categories of the dual to \cite[5.3.5.4]{Lurie2009}, coupled with the fact that Lurie's finite limits can be built up from finite products and equalizers (see \cref{remark.prodseqs} below). 
\end{proof}

\begin{remark}
    We could have defined pro-objects by means of left-filtered categories, or even left-filtered (2,1)-categories. However, since we are essentially interested only in taking limits along left-filtered diagrams, the argument dual to \cite[5.3.1.16]{Lurie2009} shows that there is no loss of generality in using only left-filtered posets. 
\end{remark}    
\begin{remark}\label{remark.prodseqs}
    One of the advantages of working with groupoids rather than arbitrary categories is that $\cat{Gpd}$ is a (2,1)-category rather than merely a 2-category; we can thus make full use of the fact that, just as for an ordinary category, (finite) limits in a (2,1)-category are constructed from (finite) products and equalizers (dual version of \cite[4.4.3.2]{Lurie2009}). Compare this to the more general setting of 2-categories, where first of all one has to bring in so-called weighted or indexed limits (\cite[1.12]{Street1980}), and secondly, more work is needed to check the existence and preservation of (weighted) limits: one can adapt the argument in \cite[Proposition 4.4]{Kelly1989} to show that products, equalizers and equifiers will do; alternatively, products, equalizers and weighted powers suffice, as in \cite[1.24]{Street1980}, corrected in \cite{Street1987} (Street refers to weighted powers as ``cotensor biproducts''). Since both \cat{Gpd} and $\cat{Gpd}_f$ are finitely complete (see the examples below), the pro-representability of $\pi_1$ and $\pi_1^f$ is reduced to checking that the two functors  preserve finite products and equalizers.  
\end{remark}

	We now recall the examples of limits in $\cat{Gpd}$ and $\cat{Gpd}_f$ which will be relevant to the discussion below.

\begin{example}[Products]
    Let $G_i$, $i\in I$ be a set of (small) groupoids. The \define{product} $\displaystyle G=\prod_IG_i$ is the groupoid, uniquely defined up to equivalence, with maps $\pi_i:G\to G_i$, universal in the sense that they induce a natural equivalence of groupoids
    \[
        \Hom_{\cat{Gpd}}(-,G) = \prod_I\Hom_{\cat{Gpd}}(-,G_i). 
    \]
It is easy to check that (up to equivalence) this is simply the familiar product: $G_0$ is the usual cartesian product of sets $\displaystyle \prod_I(G_i)_0$, and similarly for arrows. When the indexing family $I$ happens to be empty, one gets \define{the terminal object} of $\cat{Gpd}$: it is just the groupoid $1$, with one object and one morphism. $\cat{Gpd}$ has products indexed by any set $I$, while $\cat{Gpd}_f$ has finite products (i.e.\ those indexed by finite $I$). 
\end{example}
\begin{example}[Equalizers]\label{eg.equalizers}
    These are the limits (in the 2-categorical sense, as always) of diagrams of the form
    \[\tikz[anchor=base]{
      \path (0,0) node (G) {$G$} +(2,0) node (H) {$H$};
      \draw[transform canvas={yshift=0.8mm},->] (G) -- (H) node[pos=.5,auto] {$\scriptstyle \varphi$};
      \draw[transform canvas={yshift=-0.8mm},->] (G) -- (H) node[pos=.5,auto,swap] {$\scriptstyle \psi$}; 
 }\]
in $\cat{Gpd}$ or $\cat{Gpd}_f$. 

The strict analogue of this kind of limit, when one works in a strict 2-category and the limits are required to be defined up to \define{isomorphism} rather than only up to equivalence, would be referred to as an \define{iso-inserter} (\cite[Section 4]{Kelly1989}). For us, the limit would have to be a groupoid $K$ with a functor $\xi:K\to G$ and a natural transformation $p:\varphi\circ\xi\Rightarrow \psi\circ\xi$ which is universal in the appropriate sense. We describe the construction, since it will be of central importance below. 
    
	The objects of $K$ are pairs $(s,h)$, where $s\in G_0$ and $h:\varphi(s)\to\psi(s)$ is an arrow in $H$. Given such pairs $(s,h)$ and $(s',h')$, an arrow $(s,h)\to(s',h')$ in $K$ is an arrow $g$ of $G$ for which the diagram
    \[\tikz[anchor=base]{
         \path (0,0) node (phis) {$\varphi(s)$} +(2,0) node (phis') {$\varphi(s')$} +(0,-2) node (psis) {$\psi(s)$} +(2,-2) node (psis') {$\psi(s')$};
      \draw[->] (phis) -- (phis') node[pos=.5,auto] {$\scriptstyle \varphi(g)$};
      \draw[->] (psis) -- (psis') node[pos=.5,auto] {$\scriptstyle \psi(g)$};
      \draw[->] (phis) -- (psis) node[pos=.5,auto,swap] {$\scriptstyle h$}; 
      \draw[->] (phis') -- (psis') node[pos=.5,auto] {$\scriptstyle h'$};
 }\]
commutes. The natural transformation $p$ is defined by sending the object $(s,h)\in K_0$ to the arrow $h\in H$. 
\end{example}
\begin{example}[Equifiers]
	Although, according to \cref{prop.prorep}, in order to prove the pro-representability of a functor between (2,1)-categories all one needs to do is check preservation of finite products and equalizers, it is nevertheless instructive to study another kind of 2-limit, namely equifiers. In fact, the failure of equifier preservation, which is very easy to verify, will prompt us to impose additional conditions on our commutative 2-ring $A$ in order to prove the pro-representability of $\pi_1(A)$ in \cref{thm.Gal} (\cref{defn.good2rings}).  
    
	An equifier in a (2,1)-category is simply a limit of a diagram of the form
    \[ \tikz[anchor=base]{
    	\path (0,0) node (G) {$G$} +(1,0) node {$\Downarrow$ $\scriptstyle \xi$} +(2,0) node {$\Downarrow$ $\scriptstyle \xi'$}  +(3,0) node (H) {$H$};
    	\draw[->] (G) .. controls +(1,1) and +(-1,1) .. (H) node[pos=.5,auto] {$\scriptstyle \varphi$};
    	\draw[->] (G) .. controls +(1,-1) and +(-1,-1) .. (H) node[pos=.5,auto,swap] {$\scriptstyle \psi$};
  } \]
It consists of an arrow $j:K\to G$ such that $\xi.j=\xi'.j$, which is universal in the appropriate sense (where the lower dot denotes ``horizontal'' composition of 1- and 2-morphisms). If the diagram above is in $\cat{Gpd}$ (or $\cat{Gpd}_f$), $j:K\to G$ is nothing but the inclusion of the full subgroupoid of $G$ generated by those objects on which the two natural transformations $\xi,\xi'$ are equal.      

    Note that for this description of equifiers, it is crucial that we work with (2,1)-categories rather than just 2-categories: indeed, equifiers are usually described as weighted limits (cf.\ \cref{remark.prodseqs}), and in the more general setting, the usual definition need not coincide with ours. 
\end{example}
\begin{remark}
    Equifiers are finite limits in the sense of \cite{Lurie2009}: when we view a (2,1)-category $C$ as a simplicial set, an equifier is a limit of a diagram $K\to C$ for a finite simplicial set $K$. It follows from the dual of \cite[4.4.3.2]{Lurie2009} that they can be built up from products and equalizers, and  hence are preserved by any finitely continuous functor as in \cref{defn.finitely}. 
\end{remark}

	We make one final observation, this time on colimits in the 2-category $\cat{CAT}$ of possibly large but locally small categories, which will be useful later. Let $J$ be a right-directed poset, and $x:J\to\cat{CAT}$ a \define{strict} 2-functor; in other words, $x$ is a 1-functor from $J$ to $\cat{CAT}$, when both of these are regarded as ordinary 1-categories. Then $x$ has a 2-colimit $\varinjlim x$, and also a colimit $\varinjlim' x$ when regarded as a 1-functor. The universal property of $\varinjlim$ then gives us a functor $\varinjlim x\to \varinjlim' x$. Our observation, whose proof we leave to the reader, is

\begin{proposition}\label{prop.strictcolim}
    In the setting outlined above, the canonical functor $\varinjlim x\to \varinjlim' x$ is an equivalence.\qedhere
\end{proposition}


\subsection{\texorpdfstring{Does $\pi_1(A)$ preserve limits?}{Does pi1(A) preserve limits?}}

	We now come to the main question that we have to answer in order to prove that $\pi_1(A)$ and $\pi_1(A)_f$ are pro-representable (see \cref{prop.prorep}): are they finitely continuous? This question can also be phrased differently, as explained below. 

	Let $G$ be a small category.  As we remarked in \cref{eg.2rings}, $^G\cat{Set} = \Fun(G, \cat{Set}) = \sh(G^\op)$ is a commutative cartesian 2-ring.  Given a functor $f: G \to H$, the pullback map $f^*: {^G\cat{Set}} \from {^H\cat{Set}}$ has both right- and left- adjoints $f_*$ and $f_!$, given by the right and left Kan extensions.  Hence $f^*$ is both cocontinuous and continuous, and in particular respects finite products.  Thus it is a morphism of commutative 2-rings, and similarly natural transformations of functors $G \to H$ induce 2-morphisms of 2-rings.  Therefore the construction $G \mapsto \sch G$ extends naturally to a 2-functor $\sch\Box: \cat{Cat} \to \cat{2AfSch}$, where $\cat{Cat}$ is the 2-category of small categories, functors, and natural transformations. We will soon specialize to $\sch G$ for small \define{groupoids} $G$, but we denote the restriction of $\sch\Box$ to $\cat{Gpd}$ by the same symbol. 

	Since $\pi_1(A)$ is the composition
\begin{equation}\label{eq.scheming}
	G\mapsto \sch G 
  {\!{\tikz[anchor=base,auto,baseline=(a.base)] \draw[->] node (a) {\phantom{A}} ++(3cm,0) node (b) {\phantom{A}\!\!\!} (a.mid) -- node{$\scriptstyle \Hom_{\cat{2AfSch}}(\Spec(A),-)$} (b.mid);}}
	\Hom_{\cat{2AfSch}}(\Spec(A),\sch G),    
\end{equation}
it is finitely continuous as soon as $\sch\Box$ is. We begin with a result which applies to arbitrary categories, so we phrase it using the $\cat{Cat}$ version of $\sch\Box$.  
   
\begin{proposition}\label{prop.products}
  The functor $\sch\Box: \cat{Cat} \to \cat{2AfSch}$ preserves finite products.
\end{proposition}
\begin{proof}
  This follows from the equivalence $^{G\times H}\cat{Set}={^G}\cat{Set}\boxtimes{^H}\cat{Set}$ of commutative 2-rings (noted in \cref{eg.GH_boxtimes} at the level of abelian 2-groups) and \cref{remark.boxtimes_is_coprod}, saying that $\boxtimes$ is the binary coproduct of commutative 2-rings.  
\end{proof}

	For the remaineder of the paper, we specialize from \cat{Cat} to \cat{Gpd}.  Nevertheless:
\begin{lemma}
  The functor $\sch\Box: \cat{Gpd} \to \cat{2AfSch}$ does not preserve infinite products.
\end{lemma}

	As we use the word \define{continuous} to describe functors that preserve categorical limits, we prefer not to use the word ``continuous homomorphism'' when we mean that it respects some extra topology.  Rather, given topological categories $X,Y$, we say that a functor $f: X\to Y$ is \define{topological} if it is ``continuous'' in the non-categorical sense.

\begin{proof}
  Let $\Fp$ is the field with $p$ elements and $G$ a groupoid in \cat{Set}.  We will show in \cref{prop.pi1_of_fields} that $\Hom_{\cat{2AfSch}}(\Spec({^{\Fp}\cat{Vect}}), \sch G)$ is naturally equivalent to the groupoid whose objects are topological homomorphisms $\hat \ZZ \to G$, where $\hat\ZZ = \Gal\left(\overline{\Fp} / \Fp\right)$ is given its profinite topology and $G$ has the discrete topology.  (The morphisms are natural transformations of such functors; since $\hat\ZZ$ has only one object, one need not think about topology when defining the natural transformations.)  In particular, $\Hom_{\cat{2AfSch}}(\Spec({^{\Fp}\cat{Vect}}), \sch \Box) = \Hom_{\cat{TopGpd}}(\hat\ZZ, \Box)$ respects all finite limits of groupoids.
  
  But now set $G = \hat\ZZ$ with its discrete topology.  If $\Hom_{\cat{2AfSch}}(\Spec({^{\Fp}\cat{Vect}}), \sch \Box)$ respected all limits, then $\id: \hat\ZZ \to \hat\ZZ$ would be in $\Hom_{\cat{2AfSch}}(\Spec({^{\Fp}\cat{Vect}}), \sch {\hat\ZZ})$, but $\id: \text{profinite} \to \text{discrete}$ is not topological.  Thus the functor $\Hom_{\cat{2AfSch}}(\Spec({^{\Fp}\cat{Vect}}), \sch \Box)$ cannot respect infinite products, and hence neither can $\sch\Box$.
\end{proof}

\begin{remark}  
  The problem is that products of categories correspond to tensor products.  But we do not demand of a commutative 2-ring that it be able to make sense of infinite tensor products.  The reader is invited to consider also the case of an infinite product of the groupoid with two objects (and only identity morphisms).
\end{remark}

	Regarding limits other than products, we have the following negative result:

\begin{lemma}\label{eg.equifiers}
  The functor $\sch\Box: \cat{Gpd} \to \cat{2AfSch}$ does not preserve equifiers.
\end{lemma}
\begin{proof}
  Let $G,H$ be finite groups.  Pick a homomorphism $f: G \to H$ and a non-identity element $h\in H$ that centralizes the subgroup $f(G) \subseteq H$.  Then $h$ determines a nontrivial natural automorphism of the functor $f$, where $G,H$ are thought of as one-object groupoids.  The equifier of the diagram
  \[ \tikz[anchor=base]{
    \path (0,0) node (G) {$G$} +(1,0) node {$\Downarrow$ $\scriptstyle \id$} +(2,0) node {$\Downarrow$ $\scriptstyle h$}  +(3,0) node (H) {$H$};
    \draw[->] (G) .. controls +(1,1) and +(-1,1) .. (H) node[pos=.5,auto] {$\scriptstyle f$};
    \draw[->] (G) .. controls +(1,-1) and +(-1,-1) .. (H) node[pos=.5,auto,swap] {$\scriptstyle f$};
  } \]
  is the empty groupoid $\emptyset$:
  \[ \tikz[anchor=base]{
    \path (0,0) node (G) {$G$} +(-2,0) node (O) {$\emptyset$} +(1,0) node {$\Downarrow$ $\scriptstyle \id$} +(2,0) node {$\Downarrow$ $\scriptstyle h$}  +(3,0) node (H) {$H$};
    \draw[->] (G) .. controls +(1,1) and +(-1,1) .. (H) node[pos=.5,auto] {$\scriptstyle f$};
    \draw[->] (G) .. controls +(1,-1) and +(-1,-1) .. (H) node[pos=.5,auto,swap] {$\scriptstyle f$};
    \draw[->] (O) -- (G);
  } \]
  
  However, at the level of commutative 2-rings, the diagram
  \begin{equation} \label{eqn.equifier} \tikz[anchor=base,baseline=(G.base)]{
    \path (0,0) node (G) {$^G\cat{Set}$} +(1,0) node {$\Downarrow$ $\scriptstyle \id$} +(2,0) node {$\Downarrow$ $\scriptstyle h$}  +(3,0) node (H) {$^H\cat{Set}$};
    \draw[<-] (G) .. controls +(1,1) and +(-1,1) .. (H) node[pos=.5,auto] {$\scriptstyle f^*$};
    \draw[<-] (G) .. controls +(1,-1) and +(-1,-1) .. (H) node[pos=.5,auto,swap] {$\scriptstyle f^*$};
  } \end{equation}
  does not have $^\emptyset\cat{Set} = 0$ as its coequifier.  Indeed, consider the 2-ring $\{0 \to 1\}$ with two objects, a unique nontrivial morphism, and $\otimes = \times$.  There is a unique 1-morphism of 2-rings $^G\cat{Set} \to \{0 \to 1\}$, given by sending all non-initial objects of $^G\cat{Set}$ to the object $1\in \{0\to 1\}$.  This homomorphism does not factor through the zero 2-ring, but upon whiskering by this 1-morphism, the two 2-morphisms $\id,h$ of 1-morphisms $^H\cat{Set} \to \{0 \to 1\}$ become equal.  So the coequifier of the diagram \cref{eqn.equifier}, if it exists, has a 1-morphism to $\{0\to 1\}$, and thus cannot be the zero 2-ring.
\end{proof}

	As explained in the proof of \cref{eg.equifiers}, the problem for equifiers is that posets like $\{0\to 1\}$ are bad. In order to rule out such examples, we will restrict our attention to a certain full sub-2-category of $\cat{Com2Ring}$.

\begin{definition}\label{defn.goodness_cond}
	Let $A$ be a cartesian 2-ring. We say that $A$ has \define{disjoint coproducts} if whenever we have an expression $\displaystyle a=\coprod_{i\in I}a_i$ in $A$ of an object $a$ as a coproduct, the canonical maps $a_i\to a$ are monomorphisms, and the pullbacks $a_i\times_a a_j$ are the initial object for all $i\ne j\in I$. 

	We say that coproducts in $A$ are \define{stable} if for any map $\displaystyle b\to a=\coprod_{i \in I}a_i$ in $A$, the canonical maps $b_i=b\times_a a_i\to b$ make $b$ the coproduct $\displaystyle\coprod_{i\in I}b_i$. 
\end{definition}

\begin{definition}\label{defn.good2rings}
	A cartesian 2-ring is said to be \define{good} if its coproducts are disjoint and stable. A commutative 2-ring $A$ is good if $\cat{Cog}(A)$ is a good cartesian 2-ring (see \cref{prop.A'_cartesian}).     
  
  We denote by $\cat{Good2Rings}$ the 2-category of good commutative 2-rings with the usual 1- and 2-morphisms, and by $\cat{Good2AfSch}$ its 1-opposite. $\cat{GoodCart}$ will be the category of good cartesian 2-rings.  
\end{definition}

\begin{example}
	Notice that the notions of disjoint and stable coproducts are precisely the ones used in the topos literatue (e.g. \cite[Section 1 of Appendix]{MacLane1994}): all Grothendieck topoi are (cartesian) good 2-rings. In particular, $^G\cat{Set}$ is a good 2-ring for any small category $G$. We will henceforth regard the functor $\sch\Box:\cat{Gpd}\to\cat{2AfSch}$ as taking values in $\cat{Good2AfSch}$.  
\end{example}
\begin{example}
	Although we will not prove this here, every abelian commutative 2-ring (i.e.\ commutative 2-ring which is also an abelian category) is good. In particular, we can apply any results we obtain for good 2-rings to, for example, the category of quasicoherent sheaves on a scheme. 
\end{example}

	The merit of \cref{defn.good2rings} is that it allows us to prove the following result: 

\begin{proposition}\label{prop.equalizers}
	The functor $\sch\Box:\cat{Gpd}\to\cat{Good2AfSch}$ preserves equalizers. 
\end{proposition}

Before going into the proof, we rephrase the problem slightly. Let
	\begin{equation}\label{eq.equalizer_diag}
  	\tikz[baseline=(current  bounding  box.center)]{
    	\path (0,0) node (K) {$K$} +(1.5,.5) node (G) {$G$} +(1.5,-.5) node (Gbis) {$G$} +(1.5,0) node {$\Downarrow$ $\scriptstyle p$} +(3,0) node (H) {$H$};
      \draw[->] (K) -- (G) node[pos=.7,auto] {$\scriptstyle \xi$};
      \draw[->] (K) -- (Gbis) node[pos=.7,auto,swap] {$\scriptstyle \xi$};
      \draw[->] (G) -- (H) node[pos=.3,auto] {$\scriptstyle \varphi$};
      \draw[->] (Gbis) -- (H) node[pos=.3,auto,swap] {$\scriptstyle \psi$};
    } 
  \end{equation}
be an equalizer diagram in $\cat{Gpd}$ as described in \cref{eg.equalizers}, and 
	\[\tikz[anchor=base]{
  	\path (0,0) node (A) {$A$} +(2,.5) node (G) {$^G\cat{Set}$} +(2,-.5) node (Gbis) {$^G\cat{Set}$} +(2,0) node {$\Downarrow$ $\scriptstyle q$} +(4,0) node (H) {$^H\cat{Set}$};
    \draw[<-] (A) -- (G) node[pos=.7,auto] {$\scriptstyle \alpha$};
    \draw[<-] (A) -- (Gbis) node[pos=.7,auto,swap] {$\scriptstyle \alpha$};
    \draw[<-] (G) -- (H) node[pos=.3,auto] {$\scriptstyle \varphi^*$};
    \draw[<-] (Gbis) -- (H) node[pos=.3,auto,swap] {$\scriptstyle \psi^*$};
  } \]
a diagram in $\cat{Good2Rings}$. We have to prove that the latter factors as   
    \[\tikz[anchor=base]{
    \path (0,0) node (K) {$^K\cat{Set}$} +(2,.5) node (G) {$^G\cat{Set}$} +(2,-.5) node (Gbis) {$^G\cat{Set}$} +(2,0) node {$\Downarrow$ $\scriptstyle p^*$} +(4,0) node (H) {$^H\cat{Set}$} +(-2,0) node (A) {$A$};
    \draw[<-] (K) -- (G) node[pos=.8,auto] {$\scriptstyle \xi^*$};
    \draw[<-] (K) -- (Gbis) node[pos=.8,auto,swap] {$\scriptstyle \xi^*$};
    \draw[<-] (G) -- (H) node[pos=.3,auto] {$\scriptstyle \varphi^*$};
    \draw[<-] (Gbis) -- (H) node[pos=.3,auto,swap] {$\scriptstyle \psi^*$};
    \draw[<-] (A) -- (A-|K.west) node[pos=.5,auto] {$\scriptstyle \beta$};
  } \]
up to the obvious notion of natural equivalence, and that this factorization is unique, again up to the appropriate notion of equivalence. 
 
	First, according to \cref{cor.cartify}, we may as well assume that $A$ is cartesian, and everything takes place inside $\cat{GoodCart}$. Secondly, by \cref{prop.GSet -> A are torsors}, we can recast everything in terms of torsors in $A$. In that language, $\alpha$ has a corresponding right $G$-torsor $\hat x=\hat x(\alpha)$ in $A$, while the natural isomorphism $q$ is an isomorphism $q:\hat x\boxtimes_\varphi H\to \hat x\boxtimes_\psi H$ of right $H$-torsors in $A$. We now elaborate on this notation. 

	The groupoid morphism $\varphi:G\to H$ allows us to construct the $H$-torsor $_\varphi H=\varphi^*(H^H)$ in $^G\cat{Set}$. It now makes sense to talk about the right $H$-representation $\hat x\boxtimes_G(_\varphi H)$, which we denote by $\hat x\boxtimes_\varphi H$ to keep things simple. It can now be checked that $\hat x\boxtimes_\varphi H$ is actually a right $H$-torsor in $A$. 

	Regarding $\boxtimes_\varphi H$ and $\boxtimes_\psi H$ as functors from $\cat{Tors}(G,A)$ to $\cat{Tors}(H,A)$, let $\cat{Tors}(G,A;\varphi,\psi)$ be the equalizer of $\boxtimes_\varphi H$, $\boxtimes_\psi H:\cat{Tors}(G,A)\to \cat{Tors}(H,A)$. Its objects are right $G$-torsors $\hat x$ in $A$ together with a morphism $q:\hat x \boxtimes_\varphi H \to \hat x \boxtimes_\psi H$ (which will then automatically be an isomorphism by \cref{thm.iso}), and a map from $(\hat x,q)$ to $(\hat x',q')$ is a morphism $f:\hat x\to \hat x'$ in $\cat{Tors}(G,A)$ commuting with the $q$'s in the obvious sense. 

	For a right $K$-torsor $\hat y$ in $A$, one can similarly talk about the right $G$-torsor $\hat y\boxtimes_\xi G$. Note that there is an isomorphism $(\hat y\boxtimes_\xi G)\boxtimes_\varphi H\cong \hat y\boxtimes_{\varphi\circ\xi}H$ (and similarly for $\psi$). Furthermore, since by the universal property of the equalizer $\xi:K\to G$ we have an isomorphism $_{\varphi\circ\xi}H\cong {_{\psi\circ\xi}H}$, there is an isomorphism $\hat y\boxtimes_{\varphi\circ\xi}H\cong \hat y\boxtimes_{\psi\circ\xi}H$, which we will denote by $p_{\hat y}$. There is a functor $T:\cat{Tors}(K,A)\to\cat{Tors}(G,A;\varphi,\psi)$ defined by sending the right $K$-torsor $\hat y$ to $\hat y\boxtimes_\xi G$, together with the isomorphism 
	\[
 		(\hat y\boxtimes_\xi G)\boxtimes_\varphi H\cong \hat y\boxtimes_{\varphi\circ\xi}H\verylongarrow{\id_{\hat y}\boxtimes p} \hat y\boxtimes_{\psi\circ\xi}H \cong (\hat y\boxtimes_\xi G)\boxtimes_\psi H.
  \]
  
\Cref{prop.equalizers} can now be rephrased as 

\begin{proposition}\label{prop.equalizers bis}
    Given an equalizer diagram \cref{eq.equalizer_diag} and a good cartesian 2-ring $A$, the functor \[T:\cat{Tors}(K,A)\to\cat{Tors}(G,A;\varphi,\psi)\] defined above is an equivalence. 
\end{proposition}

	We work our way towards the proof by first showing:

\begin{lemma}
	\Cref{prop.equalizers bis} holds if and only if it holds in the case when $G$ is a connected groupoid. 
\end{lemma}
\begin{proof}
	We only prove the non-obvious implication connected $\Rightarrow$ general. Let $G_i$, $i\in I$ be the connected components of $G$. Then the equalizer $\xi:K\to G$ breaks up as a disjoint union of $K_i=\xi^{-1}(G_i)$, which are the equalizers of the restrictions of $\varphi$ and $\psi$ to the $G_i$. 
	
	Now let $\mathcal E=\{e_i\ |\ i\in I\}$ be a \define{partition of unity} of $A$, i.e. a decomposition $1_A=\coprod e_i$ for objects which are idempotent via the diagonal morphism: $e_i\simeq e_i\times e_i$. Define $\cat{Tors}_{\mathcal E}(K,A)$ to be the category of 2-ring morphisms (\cref{prop.GSet -> A are torsors}) $\alpha:{^K\cat{Set}}\to A$ which send $^KK_i$ to $A_i=e_iA=$ the subcategory of $A$ consisting of objects of the form $e_i\times a$. Similarly, let $\cat{Tors}_{\mathcal E}(G,A;\varphi,\psi)$ be the subcategory of $\cat{Tors}(G,A;\varphi,\psi)$ consisting of those objects whose underlying 2-ring morphism $^G\cat{Set}\to A$ lands $^GG_i$ in $A_i$. 
	
	It now remains to observe that (a) the categories $\cat{Tors}(K,A)$ and $\cat{Tors}(G,A;\varphi,\psi)$ break up as the possibly large coproducts (i.e. disjoint unions)
	\[
		\cat{Tors}(K,A)=\coprod_{\mathcal E}\cat{Tors}_{\mathcal E}(K,A),\qquad \cat{Tors}(G,A;\varphi,\psi)=\coprod_{\mathcal E}\cat{Tors}_{\mathcal E}(G,A;\varphi,\psi) 
	\] 
(for $\alpha:{^K\cat{Set}\to A}$, for example, let $e_i$ be the colimit of the $K^\op$-diagram $\alpha({^KK_i^K})$), (b) the $A_i$ are good cartesian 2-rings in their own right with units $e_i$, and (c) the functor $T$ from the statement of \cref{prop.equalizers bis} breaks up as a coproduct of the corresponding functors 
	\[
		T_{\mathcal E}:\cat{Tors}(K_i,A_i)\to\cat{Tors}_{\mathcal E}(G_i,A_i;\varphi,\psi),	
	\]
(with $\alpha({^KK_i^K})$ being a $K_i$-torsor in $A_i$ for any $\alpha\in\cat{Tors}(K,A)$, and so on). All of these remarks are more or less routine.  
\end{proof}

	In view of this result, for the remainder of this section up to the statement of \cref{prop.limits} we will assume that $G$ is a group. 
    
	Now recall from \cref{lemma.Galoiscond} that, denoting by $x$ the underlying object of a $G$-torsor $\hat x$, we have an isomorphism $\tau_x:x\boxtimes G\to x\boxtimes \hat x$ in $A^G$ ($G_0$ is a singleton now, so we no longer need to worry about $\boxtimes_{G_0}$, etc.). If we have a map $q:\hat x\boxtimes_\varphi H\to \hat x\boxtimes_\psi H$, we can apply the functor $x\boxtimes:A^H\to A^H$ to get the left hand vertical arrow in 
	\begin{equation}\label{eq.Hx}
  	\tikz[baseline=(current  bounding  box.center)]{
    	\path (0,0) node (1) {$x\boxtimes\hat x\boxtimes_\varphi H$} +(5,0) node (2) {$x\boxtimes G\boxtimes_\varphi H$} +(0,-2) node (3) {$x\boxtimes\hat x\boxtimes_\psi H$} +(5,-2) node (4) {$x\boxtimes_ G\boxtimes_\psi H$,};    
      \draw[->] (1) -- (2) node[pos=.5,auto] {$\scriptstyle \tau_x\boxtimes_\varphi H$};
      \draw[->] (3) -- (4) node[pos=.5,auto] {$\scriptstyle \tau_x\boxtimes_\psi H$};
      \draw[->] (1) -- (3) node[pos=.5,auto,swap] {$\scriptstyle x\boxtimes q$};
      \draw[->] (2) -- (4) node[pos=.5,auto,swap] {$\scriptstyle \cong$};
    } 
	\end{equation}
where the other vertical arrow is chosen so as to make the diagram commutative. In the sequel, we will abuse notation and write this vertical arrow as $x\boxtimes q:x\boxtimes_\varphi H\to x\boxtimes_\psi H$, where the latter are just $x\boxtimes({_\varphi H})$ and $x\boxtimes({_\psi H})$, and are isomorphic to the upper and lower right corners of the square, respectively.

\begin{proof of equalizers bis}
	We actually  a pseudo-inverse $S:\cat{Tors}(G,A;\varphi,\psi)\to \cat{Tors}(K,A)$ in several stages, using the notations from \cref{eg.equalizers} throughout the proof.  

\begin{defnS}
  We only describe the underlying object function of $S$, since the behavior of $S$ on arrows and the functoriality will then be clear. Let $\hat x$ be a right $G$-torsor in $A$, and $q:\hat x\boxtimes_\varphi H\to \hat x\boxtimes_\psi H$ an isomorphism of right $H$-torsors. As in the discussion above, apply the functor $x\boxtimes:A^H\to A^H$ to get an isomorphism $x\boxtimes q:x\boxtimes_\varphi H\to x\boxtimes_\psi H$ in $A^H$. 
    
	We have to describe a $K$-torsor $\hat y=S(\hat x)$ in $A$, so in particular a functor $\hat y:K^\op\to A$. For $h\in K_0$ (i.e. $h:\varphi(*)\to\psi(*)$ is an arrow in $H$ for $G_0=\{*\}$; see \cref{eg.equalizers}), $\hat y(h)$ is defined by the following pullback diagram in $A$:
    \begin{equation}\label{eq.xth}
        \tikz[baseline=(current  bounding  box.center)]{
        \path (0,0) node (1) {$\hat y(h)$} +(6,0) node (2) {$x(*)$} +(6,-1.5) node (3) {$x(*)\cdot \{\id_{\varphi(*)}\}$} +(6,-3) node (4) {$x(*)\cdot H(\varphi(*),\varphi(*))$} +(6,-4.5) node (5) {$x(*)\cdot H(\psi(*),\varphi(*))$} +(0,-4.5) node (6) {$x(*)\cdot \{h\}$};    
        \draw[->] (1) -- (2);
        \draw[->] (2) -- (3) node[pos=.5,auto] {$\scriptstyle \cong$};
        \draw[->] (3) -- (4);
        \draw[->] (4) -- (5) node[pos=.5,auto] {$x\boxtimes q$};
        \draw[->] (6) -- (5);
        \draw[->] (1) -- (6);   
      } 
  \end{equation}
The arrow labelled $x\boxtimes q$ in this diagram isn't, strictly speaking, the $x\boxtimes q$ we had before, but rather its restriction above $\varphi(*)\in H_0$. 

	If we want the $\hat y(h)$ to make up the object function part of a functor $K^\op\to A$, we have to describe how an arrow $g\in G$ induces a map from $\hat y(h)$ to $\hat y(\psi(g^{-1})h\varphi(g))$ (see \cref{eg.equalizers} for a description of the groupoid $K$). This is somewhat notationally cumbersome, but straightforward enough. One simply has to notice that $x\boxtimes q:x\boxtimes\hat x\boxtimes_\varphi H\to x\boxtimes \hat x\boxtimes_\psi H$ commutes with the left $G$-action on $x$, and then run through the identifications made in diagram \cref{eq.Hx}. We leave the details to the reader.  
    
	Now notice that the canonical maps $\hat y(h)\to x(*)$ (the horizontal maps in \cref{eq.xth}) actually identify $x(*)$ with the coproduct $\displaystyle\coprod_h \hat y(h)$ (with $h$ running through the arrows $\varphi(*)\to\psi(*)$ in $H$). This follows from the condition that coproducts are stable (\cref{defn.goodness_cond}): as $h$ runs through the set indicated above, the bottom horizontal maps in \cref{eq.xth} are the structure maps of a coproduct, and they are pulled back through the right-hand vertical map. One can thus think of this construction of a right $K$-representation in $A$ as a ``breaking up'' of $\hat x(*)$ into pieces $\hat y(h)$, and the right $G$-action sketched in the previous paragraph is just the right action of $G$ on $\hat x$, restricted to these individual pieces. Rigorously, what we mean is that for an arrow $g\in G$, when identifying $\hat x(*)$ with $\displaystyle\coprod_h \hat y(h)$, we have an equality of arrows:
	\begin{equation}\label{eq.breakup}
  	\hat x(*)\verylongarrow{g} \hat x(*)\qquad=\qquad\coprod_h \hat y(h) \verylongarrow{\coprod_h g} \coprod_h \hat y(\psi(g^{-1})h\varphi(g)).
	\end{equation} 
The coproducts run through the arrows $h:\varphi(*)\to\psi(*)$ in $H$. 
    
	This view of $\hat y$ as a decomposition of $\hat x$ into pieces which are permuted by the maps $g\in G$ makes it clear why one of the two requirements for being a torsor if satisfied (\cref{defn.torsor}), namely the one asking that the colimit of the functor $\hat y:K^\op\to A$ be the monoidal unit $1_A$. Indeed, the colimits of $\hat x$ and $\hat y$ are supposed to represent the functors $\cat{Cocones}(\hat x,-)$ and $\cat{Cocones}(\hat y,-)$ from $A$ to $\cat{Set}$, respectively. Since $\displaystyle \hat x(*)=\coprod_h\hat y(h)$, giving an arrow out of $\hat x(*)$ is the same as giving one arrow out of each $\hat y(h)$. Finally, by the identifications in \cref{eq.breakup}, giving an arrow $\hat x(*)\to a\in A$ commuting with the $G$-action (i.e.\ a cocone $\hat x\to a$) is the same as giving arrows out of the $\hat y(h)$ commuting with the right $G$-action on $\hat y$. In other words, we have a natural isomorphism $\cat{Cocones}(\hat x,-)\cong\cat{Cocones}(\hat y,-)$. Since the former is represented by $\varinjlim\hat x=1_A$ (because $\hat x$ is a $G$-torsor), so must the latter.    
    
	The final task in the construction of $\hat y$ is checking that $\tau_y:y\boxtimes_{K_0} K\to y\boxtimes_{K_0}^{\pi_0}\hat y$ is an isomorphism in $^{K_0}A^K$ (where $\pi_0=\pi_0(K)$; cf. \cref{subsection.coends}). Because $K$ is a groupoid, this is equivalent to proving that $\tau_y$ is an isomorphism when regarded as a map in $^{K_0}A^{K_0}$. Now fix $h,h'\in K_0$. We have components $\tau_y^{h,h'}:\coprod y(h)\to y(h)\times y(h')$ of $\tau_y:$, where the coproduct is indexed by the set of arrows $g\in G$ making the diagram
 \[
     \tikz[baseline=(current  bounding  box.center)]{
       \path (0,0) node (1) {$\varphi(*)$} +(2,0) node (2) {$\psi(*)$} +(2,-2) node (3) {$\psi(*)$} +(0,-2) node (4) {$\varphi(*)$};    
       \draw[->] (1) -- (2) node[pos=.5,auto] {$\scriptstyle h$};
       \draw[->] (4) -- (3) node[pos=.5,auto] {$\scriptstyle h'$};
       \draw[->] (4) -- (1) node[pos=.5,auto] {$\scriptstyle \varphi(g)$};
       \draw[->] (3) -- (2) node[pos=.5,auto,swap] {$\scriptstyle \psi(g)$};   
     }        
 \]
commutative, and we need to prove that $\tau^{h,h'}$ is an isomorphism whenever $(h,h')\in\pi_0$. Now, upon identifying $x(*)$ with $\coprod y(h)$ as discussed above, the coproduct \[\coprod_{h,h'}\tau_y^{h,h'}:x(*)\cdot G\to x(*)\times x(*)\] is nothing but $\tau_x$. Since this is an isomorphism (because $\hat x$ is a torsor), the conclusion follows from \cref{lemma.usegoodness} below.       
\end{defnS}    

\begin{remark}\label{remark.y_is_equalizer}
    Note that $x\boxtimes q:x\boxtimes\hat x\boxtimes_\varphi H\to x\boxtimes\hat x\boxtimes_\psi H$ is a map over $x$, in the sense that for all $s\in H_0$, the triangle
    \[
        \tikz{
            \path (0,0) node (1) {$x(*)\times(\hat x\boxtimes_\varphi H)(s)$} +(5,0) node (2) {$x(*)\times(\hat x\boxtimes_\psi H)(s)$} +(2.5,-1) node (3) {$x(*)$}; 
            \draw[->] (1) -- (2) node[pos=.5,auto] {$\scriptstyle x\boxtimes q$};
            \draw[->] (2) -- (3);
            \draw[->] (1) -- (3);
        }
    \]
commutes. After identifying $x\boxtimes \hat x\boxtimes_\varphi H\cong x\boxtimes_\varphi H$ (and analogously for $\psi$), this implies that \cref{eq.xth} is actually a diagram in the comma category $A\downarrow x(*)$. Therefore the two maps out of $\hat y(h)$ actually coincide (modulo the identification $x(*)\cdot\{h\}\cong x(*)$).    
\end{remark}

\begin{lemma}\label{lemma.usegoodness}
    Let $A$ be a good cartesian 2-ring, and $f_i:a_i\to b_i$, $i\in I$ maps in $A$. Let $f:a\to b$ be the coproduct of all the $f_i$'s. If $f$ is an isomorphism then so are the $f_i$'s.  
\end{lemma}    
\begin{proof}
    Fix $i\in I$, and consider the commutative diagram
    \begin{equation}\label{eq.usegoodness}
     \tikz[baseline=(current  bounding  box.center)]{
       \path (0,0) node (1) {$a_i$} +(2,0) node (2) {$b_i$} +(2,-2) node (3) {$b$} +(0,-2) node (4) {$a$};    
       \draw[->] (1) -- (2) node[pos=.5,auto] {$\scriptstyle f_i$};
       \draw[->] (2) -- (3);
       \draw[->] (4) -- (3) node[pos=.5,auto] {$\scriptstyle f$};
       \draw[->] (1) -- (4);   
     }        
 \end{equation}
Since the bottom arrow is an isomorphism, it is enough to prove that this is a pullback square. So let $c\in A$ be an object, and let 
    \[
     \tikz[baseline=(current  bounding  box.center)]{
       \path (0,0) node (1) {$a_i$} +(2,0) node (2) {$b_i$} +(2,-2) node (3) {$b$} +(0,-2) node (4) {$a$} +(-1,1) node (5) {$c$};    
       \draw[->] (1) -- (2) node[pos=.5,auto] {$\scriptstyle f_i$};
       \draw[->] (2) -- (3);
       \draw[->] (4) -- (3) node[pos=.5,auto] {$\scriptstyle f$};
       \draw[->] (1) -- (4);
       \draw[->] (5) .. controls +(1,0) and +(-.5,1) .. (2) node[pos=.5,auto] {$\scriptstyle f'$};
       \draw[->] (5) .. controls +(0,-1) and +(-1,.5) .. (4);          
     }        
 \]
commute. Then, by coproduct stability, we have $\displaystyle c\cong\coprod_I c\times_a a_i\cong\coprod_I c\times_b a_i$. For $j\ne i$, we have a map $f'\times_b f_j: c\times_b a_j\longrightarrow b_i\times_b b_j=0_A$ (by the property that coproducts are disjoint), so $c\times_b a_j=0_A$ (simple exercise: if an object in a good cartesian 2-ring admits a map into the initial object, then it is initial). But this means that $c\cong c\times_b a_i\cong c\times_a a_i$, so the map $c\to a$ factors through $a_i$, as in the left hand triangle below:
  \[
     \tikz[baseline=(current  bounding  box.center)]{
       \path (0,0) node (1) {$a_i$} +(2,0) node (2) {$b_i$} +(2,-2) node (3) {$b$} +(0,-2) node (4) {$a$} +(-1,1) node (5) {$c$};    
       \draw[->] (1) -- (2) node[pos=.5,auto] {$\scriptstyle f_i$};
       \draw[->] (2) -- (3);
       \draw[->] (4) -- (3) node[pos=.5,auto] {$\scriptstyle f$};
       \draw[->] (1) -- (4);
       \draw[->] (5) .. controls +(1,0) and +(-.5,1) .. (2) node[pos=.5,auto] {$\scriptstyle f'$};
       \draw[->] (5) .. controls +(0,-1) and +(-1,.5) .. (4);   
       \draw[->] (5) -- (1);        
     }        
 \]
The map $c\to a_i$ is unique because $a_i\to a$ is a monomorphism, and the upper triangle commutes because $b_i\to b$ is similarly mono.  
\end{proof}
\begin{remark}\label{remark.usegoodness}
    Below we will make use of the observation that in fact, \cref{eq.usegoodness} is a pullback even without the assumption that $f$ is an isomorphism. The argument given above can easily be adapted to prove this strengthened version of the lemma. 
\end{remark}

\begin{TS=1}
    We keep the notations we've been using: $\hat x$ is a right $G$-torsor in $A$, $q:\hat x\boxtimes_\varphi H\to \hat x\boxtimes_\psi H$ is an isomorphism of right $H$-torsors, and $\hat y=S(\hat x)$. Recall that the underlying $G$-torsor of $T(\hat y)$ is simply $\hat y\boxtimes_\xi G$. We have, by definition, 
    \[
    	(\hat y\boxtimes_\xi G)(*)=\int^{h\in K_0}\hat y(h)\cdot G.
    \]
This is nothing but the coproduct of all $\hat y(h)$ for arrows $h:\varphi(*)\to\psi(*)$ in $H$, which, as seen above in the construction of $S$, can be identified with $\hat x(*)$. Furthermore, in that identification, the right action of $g\in G$ on $\hat x(*)$ is the coproduct of its actions on the ``summands'' $\hat y(h)$. We can thus identify $\hat y\boxtimes_\xi G$ with $\hat x$ as right $G$-torsors. 

	Now recall the isomorphism $p_{\hat y}:\hat y\boxtimes_{\varphi\circ\xi}H\to \hat y\boxtimes_{\varphi\circ\psi}H$. We have to prove that if the vertical arrows are the identifications \[\hat y\boxtimes_{\varphi\circ\xi} H\cong (\hat y\boxtimes_\xi G)\boxtimes_\varphi H\cong \hat x\boxtimes_\varphi H\] (and similarly with $\psi$ instead of $\varphi$) coming from $\hat y\boxtimes_\xi G\cong \hat x$, then
    \[
         \tikz[baseline=(current  bounding  box.center)]{
           \path (0,0) node (1) {$\hat y\boxtimes_{\varphi\circ\xi}H$} +(4,0) node (2) {$\hat y\boxtimes_{\psi\circ\xi}H$} +(4,-1) node (3) {$\hat x\boxtimes_\psi H$} +(0,-1) node (4) {$\hat x\boxtimes_\varphi H$};    
           \draw[->] (1) -- (2) node[pos=.5,auto] {$\scriptstyle p_{\hat y}$};
           \draw[->] (2) -- (3) node[pos=.5,auto] {$\scriptstyle \cong$};
           \draw[->] (4) -- (3) node[pos=.5,auto] {$\scriptstyle q$};
           \draw[->] (1) -- (4) node[pos=.5,auto,swap] {$\scriptstyle \cong$};
         }        
    \]
commutes. Since all maps are in $A^H$, it is sufficient to fix an object $h\in K_0$, and prove the commutativity of the outer square in
    \[
         \tikz[baseline=(current  bounding  box.center)]{
           \path (0,0) node (1) {$(\hat y\boxtimes_{\varphi\circ\xi}H)(\varphi(*))$} +(4,0) node (2) {$(\hat y\boxtimes_{\psi\circ\xi}H)(\varphi(*))$} +(4,-1.5) node (3) {$(\hat x\boxtimes_\psi H)(\varphi(*))$} +(0,-1.5) node (4) {$(\hat x\boxtimes_\varphi H)(\varphi(*))$} +(0,1.5) node (5) {$\hat y(h)\cdot\{\id_{\varphi(*)}\}$} +(4,1.5) node (6) {$\hat y(h)\cdot\{h\}$};    
           \draw[->] (1) -- (2) node[pos=.5,auto] {$\scriptstyle p_{\hat y}$};
           \draw[->] (2) -- (3) node[pos=.5,auto] {$\scriptstyle \cong$};
           \draw[->] (4) -- (3) node[pos=.5,auto] {$\scriptstyle q$};
           \draw[->] (1) -- (4) node[pos=.5,auto,swap] {$\scriptstyle \cong$};
           \draw[->] (5) -- (6) node[pos=.5,auto] {$\scriptstyle \cong$};
           \draw[->] (5) -- (1);
           \draw[->] (6) -- (2); 
         }        
    \]
The top arrow is just the identity on $\hat y(h)$. After eliminating it, the outer square becomes the boundary of the diagram     
    \begin{equation}\label{eq.heptagon}
         \tikz[baseline=(current  bounding  box.center)]{
           \path (0,0) node (1) {$\hat y(h)$} +(2,-1.5) node (2) {$\hat x(*)\cdot \{h\}$} +(2,-3) node (3) {$(x\boxtimes_\psi H)(\varphi(*))$} +(2,-4.5) node (4) {$(\hat x\boxtimes_\psi H)(\varphi(*))$} +(-2,-4.5) node (5) {$(\hat x\boxtimes_\varphi H)(\varphi(*))$} +(-2,-3) node (6) {$(x\boxtimes_\varphi H)(\varphi(*))$} +(-2,-1.5) node (7) {$\hat x(*)\cdot \{\id_{\varphi(*)}\}$};    
           \draw[->] (1) -- (2);
           \draw[->] (1) -- (7);
           \draw[->] (2) -- (3);
           \draw[->] (3) -- (4);
           \draw[->] (7) -- (6);
           \draw[->] (6) -- (5);
           \draw[->] (6) -- (3) node[pos=.5,auto] {$\scriptstyle x\boxtimes q$};
           \draw[->] (5) -- (4) node[pos=.5,auto] {$\scriptstyle q$};           
         }        
    \end{equation}     
Here, the vertical map $(x\boxtimes_\varphi H)(\varphi(*))\to(\hat x\boxtimes_\varphi H)(\varphi(*))$ is the obvious one,
    \[
        \hat x(*)\cdot H(\varphi(*),\varphi(*))\longrightarrow\int^{G}\hat x(*)\cdot H(\varphi(*),\varphi(*)),
    \]
and similarly for $\psi$; the diagonal maps are both the usual map $\hat y(h)\to\hat x(*)$, once we've identified the copower of $\hat x(t)$ by a singleton with $\hat x(t)$. The bottom square in \cref{eq.heptagon} is clearly commutative, while the top pentagon is essentially \cref{eq.xth} (see \cref{remark.y_is_equalizer}). It follows that \cref{eq.heptagon} is commutative, and we are done.   
       
	This argument shows that one can identify $T(\hat y)\in\cat{Tors}(G,A;\varphi,\psi)$ with $\hat x$ together with the isomorphism $q:\hat x\boxtimes_\varphi H\to\hat x\boxtimes_\psi H$ we started out with. Finally, checking the functoriality of this identification is a simple diagram chase.   
\end{TS=1}

\begin{ST=1}
	Let $\hat y$ be a right $K$-torsor in $A$, and $\hat x=\hat y\boxtimes_\xi G$ the underlying right $G$-torsor of $T(\hat y)$. As in the previous proof, we then have $\displaystyle \hat x(*)=\coprod_h\hat y(h)$, where $h$ ranges through the arrows $\varphi(*)\to\psi(*)$ in $h$. Furthermore, the right action of $g\in G$ on $x(*)$ is just the coproduct of the right actions of the same $g$ on the summands $\hat y(h)$ of $\hat x(*)$.  
    
	It is now almost tautological that $\hat y(h)$ fits in a diagram \cref{eq.xth}, with $q=p_{\hat y}:\hat y\boxtimes_{\varphi\circ\xi}H\to\hat y\boxtimes_{\psi\circ\xi}H$. But by \cref{remark.usegoodness}, this means that $y(h)$ is uniquely determined by $(\hat x,q)\in\cat{Tors}(G,A;\varphi,\psi)$. Now let $g\in G$, regarded also as an arrow $h\to h'$ in $K$. Since the right action of $g$ on $y(h)$ fits in the square
    \[
        \tikz{
            \path (0,0) node (1) {$\hat y(h)$} +(2,0) node (2) {$\hat y(h')$} +(2,-2) node (3) {$\hat x(*)$} +(0,-2) node (4) {$\hat x(*)$};
            \draw[->] (1) -- (2) node[pos=.5,auto] {$\scriptstyle g$};
            \draw[->] (4) -- (3) node[pos=.5,auto,swap] {$\scriptstyle g$};
            \draw[->] (2) -- (3);
            \draw[->] (1) -- (4);
        }
    \]
and the maps $\hat y(h)\to\hat x(*)$ are mono, the right action of $g$ on $\hat y(h)$ is also uniquely determined by the object $(\hat x,q)\in\cat{Tors}(G,A;\varphi,\psi)$. 

    In summary, for every right $K$-torsor $\hat y$ in $A$, we have an isomorphism $\hat y\cong ST(\hat y)$.
\end{ST=1}
\end{proof of equalizers bis}
   
	We have now proven \cref{prop.equalizers}. Together with \cref{prop.products} and the observation that $\cat{Gpd}_f\to\cat{Gpd}$ preserves finite limits, it implies (cf.\ \cref{defn.finitely}):
\begin{proposition}\label{prop.limits}
  The functor $\sch\Box: \cat{Gpd} \to \cat{Good2AfSch}$ and its restriction to $\cat{Gpd}_f$ are finitely continuous.\qedhere
\end{proposition}

	Finally, this allows us to prove the main result of the paper:

\begin{theorem}\label{thm.Gal}
  If $A$ is a good 2-ring, the 2-functors $\pi_1(A):\cat{Gpd}\to\cat{Gpd}$ and $\pi_1^f:\cat{Gpd}_f\to\cat{Gpd}$ are pro-representable.
\end{theorem}
\begin{proof}
  We noticed above in \cref{eq.scheming} that $\pi_1(A)$ is nothing but $\Hom_{\cat{2AfSch}}(\Spec(A),-)$ composed with $\sch\Box$. Since the former is clearly finitely continuous while the latter is finitely continuous by \cref{prop.limits}, the conclusion follows from \cref{prop.prorep}. The analogous argument applies to $\pi_1^f$.  
\end{proof}

\Section{\texorpdfstring{Examples and connections with prior definitions of $\pi_1$}{Examples and connections with prior definitions of pi1}} \label{section.examplesandconnections}

	In this section, we compute $\pi_1(A)$ and/or $\pi_1^f(A)$ for some commutative 2-rings $A$, and show that the construction agrees (in the appropriate sense) with some other definitions of $\pi_1$ which have appeared in the literature. As the results are by nature somewhat eclectic, our references to this literature will be rather extensive.

	We begin with the following useful observation:

\begin{proposition}\label{prop.pi1ofcat}
  Let $H$ be a small groupoid.  Then $\pi_1({^H}\cat{Set})$ is equivalent to $H$.  More generally, let $C$ be a small category.  Then $\pi_1({^H}\cat{Set})$ is the \define{universal enveloping groupoid} of $C$, i.e.\ it is the result of applying to $C$ the left adjoint of the functor $\Forget: \cat{Gpd} \to \cat{Cat}$.
\end{proposition}
 
\begin{proof}
  Let $C$ be a small category and $G$ a small groupoid.  Then a morphism ${^{G}\cat{Set}} \to {^{C}\cat{Set}}$ of commutative 2-rings is, by \cref{prop.GSet -> A are torsors}, precisely a $G$-torsor in ${^{C}\cat{Set}}$.  But this is nothing more nor less than a functor from $C$ to the category of $G$-torsors in $\cat{Set}$.  We have shown in \cref{thm.iso,prop.smallness} that the category of $G$-torsors in any commutative 2-ring is a small groupoid; in this case, it is well known that $\cat{Tors}(G,\cat{Set})$ is equivalent to $G$.  Thus $\Hom_{\cat{Com2Ring}}({^{G}\cat{Set}} \to {^{C}\cat{Set}}) \cong \Hom_{\cat{Cat}}(C,G)$. Letting $G$ range over small groupoids, the functor $\Hom_{\cat{Cat}}(C,-)$ is by definition represented by the universal enveloping groupoid of $C$
\end{proof}

\begin{corollary}\label{lemma.scheming_fully_faithful}
  The functor $\Spec({^{\Box}\cat{Set}})$ is a fully faithful embedding $\cat{Gpd} \mono \cat{Af2Sch}$.
\end{corollary}

\begin{remark}
  One can equivalently study torsors in terms of right-principal bibundles, and these results follow from \cite{Blohmann2008} after dropping all references to a differentiable structure.
\end{remark}

\begin{remark}
  The proof of \cref{prop.pi1ofcat} depends on the fact that $G$ is a groupoid.  In \cref{section.nongroupoids} we will discuss what happens when $G$ is allowed to range over categories.  Things will immediately fail: we will show that there are generically more 2-ring morphisms ${^{G}\cat{Set}} \to {^{C}\cat{Set}}$ than functors $C \to G$ when $G$ is not a groupoid.  Nevertheless, $\Spec({^{\Box}\cat{Set}}) : \cat{Cat} \to \cat{Af2Sch}$ is faithful on 1-morphism and fully faithful on 2-morphisms.
\end{remark}

\subsection{Galois theory for rings and fields}

We next turn to commutative 2-rings of the form $\cat{Mod}^R$ for $R$ a commutative ring. We will be interested mostly in the profinite $\pi_1^f(R)=\pi_1^f(\cat{Mod}^R)$, as this is the construction that has a counterpart in the literature to which we can compare our version.  The Galois theory of arbitrary commutative rings was first described in \cite{Magid1974}; see also \cite{Borceux2001} for a modern treatment.

	Recall that to any commutative ring $R$, one can associate a profinite (i.e.\ compact, Hausdorff, totally disconnected) topological space $\PSpec(R)$, called its \define{Pierce spectrum} (see \cite[Definition II.1]{Magid1974}, where it is called the \define{Boolean spectrum}, or \cite[Definition 4.2.4]{Borceux2001}): $\PSpec(R)$ is the set of connected components of the Zariski spectrum $\spec(R)$ of $R$, endowed with the final topology via the map $\spec(R)\to\PSpec(R)$. $\PSpec$ is a functor from the category $\cat{AfSch}$ of (1-)affine schemes to the category $\cat{PfTop}$ of profinite topological spaces. 

	Magid \cite[Definition IV.16 and Theorem IV.20]{Magid1974} associates to any commutative ring $R$ a \define{separable closure} $S$, which is an $R$-algebra uniquely determined up to isomorphism of $R$-algebras.  Consider the \define{kernel pair} groupoid $S\otimes_R S \rightrightarrows S$
in \cat{AfSch}; by \cite[Corollary IV.28]{Magid1974}, it maps under $\PSpec$ to a groupoid $\Pi(R)$ in \cat{PfTop}, which Magid calls the \define{fundamental groupoid of $R$}. Our result is that $\pi_1^f(R)=\Pi(R)$, in the following sense:

\begin{theorem}\label{prop.Magid}
    Let $R$ be a commutative ring. Regard the 2-category $\cat{Gpd}_f$ as a sub-2-category of $\cat{Gpd}(\cat{PfTop})$ (the latter with topological functors and topological natural transformations as 1- and 2-morphisms respectively). Then the 2-functor $\cat{Gpd}_f\to\cat{Gpd}$ represented by $\Pi=\Pi(R)\in\cat{Gpd}(\cat{PfTop})$ is equivalent to $\pi_1^f(\cat{Mod}^R)$. 
\end{theorem} 

We will need some preliminary results: 

\begin{lemma}\label{lemma.torsors_are_split}
    Let $R$ be a commutative ring, $G$ a finite groupoid, and $\hat x$ a right $G$-torsor in $\cat{Mod}^R$. For every $s\in G_0$, $\hat x(s)$ is projective finitely generated over $R$. Moreover, the duals $\hat x(s)^*=\Hom_R(\hat x(s),R)$ are separable algebras over $R$. 
\end{lemma}
\begin{sketch}
    By \cref{lemma.main} and \cref{prop.GSet -> A are torsors} (and their proofs), $\hat x$ is the image through a monoidal functor $^G\cat{Mod}^R\to A=\cat{Mod}^R$ of the object $G\cdot R\in {^GA^G}$ defined by setting $(G\cdot R)(t,s)$ ($s,t\in G_0$) to be the free $R$-module on the set $G(t,s)$. Since $G$ is finite, $G\cdot R$ is a dualizable object in $^GA^G$, and so  $\hat x\in A^G$ must also be dualizable (and hence so must all the $\hat x(s)\in A$ for $s\in G_0$). The first statement now follows from the well-known fact that the dualizable objects of $\cat{Mod}^R$ are exactly the projective finitely generated modules. 
    
    The second statement  when $G$ is a group is essentially \cite[Corollary 2.4]{Nuss2006} applied to the algebra map $R\to\hat x$; the extension of the argument to arbitrary finite groupoids is simple enough conceptually, but notationally cumbersome, and so we leave it to the reader. 
\end{sketch}

Now consider the cartesian monoidal category $B=\cat{AfSch}_R$ of affine $R$-schemes. Let $G$ be a finite groupoid, and $\hat x$ a right $G$-torsor in $B$ (\cref{defn.torsor} applies verbatim to any symmetric monoidal category). It can be shown that the commutative $R$-algebras $\hat x(s)$, $s\in G_0$ are projective finitely generated: when $G$ is a group, this is for example a consequence of \cite[Theorem 1.7 (1)]{Kreimer1981}, and the general case follows by regarding a finite groupoid as a finite disjoint union of finite groups. 

Since the dualization functor $\Hom_R(-,R)$ implements an equivalence between the projective finitely generated commutative $R$-algebras and projective finitely generated cocommutative $R$-coalgebras, this argument, together with the first statement in \cref{lemma.torsors_are_split} and the identification \[\cat{Tors}(G,\cat{Cog}(\cat{Mod}^R)) \simeq \cat{Tors}(G,\cat{Mod}^R) \simeq \pi_1^f(R)(G)\] (cf.\ \cref{prop.GSet -> A are torsors} and \cref{cor.cartify}), proves:

\begin{lemma}\label{lemma.torsors_in_coalg}
	With the notations used above, the functor $\Hom_R(-,R)$ implements an equivalence $\cat{Tors}(G,\cat{Alg}_R^\op) \simeq \pi_1^f(R)(G)$. Moreover, this equivalence is natural in $G\in\cat{Gpd}_f$.\qedhere 
\end{lemma}

We now recall some more terminology, enough to be able to phrase the proof. We refer the reader to \cite{Borceux2001} for the missing details. For a ring homomorphism $f:R\to S$, Borceux and Janelidze introduce the notion of $R$-algebra \define{split} by $f$ (\cite[Definition 4.5.1]{Borceux2001}).  It will not be important for us to know precisely what this means; we only need to know that when $f:R\to S$ is the separable closure of $R$, the separable projective finitely generated $R$-algebras are certainly split by $f$ (this follows from \cite[A.1, remark on page 309]{Borceux2001}). We write $\cat{Split}_R$ for the category of $R$-algebras split by the separable closure $R\to S$, and refer to the objects of $\cat{Split}_R$ as \define{split $R$-algebras}.  

	Combining (the proof of) \cite[Theorem 4.7.15]{Borceux2001} with the fact that the separable closure $R\to S$ is of Galois descent in the sense of \cite[Definition 4.5.2]{Borceux2001} (noted in passing right after that definition), we have:

\begin{proposition}[Borceux, Janelidze, Magid]\label{prop.BJM}
  Let $R$ be a commutative ring, and $\Pi=\Pi(R)$ its profinite fundamental groupoid in Magid's sense. Let $^\Pi\cat{PfTop}$ be the category of topological representations of $\Pi$ in $\cat{PfTop}$.  The Pierce spectrum functor implements an equivalence\\[-1ex]
    
    \mbox{} \hfill $\displaystyle
        \cat{Split}_R^\op\simeq{^\Pi}\cat{PfTop}.
    $
\qedhere
\end{proposition}    

We are now ready to prove our result. 

\begin{proof of Magid}
    By \cref{lemma.torsors_are_split,lemma.torsors_in_coalg} and the observation made above that separable projective finitely generated $R$-algebras are split, we have an equivalence \[\pi_1^f(R)(G) \simeq \cat{Tors}(G,\cat{Split}_R^\op).\] Composing further with the equivalence \[ \cat{Tors}(G,\cat{Split}_R^\op) \simeq \cat{Tors}(G,{^\Pi}\cat{PfTop})\] obtained from \cref{prop.BJM}, we get \[ \pi_1^f(R)(G) \simeq  \cat{Tors}(G,{^\Pi}\cat{PfTop}).\]
    
    Now notice, just as in the proof of \cref{lemma.scheming_fully_faithful}, that the objects of $\cat{Tors}(G,{^\Pi}\cat{PfTop})$ are the right-principal $\Pi$--$G$ bibundles in $\cat{PfTop}$. By the appropriate analogue of \cite[Proposition 3.7]{Blohmann2008}, we will be able to conclude that $\cat{Tors}(G,{^\Pi}\cat{PfTop}) \simeq \Hom_{\cat{Gpd}(\cat{PfTop})}(\Pi,G)$ if for every such bibundle $X$, the map $X\to \Pi_0$ admits a section. But the map $X\to \Pi_0$ is surjective by one of the torsor conditions, and by \cite[Theorem IV.16 i)]{Magid1974}, $\Pi_0$ is an extremely disconnected compact Hausdorff space (\cite[Definition I.15]{Magid1974}), and hence \cite[Theorem I.19]{Magid1974} surjections onto it in $\cat{PfTop}$ split.
\end{proof of Magid}

	Now let us specialize to the case when $R$ is a field $k$. First, it is not hard to show that Magid's $\Pi(k)$ is just the usual absolute Galois group $\Gal(k)$ (by applying \cite[Corollary 4.7.16]{Borceux2001} to a separable closure $k\to k_s$ of $k$, for example). In particular, $\pi_1^f(k)=\Gal(k)$ in the sense of \cref{prop.Magid}.  But in the field case, we can drop all finiteness conditions.  This result is probably well known, but we have not found a direct reference: 

\begin{proposition}\label{prop.pi1_of_fields}
    Let $k$ be a field. Then, $\pi_1(k)=\pi_1(\cat{Vect}^k)$ is represented by $\Gal(k)$. 
\end{proposition} 
\begin{proof}
    Let $G$ be a small groupoid. Because $k$ is a connected ring (i.e.\ it has no non-trivial idempotents), a right $G$-torsor in $\cat{Vect}^k$ is actually a right $G_i$-torsor for some connected component $G_i$ of $G$. Similarly, $\Hom(\Gal(k),-)$ commutes with disjoint unions of groupoids because $\Gal(k)$ has one object. So it suffices to focus on $G$-torsors over $k$ for \define{groups} $G$. In this case, a right $G$-torsor in $\cat{Vect}^k$ is nothing more than a cocommutative right $G$-module coalgebra $C$, satisfying the conditions that (a) the map $C\otimes k[G]\to C\otimes C$ defined by $c\otimes g\mapsto c_{(1)}\otimes c_{(2)}g$ is an isomorphism, and (b) the counit $\varepsilon:C\to k$ is the quotient of $C$ by $G$ (cf.\ \cref{defn.torsor}). Using (a), one sees easily that $C$ is trivial (i.e.\ isomorphic to $k[G]$ as right $G$-module coalgebra) if and only if it has a grouplike element. 
    
    For a finite Galois extension $k\subseteq K$, let $\cat{Tors}(G,k\backslash K)$ be the category of right $G$-torsors over $k$ which become trivial after extending scalars to $K$, with the usual torsor maps as morphisms. Let $\cat{Triv}(G,K/k)$ be the category whose objects are left $\Gal(K/k)$-module structures on $K[G]$ making both the counit $K[G]\to K$ and the comultiplication $K[G]\to K[G]\otimes_KK[G]$ $\Gal(K/k)$-equivariant, and commuting with the right $G$-action, and whose morphisms are the $\Gal(K/k)$-equivariant torsor maps $K[G]\to K[G]$. By Galois descent, tensoring with $K$ gives an equivalence:
    \begin{equation}\label{eq.fields1}
        \cat{Tors}(G,k\backslash K) \simeq \cat{Triv}(G,K/k).
    \end{equation}    
        
    Since cocommutative coalgebras over separably closed fields always have grouplike elements, it follows that extending scalars from $k$ to its separable closure $k_s$ makes any $G$-torsor $C$ trivial. But a grouplike of $C\otimes k_s$ lives in some $C\otimes K$ for $k\subseteq K$ as above, and hence we have a filtered union $\cat{Tors}(G,k)=\bigcup_K \cat{Tors}(G,k\backslash K)$ of categories. Combining this with \cref{eq.fields1} and with the fact that filtered unions of categories are equivalent to the (2-)colimits of the corresponding diagrams in $\cat{CAT}$ (\cref{prop.strictcolim}), we get:
    \begin{equation}\label{eq.fields2}
        \cat{Tors}(G,k) \simeq \varinjlim_K \cat{Triv}(G,K/k).  
    \end{equation}
    
    Now notice that the subset $G\subset K[G]$ is invariant under the action of $\Gal(K/k)$ on $K[G]$ as in the definition of $\cat{Triv}(G,K/k)$ (as $G$ is precisely the set of grouplikes in $K[G]$), and sending an object of $\cat{Triv}(G,K/k)$ to $G$ with the corresponding left $\Gal(K/k)$-action and its usual right $G$-action implements an equivalence between $\cat{Triv}(G,K/k)$ and the category of right principal $\Gal(K/k)$--$G$ bibundles. But by an argument analogous to the one used in the proof of \cref{prop.Magid} (and based on \cite[Proposition 3.7]{Blohmann2008}), the latter is equivalent to $\Hom(\Gal(K/k),G)$ (topological homomorphisms and natural transformations). Plugging this into \cref{eq.fields2} gives:
    \[
        \cat{Tors}(G,k) \simeq \varinjlim_K \Hom(\Gal(K/k),G) \simeq \Hom(\Gal(k),G).
    \]
The first equivalence is what we've just shown, while the second one follows from \cref{prop.strictcolim} again. 
\end{proof}
\begin{remark}
    We do not know if \cref{prop.pi1_of_fields} holds for arbitrary commutative rings. 
\end{remark}

\subsection{The \'etale fundamental pro-group of a connected site} \label{section.comparisonArtin}

	Our final comparison is to the \'etale fundamental pro-group $\pi_1^\et$ of a connected site in the sense of \cite{Artin1986}. Not surprisingly, that notion and ours agree in the appropriate sense. The details will again be rather sketchy. We recall only that the authors of \cite{Artin1986} introduce, for any locally connected site $C$, a pro-object in the homotopy category of simplicial sets. When the site is connected and pointed, the homotopy group functors $\pi_n$ yield, when applied to this construction, the \'etale homotopy pro-groups $\pi_n^\et(C)$ (cf.\ \cite[$\S$9]{Artin1986}).  

	The pro-object $\pi_1^\et$ of the (1-)category $\cat{Gp}$ of groups has a property very similar to the one we are after (cf.\ \cite[Corollary 10.7]{Artin1986}, and also \cite[Proposition 5.6]{Friedlander1982}): 

\begin{proposition}
	Let $C$ be a pointed connected site, and $\sh(C)$ the topos of sheaves on $C$. For any group $G$ we have a bijection of sets:
  	\begin{equation}\label{eq.AM}
    	\{\text{isomorphism classes in }\cat{Tors}(G,\sh(C))\} \cong \Hom_{\cat{Pro-Gp}}(\pi_1^\et(C),G)/G,
    \end{equation}
where we are quotienting by the conjugation action of $G$.     
\qedhere
\end{proposition}

\begin{remark}
	There seems to be an error in the statement of \cite[Corollary 10.7]{Artin1986}, where the quotienting by $G$ is omitted. Friedlander, on the other hand, states the result in terms of pointed torsors, so there is no need for quotienting in \cite[Proposition 5.6]{Friedlander1982}. We thank the referee for pointing all of this out. 
\end{remark}

	Note that the right hand side of \cref{eq.AM} is just the colimit $\pi_0\left(\varinjlim_i \limits \Hom_{\cat{Gpd}}(G_i,G)\right)$ in $\cat{Set}$, where $G_i$ are the groups appearing in a diagram representing $\pi_1^\et(C)\in\cat{Pro-Gp}$ (here regarded as groupoids). One can now refine the proof of this result. First, one can replace the bijection of sets by an equivalence
	\begin{equation}\label{eq.AMbis}
  	\cat{Tors}(G,\sh(C))\simeq \varinjlim' \nolimits \Hom_{\cat{Gpd}}(G_i,G)
  \end{equation}
of categories. Here, the right hand side is the strict colimit of categories, as in \cref{prop.strictcolim} and the discussion preceding it. One can then use the connectedness of $C$ (this is simply the condition that the final object of $\sh(C)$ cannot be written as a non-trivial binary coproduct; cf.\ \cite[$\S$9]{Artin1986}) to extend the result from groups $G$ to arbitrary groupoids, simply by noting that both sides in \cref{eq.AMbis} distribute over disjoint unions of groups. Finally, applying \cref{prop.strictcolim} to replace $\varinjlim'$ by $\varinjlim$, we find:

\begin{proposition}\label{prop.pi1=etpi1}
    For any connected site $C$, $\pi_1(\sh(C))$ and the \'etale fundamental pro-group $\pi_1^\et(C)$ of \cite{Artin1986} are represented by the same left-filtered diagram in $\cat{Gpd}$. 
\qedhere  
\end{proposition}


\Section{Generalizations and directions for further research}\label{section.furtherresearch}

In this paper we have proposed a categorification of affine algebraic geometry; our categorification is rich enough to allow for an ``affine'' language for arbitrary locales, schemes, etc.  For every affine 2-scheme $X$, we defined a functor $\pi_1(X): \cat{Gpd} \to \cat{Gpd}$ which controls the torsors over $X$; we also gave conditions (satisfied in all examples) for this functor to be pro-representable, and we declared the representing pro-object the \define{fundamental pro-groupoid of $X$}.  In this final section, we will briefly outline three possible generalizations to be pursued in future work: first, one can try to enrich $\pi_1(X)$ to a pro-groupoid in a category richer than \cat{Set}; second, one can try to control ``torsors'' for non-groupoids; and third, one can try to apply our machinery towards even higher categorifications.


\subsection{Fundamental pro-ADJECTIVE groupoids}

Given an affine 2-scheme $X$, we defined $\pi_1(X)$ as the representing object of $G \mapsto \Hom_{\cat{Af2Sch}}(X,G)$, where $G$ ranges over groupoids in \cat{Set}, and where we think of it as an affine 2-scheme via $G = \sch G$, c.f.\ \cref{lemma.scheming_fully_faithful}.  This approach is justified because, for $G$ a groupoid in \cat{Set} and $X$ an arbitrary affine 2-scheme, we showed in \cref{thm.iso,prop.smallness} that $\Hom_{\cat{Af2Sch}}(X,G)$ is an essentially small groupoid, and so it stands a chance of being a groupoid of the form $\Hom_{\cat{ProGpd}}(\pi_1(X),G)$; we also gave a ``geometric'' description of the objects in $\Hom_{\cat{Af2Sch}}(X,G)$ as ``$G$-torsors over $X$'' (\cref{prop.GSet -> A are torsors}).  But in (1-)algebraic geometry, we have good notions of ``torsor'' for groupoids in categories richer than \cat{Set}.  To what extent can we extend that theory here?  Is there an enriched $\pi_1$?

\begin{example}\label{eg.Hopfs}
  Two well-known sources of symmetric monoidal categories are commutative Hopf algebras (take the category of comodules) and cocommutative Hopf algebras (category of modules).  Commutative Hopf algebras are precisely \define{(1-)affine algebraic groups}. Cocommutative Hopf algebras deserve to be called \define{1-affine coalgebraic groups}: a cocommutative Hopf algebra $H$ over a commutative ring $R$ is precisely a group object in $\cat{Cog}(\cat{Mod}^R)$.
  
  If one includes infinite-dimensional representations, then, working over any commutative ring $R$, these categories are examples of commutative 2-rings, and indeed commutative $\cat{Mod}^R$-2-algebras.  The classical Tannakian theory of \cite{MR654325,MR1106898} says that one can recover an affine algebraic group $\spec(H)$, where $H$ is a commutative Hopf algebra over $R$, from the groupoids of commutative $\cat{Mod}^R$-2-algebra homomorphisms $^H\cat{Comod}^R \to \cat{Mod}^S$, as $S$ ranges over $R$-algebras.  Similar results hold for more general (co)algebraic groupoids.
\end{example}

	There are myriad variations on the following question.  For definiteness, we focus on algebraic groupoids, and ask:

\begin{example}\label{eg.question}
  Fix a commutative (1-)ring $R$, or even a field.  Let $A$ be a commutative $\cat{Mod}^R$-2-algebra.  Let $G$ range over, say, groupoid objects in the category of schemes over $\spec R$, and write $\QCoh(G)$ for the commutative $\cat{Mod}^R$-2-algebra of quasicoherent sheaves (of $R$-modules) on the corresponding stack.  Is the functor
  \[ G \mapsto \Hom_{\cat{Mod}^R\cat{-com2alg}}(\QCoh(G), A)\]
  represented by a pro-object in the category of groupoids over $\spec R$?
  
  Of course, by \cref{prop.lurieTK} one would expect that when $A$ is itself $\QCoh($a geometric stack over $R)$, then $\pi_{1}(A)$ does exist and is precisely this stack.
\end{example}

	We do not know how to answer the question in \cref{eg.question}.  The main problem is that in the $R$-linear setting, we do not know of any simplification akin to the passage $A \leadsto \cat{Cog}(A)$ that we introduced in \cref{section.cartesian}; this simplification was one of our main tools for the results in \cref{section.mainresults}.

	If a commutative 2-ring $A$ is not a 2-algebra for $\cat{Mod}^R$, then one should not expect it to have a ``$\pi_1$ enriched in groupoids over $R$.''  Nevertheless, there may be classes of affine 2-schemes that are broader than \cat{Gpd} for which any good affine 2-scheme has a ``$\pi_1$ enriched in this class.''  A reasonable minimal condition on such classes is that for any test object $G$, the functor $\Hom_{\cat{Af2Sch}}(-,G)$ should be valued in essentially-small groupoids (by definition, it is valued in large categories); as we will illustrate in \cref{section.nongroupoids}, relaxing this condition pretty much destroys the possibility of \mbox{(pro-)}representability of $\pi_1$-style functors.  We include one example, which we will revisit in \cref{section.higherschemes}, to illustrate that there are very natural and important affine 2-schemes $G$ which are not groupoids in \cat{Set} but for which $\Hom_{\cat{Af2Sch}}(-,G)$ is well-behaved as a functor $\cat{Af2Sch} \to \cat{Gpd}$:

\begin{example} \label{eg.GM}
  In \cref{eg.Set[X]} we defined a commutative 2-ring $\cat{Set}[X]$ whose underlying presentable category is the category of sheaves on a countable discrete set, by equipping the countable discrete set with the monoid structure $\NN$ and its category of sheaves with the convolution product.  If we instead equip the countable discrete set with the group structure $\ZZ$, the same construction builds a commutative 2-ring $\cat{Set}[X,X^{-1}]$.

  Let $A$ be a commutative 2-ring.  Then a morphism $\Spec(A) \to \Spec(\cat{Set}[X,X^{-1}])$ of affine 2-schemes is precisely an invertible object in $(A,\otimes)$.  Moreover, a morphism in the category $\Hom_{\cat{Af2Sch}}(\Spec(A),\Spec(\cat{Set}[X,X^{-1}]))$ is an isomorphism of invertible $A$-objects.  Along with the presentability of $A$, it follows that $\Hom_{\cat{Af2Sch}}(-,\Spec(\cat{Set}[X,X^{-1}]))$ takes values in $\cat{Gpd}$.  
\end{example}

\begin{example}
  A theory of enriched fundamental groupoids does exist in the homotopy theory of topoi.  The philosophy of topoi is that they are well enough behaved to repeat most constructions (like our move to coalgebras) that work over $\cat{Set}$.  So it is not too surprising that if a topos $T$ lives over a topos $S$, then $T$ has a homotopy type valued in $S$.
\end{example}


\subsection{Non-groupoids} \label{section.nongroupoids}

	Let $C$ be a small category.  Then $^C\cat{Set} = \Fun(C, \cat{Set})$ is a commutative 2-ring satisfying $^C\cat{Set}\boxtimes A = {^CA}$ for each commutative 2-ring $A$.  One could ask to understand the category $\Hom_{\cat{Com2Ring}}({^C\cat{Set}},A)$ in terms of $A,C$.  The hope is that in this way one could define for each commutative 2-ring $A$ a ``fundamental (pro-)category'' that controls $\Hom_{\cat{Com2Ring}}({^\Box\cat{Set}},A)$, which would be more sensitive than the groupoid-valued $\pi_1(-)$ defined in this paper.  We will argue in this section that this hope is impossible, but we do not know an elementary description of $\Hom_{\cat{Com2Ring}}({^C\cat{Set}},A)$ for any examples --- even the functor $\Hom_{\cat{Com2Ring}}({^\Box\cat{Set}},\cat{Set})$ ranging over small categories is complicated (it is the identity when applied to small groupoids).

\begin{definition}
  The \define{walking idempotent} is the category $M$ with one object $\bullet$ and a non-identity morphism $m: \bullet \to \bullet$ satisfying $m^2 = m$.  The \define{walking projection} is its idempotent splitting: it is the category $P$ with two objects $\bullet,\circ$, generated by the morphisms $p: \bullet \to \circ$ and $i: \circ \to \bullet$, with the relation $\id_\circ = \{\circ \overset i \to \bullet \overset p \to \circ\}$.
\end{definition}

	The reader is invited to calculate:
\begin{lemma}
  Let $M$ denote the walking idempotent and $P$ the walking projection.  There is an equivalence of categories $P = \Hom({^M}\cat{Set}, \cat{Set})$; the canonical injection $M \mono \Hom({^M\cat{Set}}, \cat{Set})$ sends $m \mapsto \{\bullet \overset p \to \circ \overset i \to \bullet\}$. 
  
  Let $f: M \to P$ denote this map.  For any commutative 2-ring $A$, the morphism $f^*: {^PA} \to {^MA}$ is an equivalence of commutative $A$-2-algebras.
  \qedhere
\end{lemma}
Thus the 2-category \cat{Cat} of small categories does not naturally embed into \cat{Af2Sch}, in marked contrast with \cref{lemma.scheming_fully_faithful}.

Given \cref{thm.Gal}, one might still hope that there be a small (pro-)category $\pi_1(\cat{Set})$ with 
$\Hom_{\cat{Com2Ring}}({^C\cat{Set}}, \cat{Set}) \cong \Hom_{\cat{Cat}}(\pi_1(\cat{Set}) , C)$ for each small category $C$.  Our next result shows that this is impossible:
\begin{proposition}
  Let $M,P$ denote respectively the walking idempotent and walking projection.  There does not exist a pro-category $X$ with $\Hom_{\cat{ProCat}}(X, M) = P$.
\end{proposition}
\begin{sketch}
  We prove that $X$ cannot be a small category; we leave the ``pro'' case to the reader.
  
  Let $X$ be a small category.  Then there is a map $\Hom_{\cat{Cat}}(X,M) \to M$ which sends all non-isomorphisms in $\Hom_{\cat{Cat}}(X,M)$ to $m\in M$.  This map is a surjection on morphisms if $X \neq \emptyset$.  But $P$ does not surject onto $M$.
\end{sketch}

This illustrates that a ``category-valued Galois theory'' is necessarily very different from the groupoid-valued theory that has been developed.  Even accepting that $\cat{Tors}(-,A)$ is not \linebreak\mbox{(pro-)}representable, one can still ask to understand this functor in a more hands-on way.  For any category $C$ and any commutative 2-ring $A$, the category $\Hom_{\cat{Com2Ring}}({^C\cat{Set}}, A)$ has all idempotents split, but our next example shows that the functor $\pi_1(A)(-)$ does more than split idempotents:
\begin{example}
  Let $\mathbb N$ denote the category with one object, freely generated by a morphism.  The category $\Hom_{\cat{Com2Ring}}({^{\NN}\cat{Set}}, \cat{Set})$ has an object outside of the essential image of $\mathbb N \mono \Hom_{\cat{Com2Ring}}({^{\NN}\cat{Set}}, \cat{Set})$ with endomorphism monoid $\mathbb Z$.
\end{example}

Our intuition says that the objects in $\Hom_{\cat{Com2Ring}}({^C\cat{Set}}, \cat{Set})$ correspond to ``partial group\-oid\-iza\-tions of $C$.''  Studying the functor $\Hom_{\cat{Com2Ring}}(^\Box\cat{Set},\cat{Set}) = \cat{Tors}(-,\cat{Set}) : \cat{Cat} \to \cat{Cat}$ in more detail, and in particular making this intuition precise, will be the topic of future work.


\subsection{Higher affine schemes} \label{section.higherschemes}

In this final subsection we outline very briefly (and conjecturally) why we think it would be interesting to extend the discussion to a higher-categorical setting.
 
One reason is that one might want to study ``higher Picard-Brauer groups.''  We described in \cref{eg.GM} the affine 2-scheme $\Spec(\cat{Set}[X,X^{-1}])$.  Maps into this 2-scheme classify ``line bundles'' over the source 2-scheme, and hence $\Spec(\cat{Set}[X,X^{-1}])$ represents the ``Picard group.''  In fact, $\Spec(\cat{Set}[X,X^{-1}])$ is the first term in what should be an infinite sequences of higher stacks.  We will now describe the zeroth term.

Let $\ZZ$ denote the groupoid with one object, freely generated by an invertible morphism.  One can give $\sh(\ZZ)$ a non-cartesian commutative 2-ring structure  by choosing $\otimes = \times_{\ZZ}$; i.e.\ if $X,Y \in \sh(\ZZ)$, then $X \otimes Y = (X\times Y) / \langle (zx,y) = (x,zy)\; \forall x\in X,y\in Y, z\in \ZZ\rangle$.  Then $\Spec(\sh(\ZZ),\otimes = \times_{\ZZ})$ corresponds in algebraic geometry to the affine group scheme $\GL(1)$.  Indeed: if $R$ is a commutative ring, then a morphism $\Spec({^{R}\cat{Mod}}) \to \Spec(\sh(\ZZ),\otimes = \times_{\ZZ})$ is precisely an invertible element of $R$.

Unfortunately, the second term in the desired infinite sequence already does not make sense as an affine 2-scheme.  Any commutative 2-ring $A$ has a symmetric monoidal 2-category $A\cat{-2Mod}$ of modules, and the invertible $A$-2-modules are the objects of the \define{2-Picard group of $A$}.  When $A = \cat{Vect}^k$ for $k$ a field, (the group of equivalence classes of objects in) this ``2-Picard group'' is the Brauer group ${\rm H}^2(\Gal(k),\GL(1))$.  (When $A = \cat{Mod}^R$ for $R$ a ring, the usual definition of  Brauer group is slightly smaller than the corresponding cohomology group; it is the latter that we expect to equal the algebraically-defined 2-Picard group.)

So $\GL(1) = \Pic^{0}$ represents the \define{0-Picard group}, just as $\Pic^{1} = \Spec(\cat{Set}[X,X^{-1}])$ represents the \define{1-Picard group}.  In general, we expect there should be a notion of ``affine $n$-scheme'' in which all the higher Picard groups $\Pic^{n}$ are representable.  Moreover, we expect that the group of equivalence classes of objects in $\Hom(X,\Pic^{n})$ is precisely ${\rm H}^n(X,\GL(1))$ (for some appropriate notion of cohomology).

Another reason to be interested in higher affine schemes is the desire to define higher $\pi_n$s in the Tannakian framework. In algebraic geometry, there exists an \'etale homotopy theory much more general than just fundamental groups \cite{MR559531,Artin1986}.  From a Tannakian perspective, higher homotopy types should control the representation theory of higher-categorical groupoids.  Making this precise would be valuable but difficult.
  
The problem with trying to define higher homotopy types by considering Tannaka-Krein categories as in \cref{defn.TKcat} quickly becomes apparent: while a groupoid is naturally represented in a 1-category like $\cat{Vect}$, to represent an $n$-groupoid requires allowing the representations to take values in an $n$-category.

Following the translation $R \mapsto \cat{Mod}^R$, candidate examples of \define{3-abelian groups} are already available: any commutative 2-ring $A$ should have a 3-abelian group (and in fact a \define{commutative 3-ring}) of modules $A\cat{-2Mod}$; we implicitly invoked this in the almost-definition of $\Pic^{2}$ above.  It is thus tempting to try to define precisely what should be a ``presentable $n$-category,'' and to show that the representation theory of a commutative $n$-ring is an example.

Unfortunately, working one's way up the categorical ladder quickly becomes intractable: for example, a definition of ``4-category'' similar to the well-known definition of 2-category that we have employed throughout this paper spans 51 pages  \cite{Trimble2006}.  One would then need to set up a theory of presentable categories, sheaves on sketches, etc., similar to what we have used in this paper.

Because of the ground-breaking work by Lurie \cite{Lurie:ly,Lurie:zr,Lurie:2009vn,Lurie2009}, Toen and Vezzosi \cite{MR2061855,MR2137288,MR2394633,MR2483943}, and others, higher-categorical representation theory is not completely out of reach.  But much more work is required to develop a full higher-categorical representation theory.  Related work in this direction exists in the fusion category literature (\define{fusion categories} are a particularly finite categorification of rings); c.f.\ \cite{MR2046203,MR2677836} and references therein.

\end{document}